\documentclass[12pt]{article}

\usepackage{amsfonts,amssymb,latexsym,amsmath}

\usepackage[english]{babel}
\usepackage{graphicx,bm}
\usepackage[margin=1in]{geometry}
\usepackage{enumitem}
\usepackage{tikz}
\tikzset{
  font={\fontsize{8pt}{12}\selectfont}}
\usepackage{verbatim}
\usetikzlibrary{matrix}

\usepackage{color}

\newcommand{\red}[1]{{\color{red}\bf  #1}}

\newtheorem{thm}{Theorem}[section]
\newtheorem{prop}[thm]{Proposition}
\newtheorem{lem}[thm]{Lemma}
\newtheorem{cor}[thm]{Corollary}
\newtheorem{rmk}[thm]{Remark}

\newenvironment{proof}
{\par\addvspace{0.3cm}\noindent{\rm Proof. }}
{\nopagebreak\mbox{}\hfill $\Box$\par\addvspace{0.25cm}}

\newcommand{\T}{\mathbb{T}}
\newcommand{\Z}{\mathbb{Z}}
\newcommand{\R}{\mathbb{R}}

\newcommand{\cL}{\mathcal{L}}
\newcommand{\cS}{\mathcal{S}}
\newcommand{\cF}{\mathcal{F}}
\newcommand{\cJ}{\mathcal{J}}
\newcommand{\LH}{{\cL(H)}}

\newcommand{\tQC}{\widetilde{QC}}

\DeclareMathOperator*{\esssup}{ess\, sup}
\DeclareMathOperator*{\slim}{s-lim}

\DeclareMathOperator{\clos}{clos}

\numberwithin{equation}{section}

\begin{document}

\title{Finite Section Method for singular integrals with operator-valued PQC-coefficients and a flip}
\author{Torsten Ehrhardt\thanks{tehrhard@ucsc.edu}\\
   Department of Mathematics\\
        University of California\\
       Santa Cruz, CA 95064, USA
  \and   Zheng Zhou\thanks{zzho18@ucsc.edu}\\
 Department of Mathematics\\
        University of California\\
       Santa Cruz, CA 95064, USA}

\maketitle



\begin{abstract}
    We establish necessary and sufficient conditions for the stability of the finite section method for operators belonging to a certain
    $C^*$-algebra of operators acting on the Hilbert space $l^2_H(\mathbb{Z})$ of $H$-valued sequences where $H$ is a given Hilbert space.
    Identifying $l^2_H(\mathbb{Z})$ with the $L^2_H$-space over the unit circle, the $C^*$-algebra in question is the one which contains all singular integral operators with flip
    and piecewise quasicontinous $\mathcal{L}(H)$-valued generating functions on the unit circle. The result is a generalization of an older result where the same problem, but without the flip
    operator was considered. The stability criterion is obtained via $C^*$-algebra methods and says that a sequence of finite sections is stable
    if and only if certain operators associated with that sequence (via $^*$-homomorphisms) are invertible.
\end{abstract}


\section{Introduction}\label{s:1}

Let $X$ be a Banach space and $A_n$ be a sequence of bounded linear operators on $X$. The sequence  $(A_n)_{n=1}^\infty$ is said to be {\em stable}
if there exists an $n_0$ such that for all $n\ge n_0$ the operators $A_n$ are invertible on $X$ and if
$$
\sup\limits_{n\ge n_0} \|A_n^{-1}\|_{\mathcal{L}(X)} <\infty.
$$
The notion of stability is very fundamental in numerical analysis and for questions where asymptotic invertibility plays a role.
Our setting is that of $X=l^2_H(\mathbb{Z})$, the space of square-summable $H$-valued sequences $(x_n)_{n\in\Z}$, where $H$ is a given Hilbert space.
Note that bounded linear operators on $X$ can be considered as infinite matrices whose entries belong to $\mathcal{L}(H)$,
the $C^*$-algebra of all bounded linear operators on $H$.

We define the finite section operator $P_n$ acting on $l^2_H(\mathbb{Z})$ by
$$
P_n: (x_k)_{k\in \mathbb{Z}} \mapsto (y_k)_{k\in \mathbb{Z}},\qquad y_{k}=\begin{cases}
        x_k& \text{if}\enskip -n\le k<n \\
        0  &   \text{if}\enskip k<-n\enskip\text{or}\enskip k\ge n.
        \end{cases}
$$
Given an operator $A\in \mathcal{L}(l^2(\mathbb{Z}))$, the problem of approximately solving the equation $Ax=y$ via
the finite section method, $P_n A P_n x_n =P_n y$, leads to the question of the stability of the sequence $A_n=P_nAP_n$ (thought of as operators acting on the image of $P_n$).
While this problem is out of reach for general $A$, much work has been devoted to give necessary and sufficient stability conditions
for $A$ taken from specific classes. In this paper we will take $A$ from a $C^*$-algebra $\mathcal{S}^J(PQC_{\cL(H)})$, which we are going to describe now.

Consider the following projection operator and the flip operator on $l^2_H(\Z)$:
 \begin{align*}
        P&:\enskip (x_n)_{n\in \mathbb{Z}}\mapsto (y_n)_{n\in \mathbb{Z}}\enskip \text{with}\enskip y_n =\begin{cases}
        x_n& \text{if}\enskip n\ge 0\\
        0  &   \text{if}\enskip n<0,
        \end{cases}\\
        J&:\enskip (x_n)_{n\in \mathbb{Z}} \mapsto (x_{-n-1})_{n\in \mathbb{Z}}.
 \end{align*}
Furthermore, for an $\cL(H)$-valued measurable and essentially bounded functions $a\in L^\infty_{\cL(H)}(\mathbb{T})=:L^\infty_{\cL(H)}$ defined on the unit circle
$\mathbb{T}=\{\,z\in\mathbb{C}\,:\,|z|=1\,\}$ with $\cL(H)$-valued Fourier coefficients
$$
a_n=\frac{1}{2\pi}\int_{0}^{2\pi} a(e^{ix}) e^{-inx}\, dx,\qquad n\in \Z,
$$
we define the Laurent operators $L(a)$ acting on  $l^2_H(\Z)$ by
   \begin{equation}\label{laurentoperator-1}
        L(a):\enskip (x_n)_{n\in \mathbb{Z}}\mapsto \left(\sum_{k\in\Z} a_{n-k}x_k\right)_{n\in \mathbb{Z}}.
   \end{equation}
Let us remark that if we identify the space  $l^2_H(\Z)$ with the Lebesgue space $L^2_H(\T)$ in the usual way, then
$P$ becomes the Riesz projection, $2P-I=S_{\T}$ is the singular integral operators on $\T$, the flip operator $J$ maps a function
$f(t)$ to a function $t^{-1}f(t^{-1})$ where $t\in \T$, and the Laurent operators correspond to the operators of multiplication
by the function $a(t)$, i.e., $f(t)$ is sent to $a(t)f(t)$, $t\in\T$. Note also that the algebra generated by $I$, $P$, $J$, and Laurent operators
$L(a)$ contains the following operators of interest, namely, singular integral operators with flip
$$
L(a) P + L(b)Q + L(c)JP + L(d) J Q
$$
and, in particular,  Toeplitz-plus-Hankel operators
$$
PL(a)P + PL(b)J P + Q.
$$
The last two classes of operators have been studied regarding Fredholmness and, to some extend, also invertiblity in the case of finite dimensional $H$ (see, e.g.,  \cite{BaEh-13, BaEh-06Ap, RoSi2012,DiSi2016} and the references given there).

The class of Laurent operators $L(a)$ with arbitrary $a\in L^\infty_{\cL(H)}$ is still too large to be handled, and therefore we restrict ourselves in this paper to $\mathcal{L}(H)$-valued
{\em piecewise quasicontinuous} symbols. This class, $PQC_{\mathcal{L}(H)}$, is defined as the smallest closed subalgebra of $L^\infty_{\cL(H)}$
containing all piecewise continuous  $\mathcal{L}(H)$-valued functions on $\T$ (the class $PC_{\mathcal{L}(H)}$) and all
quasicontinuous $\mathcal{L}(H)$-valued functions on $\T$ (the class $QC_{\mathcal{L}(H)}^s$). The precise definitions  will be given below, but note that it has been observed in \cite{ehrhardt1996symbol} that there are several reasonable ways to define quasicontinuous functions in the operator-valued case. {}For the results of this paper, we have to take the smallest class $QC_{\mathcal{L}(H)}^s$ that arises from these possibilities.

Now we are able to define the class of operators  $\mathcal{S}^J(PQC_{\cL(H)})$ as the smallest closed subalgebra of $\cL(l^2_H(\Z))$ which contains the
identity operator $I$, the projection $P$, the flip $J$, and all Laurent operators $L(a)$ with $a\in PQC_{\cL(H)}$. Let us also define the (smaller) class of operators
$\mathcal{S}(PQC_{\cL(H)})$ which is generated by the same operators except the flip operator $J$.

In order to describe the content of this paper, we need one more definition. Below we will define $\cF$ as the set of all bounded sequences
$(A_n):=(A_n)_{n=1}^\infty$ of  bounded  linear operators $A_n\in \cL( \mathrm{Im}( P_n))$ and define algebraic relations and a norm on $\cF$ that make it into a $C^*$-algebra. 
Clearly, $\cF$ contains all sequences $(P_n A P_n)$ with $A\in \cL(l^2_H(\Z))$. In Section \ref{s:4} we will introduce
a certain set $\mathcal{J}\subseteq \cF$ whose definition is dispensable now.

Let  $\cF^J(PQC_{\cL(H)})$ be the smallest closed subalgebra of $\cF$ containing all elements of $\mathcal{J}$ and all sequences $(P_n A P_n)$ with
$A\in  \cS^J(PQC_{\cL(H)})$. Furthermore, let $\cF(PQC_{\cL(H)})$ be the smallest closed subalgebra of $\cF$ containing all elements of $\mathcal{J}$ and all sequences $(P_n A P_n)$ with $A\in  \cS(PQC_{\cL(H)})$.

The goal of this paper is to establish explicit necessary and sufficient conditions for the stability of any sequence
$(A_n)\in \cF^J(PQC_{\cL(H)})$.
It is a generalization of \cite{ehrhardt1996finite}, where stability conditions were established for sequences $(A_n)\in \cF(PQC_{\cL(H)})$.
Thus the difference is that we allow for the flip operator $J$ to be present. The main result will be of the form
that a  sequence $(A_n)$ is stable if and only if a certain collection of operators associated with this sequence consist of invertible operators only (see Theorem \ref{thm92} below).
Each operator in this collection can be obtained via a *-homomorphism from $\cF^J(PQC_{\cL(H)})$ into a certain algebra of operators
(see Theorem \ref{thm91} and Propositions \ref{prop35}--\ref{prop36}).

\bigskip
Let us give an outline of the paper.
In Section \ref{s:2} we introduce basic notation, and in Section \ref{sec:3} we  modify the notions of compact
operators and of strong convergence on $l^2_H(\Z)$ in order to be able to use them in the operator-valued setting.
The proofs of stability result rely on expressing stability as an invertibility problem in a certain $C^*$-algebra
and applying a lifting theorem. This is done in Section \ref{s:4}, which simplifies the stability problem in such a way that
a localization principle can be applied in Section \ref{s:5}. Until this point the line of proof is basically the same as in  \cite{ehrhardt1996finite}.
However, while in \cite{ehrhardt1996finite} the localization is done over the maximal ideal space $M(QC)$ of (scalar) quasicontinuous functions $q\in QC$, the localization here must be taken over $M(\widetilde{QC})$, the maximal ideal space of all {\em even quasicontinuous functions}
$q\in  \widetilde{QC}$, $q(e^{ix})=q(e^{-ix})$. The reason is that the flip operator $J$ only commutes with even functions.

We are thus required to take a closer look at $M(\widetilde{QC})$ and to put it in relation with $M(QC)$.
It is well-known $M(QC)$ that decomposes into fibers $M_\tau(QC)$ over $\tau\in \T$, and each of these fibers
decomposes further into three disjoint sets:
$$
M_\tau^0(QC)=M^+_\tau(QC)\cap M^-_\tau(QC),\quad 
M^+_\tau(QC)\setminus M_\tau^0(QC),\quad
M^-_\tau(QC)\setminus M_\tau^0(QC).
$$ 
Similarly,  the maximal ideal space  $M(\widetilde{QC})$ can be decomposed into fibers
$M_\tau(\widetilde{QC})$ over $\tau\in\overline{\T_+}:=\{\,\tau\in \T\,:\, \mbox{Im} (\tau)\ge 0 \;\}$.
In the paper  \cite{MevenQC} these fibers have been analysed and it turned out that each of them decomposes
into either two or three disjoint disjoint sets depending on whether $\tau=\pm1$ or $\tau\in\T_+
:=\{\,\tau\in \T\,:\, \mbox{Im} (\tau) > 0 \;\}$ (see formulas \eqref{dec-1}--\eqref{dec-3}).
We will recall the corresponding results in Section \ref{s:6}. 
In some cases short proofs are provided for sake of illustration.

The most difficult part is to identify the local algebras obtain from the localization done in Section \ref{s:5}.
Corresponding to the afore-mentioned decomposition of $M(\widetilde{QC})$ we are led to several cases
and we show that local quotient algebras are *-isomorphic to certain algebras of concrete operators.
In one case, the concrete identification is actually despensible in view of the stability problem.
This will be done in Section \ref{s:7}--\ref{s:9}. The results of Section \ref{s:6} are used therein.
After having done the identification, we summarize what we have obtained so far and state the main result in Section 
\ref{s:10}. As mentioned above, the stability criterion is established in Theorem \ref{thm92} and it involves *-homomorphism which 
are explicitly given  in Theorem \ref{thm91} as well as in Propositions \ref{prop35} and \ref{prop36}.

\bigskip

There are intersections of our results with previous work.
In the preprint \cite{RochPreprint}  a stability criterion is established for certain sequences which include the finite sections of operators from the algebra generated by $I$,
$P$, $J$, and $L(a)$ with $a\in PC$. It does not include quasicontinuous functions and does not cover the operator-valued setting (i.e., it corresponds to $\dim H<\infty$). On the other hand,
it is more general in the sense that it covers operators on $l^p$-spaces ($p\neq 2$) and the finite sections can be of a more general form
 $P_{nk} A P_{kn}$ for (fixed) positive integers $k$.
Another related work can be found in \cite{RochSantosSilbermann}. There the stability of the finite truncations of operators on $L^2(\mathbb{R})$ taken from the algebra generated by
the operator of multiplication $\chi_{[0,\infty)}$, the flip on the real line, and convolution operators with piecewise continuous generating functions on $\mathbb{R}$ is established.
Apart from the quasicontinuous and operator-valued ingredient, this means that our results are
the ``unit circle version'' instead of the ``real line version'' done in \cite{RochSantosSilbermann}.

The current paper also has applications. In the paper \cite{BaEh1} the asymptotics (as $n\to\infty$) of
determinants of certain Toeplitz + Hankel matrices $T_n(a)+H_n(b)$ with singular  symbols are determined. The proof requires as 
an auxiliary result the stability not of the underlying Toeplitz + Hankel matrices $T_n(a)+H_n(b)$, but of the finite sections 
$P_n T^{-1}(\psi)(T(c)+H(d))T^{-1}(\psi^{-1})P_n$, where the symbols $\psi, c,d$ are piecewise continuous  (see also \cite{BaEh2}).
For this application the operator-valued and quasicontinuous part of our result are not essential, and in fact, the needed result can be
obtained also from \cite{RochPreprint}.

Let us mention some further connections. In \cite{BoSiWi-1995} (see also \cite{BoSi-1994}) an operator-valued version of the
Szeg\"o-Widom limit theorem was established. In other words, the asymptotics of the determinants of Toeplitz matrices with
operator-valued entries is established. The stability results of the current paper could be useful to prove an operator-valued version for
the asymptotics of determinants of Toeplitz + Hankel matrices  with operator-valued entries.

Furthermore,  in \cite{bottcher1996infinite} the asymptotics of the finite truncations of Wiener-Hopf operators (with piecewise continuous symbols) were established, by identifying them with  Toeplitz matrices with operator-valued entries. Our results could be useful for
establishing an analogue for finite truncations of Wiener-Hopf-Hankel operators.

\newpage

\section{Preliminaries}\label{s:2}

In this section we define some notation and introduce basic results, which we need subsequently. Some of the results can be found in \cite{ehrhardt1996finite}.
Throughout the paper let $H$ stand for a (given) arbitrary Hilbert space.

\bigskip
Let $\mathbb{T} = \{z\in \mathbb{C}:|z| = 1\}$ be the unit circle, let $\mathbb{Z}$ (resp., $\mathbb{Z}^+$) stand for the set of integers (resp., non-negative integers), and set $\mathbb{Z}_n = \{-n, -n+1, \cdots, n-1\}$.  By $l_H^2(\mathbb{Z})$ we denote the Hilbert space of all two-sided sequences $x = (x_n)_{n\in \mathbb{Z}}$ with $x_n\in H$ for which
        \[   \|x\|_{l_H^2(\mathbb{Z})} := \left(\sum_{n\in \mathbb{Z}}\|x_n\|_H^2\right)^{1/2} < \infty. \]
Similarly, one can define the Hilbert spaces $l_H^2(\mathbb{Z^+})$ and $l_H^2(\mathbb{Z}_n)$. Further, for $M = \mathbb{R}$ (resp.,\ $M = \mathbb{R^+}$, or $M=[-1,1]$), denote by $L_H^2(M)$ the Hilbert space of all $H$-valued Lebesgue measurable functions $f$ on $M$ for which
        \[\|f\|_{L_H^2(M)} := \left(\int_M \|f(x)\|_H^2 dx\right)^{1/2} < \infty.\]
In the case when $H = \mathbb{C}$, we omit the index $H$.

\bigskip
For a Banach space $X$, let $\mathcal{L}(X)$ denote the space of all bounded linear operators on $X$, and by $\mathcal{K}(X)$ we refer to the ideal of all compact operators on $X$. In particular, $\mathcal{L}(H)$ is a $C^*$-algebra, and we denote its unit element by $e$.

Let $L_{\mathcal{L}(H)}^\infty$ stand for the $C^*$-algebra of all Lebesgue measurable and essentially bounded functions $a$ on $\mathbb{T}$ with values in $\mathcal{L}(H)$ and with the norm defined by
        \[  \|a\|_{L_{\mathcal{L}(H)}^\infty} := \esssup\limits_{t\in \mathbb{T}} \|a(t)\|_{\mathcal{L}(H)}. \]
Again, we omit the index $\cL(H)$ if $H = \mathbb{C}$. Note that the $C^*$-algebras $\mathcal{L}(H)$ and $L^\infty$ are $^*$-isomorphic to corresponding $^*$-subalgebras of $L_{\mathcal{L}(H)}^\infty$ consisting of constant functions or functions whose values are multiples of the unit element $e$, respectively.

The $\emph{Laurent operator}$ $L(a)$ of a function $a\in L_{\mathcal{L}(H)}^\infty$ is defined by
        \begin{equation}\label{laurentoperator}
        L(a):\enskip (x_n)_{n\in \mathbb{Z}}\mapsto \left(\sum_k a_{n-k}x_k\right)_{n\in \mathbb{Z}},
        \end{equation}
where $a_n\in \cL(H)$ is the $n$-th Fourier coefficient of $a$:
        \[ a_n = \frac{1}{2\pi} \int_0^{2\pi} a(e^{i\phi})e^{-in\phi}d\phi. \]

In \cite{page1970bounded}, it is shown that $L(a)$ is a bounded linear operator on $l_H^2(\mathbb{Z})$ whenever $a\in L^\infty_{\cL(H)}$. Conversely, suppose that $(a_n)_{n\in \mathbb{Z}}$ is a sequence with $a_n\in \mathcal{L}(H)$. Then the linear operator defined by \eqref{laurentoperator} is bounded only if there is a (uniquely determined) function $a\in L_{\mathcal{L}(H)}^\infty$ whose Fourier coefficients are $a_n$. Moreover, we have $\|L(a)\|_{\mathcal{L}(l_H^2(\mathbb{Z}))} = \|a\|_{L_{\mathcal{L}(H)}^\infty},\enskip L(ab) = L(a)L(b)$ and $L(a^*) = L(a)^*$. Hence, $L_{\mathcal{L}(H)}^\infty$ is $^*$-isomorphic to the $^*$-subalgebra of all Laurent operators in $\mathcal{L}(l_H^2(\mathbb{Z}))$. We will use this identification without citing, and for brevity we will frequently write $a$ instead of $L(a)$.

Furthermore, define the following bounded  linear operators on $l_H^2(\mathbb{Z})$:
        \begin{align*}
        P&:\enskip (x_n)_{n\in \mathbb{Z}}\mapsto (y_n)_{n\in \mathbb{Z}}\enskip \text{with}\enskip y_n =\begin{cases}
        x_n& \text{if}\enskip n\ge 0\\
        0  &   \text{if}\enskip n<0,
        \end{cases}\\
        J&:\enskip (x_n)_{n\in \mathbb{Z}} \mapsto (x_{-n-1})_{n\in \mathbb{Z}}.
        \end{align*}
We denote the identity mapping on $l_H^2(\mathbb{Z})$ by $I$, and define $Q := I-P$. Note that the following relations hold: $P^2 = P = P^*, \enskip J^2 = I, \enskip J^* = J$ and $JPJ = Q$.

When identifying $l^2_H(\Z)$ with $L^2_H(\T)$,  the \emph{singular integral operator} $S_{\mathbb{T}}$ on $\mathbb{T}$ corresponds to $P-Q = 2P-I$.
The \emph{Toeplitz operator} with generating function $a\in L_{\mathcal{L}(H)}^\infty$ is $T(a) = PL(a)P$, and the corresponding \emph{Hankel operator} is given by $H(a) = PL(a)JP$. For $a\in L_{\mathcal{L}(H)}^\infty$, let $\tilde{a}\in L_{\mathcal{L}(H)}^\infty$ denote the function
\begin{equation}\label{tilde}
   \tilde{a}(t) := a(1/t),\quad t\in \mathbb{T}.
\end{equation}
Note that  $JL(a)J = L(\tilde{a})$.


For $n\in \mathbb{Z}^+$,  introduce the following bounded linear operators acting on $l^2_H(\Z)$:
        \begin{align*}
        P_n&:\enskip (x_k)_{k\in \mathbb{Z}}\mapsto (y_k)_{k\in \mathbb{Z}}\enskip \text{with}\enskip y_k =\begin{cases}
        x_k& \text{if}\enskip -n\le k<n \\
        0  &   \text{if}\enskip k<-n\enskip\text{or}\enskip k\ge n
        \end{cases}\\
        W_n&:\enskip (x_k)_{k\in \mathbb{Z}}\mapsto (y_k)_{k\in \mathbb{Z}}\enskip \text{with}\enskip y_k =\begin{cases}
        x_{n-1-k}   & \text{if}\enskip 0\le k<n \\
        x_{-n-1-k}  & \text{if}\enskip -n\le k<0\\
        0           & k<-n\enskip\text{or}\enskip k\ge n
        \end{cases}\\
        V_{n}&:\enskip (x_k)_{k\in \mathbb{Z}}\mapsto (y_k)_{k\in \mathbb{Z}}\enskip \text{with}\enskip y_k =\begin{cases}
        0      & \text{if}\enskip -n\le k<n\\
        x_{k-n}& \text{if}\enskip k\ge n \\
        x_{k+n}  &   \text{if}\enskip k<-n\qquad\quad\enskip
        \end{cases}\\
        V_{-n}&:\enskip (x_k)_{k\in \mathbb{Z}}\mapsto (y_k)_{k\in \mathbb{Z}}\enskip \text{with}\enskip y_k =\begin{cases}
        x_{k+n}& \text{if}\enskip k\ge 0 \\
        x_{k-n}  &   \text{if}\enskip k<0\qquad\qquad\enskip\enskip
        \end{cases}\\
        U_n&:\enskip (x_k)_{k\in \mathbb{Z}}\mapsto (x_{k-n})_{k\in \mathbb{Z}}\\
        U_{-n}&:\enskip (x_k)_{k\in \mathbb{Z}}\mapsto (x_{k+n})_{k\in \mathbb{Z}}
        \end{align*}
Further, set $Q_n := I- P_n$.

We introduce the following $^*$-subalgebras of $L_{\mathcal{L}(H)}^\infty$. Let $PC_{\mathcal{L}(H)}$ denote the set of all $\mathcal{L}(H)$-valued piecewise continuous functions, i.e., functions $p\in L_{\mathcal{L}(H)}^\infty$ for which the one-sided limits
\begin{equation}
    p(\tau\pm 0):= \lim\limits_{t\to \tau\pm 0} p(t) =\lim\limits_{x\to\pm 0} p(\tau e^{ix})
\end{equation}
exist for all $\tau\in \mathbb{T}$, where the limit is taken in the operator norm of $\cL(H)$. It can be shown easily that
        \begin{equation}
        PC_{\mathcal{L}(H)} = \mathrm{clos}_{L_{\mathcal{L}(H)}^\infty}\left\{\sum_i a_ip_i: \enskip a_i\in \mathcal{L}(H), p_i\in PC\right\}.\label{PCLH}
        \end{equation}
By $C_{\mathcal{L}(H)}(\mathbb{T})$ we refer to the set of all $\mathcal{L}(H)$-valued continuous functions on $\mathbb{T}$, and denote by $H_{\mathcal{L}(H)}^\infty$ (resp., $\overline{H_{\mathcal{L}(H)}^\infty}$) the Hardy space consisting of all $a\in L_{\mathcal{L}(H)}^\infty$ whose Fourier coefficients $a_n$ vanish for all $n<0$ (resp., $n>0$).

The class $QC$ of (scalar) quasicontinuous functions is defined as
$QC:=(C(\T)+H^\infty)\cap (C(\T)+\overline{H^\infty})$. It is known that $QC$ is a $^*$-subalgebra of $L^\infty$ and that $a\in QC$ if and only if
both Hankel operators $H(a)$ and $H(\tilde{a})$ are compact.
As discussed in \cite{ehrhardt1996symbol},  there exist several possibilities of defining quasicontinuous functions in the $\mathcal{L}(H)$-valued setting. In particular, we introduce the following two:
        \begin{align}
        QC_{\mathcal{L}(H)} &:= (C_{\mathcal{L}(H)}(\mathbb{T}) + H_{\mathcal{L}(H)}^\infty) \cap (C_{\mathcal{L}(H)}(\mathbb{T}) + \overline{H_{\mathcal{L}(H)}^\infty}),\\
        QC_{\mathcal{L}(H)}^s &:= \mathrm{clos}_{L_{\mathcal{L}(H)}^\infty} \left\{\sum_i a_iq_i: \enskip a_i\in \mathcal{L}(H), q_i\in QC\right\}.
        \end{align}
Both $QC_{\mathcal{L}(H)}$ and $QC_{\mathcal{L}(H)}^s$ are $^*$-subalgebras of $L_{\mathcal{L}(H)}^\infty$, and $QC_{\mathcal{L}(H)}^s \subseteq QC_{\mathcal{L}(H)}$. The inclusion is proper if and only if $\dim H = \infty$. Finally, let $PQC_{\mathcal{L}(H)}$ be the smallest closed subalgebra of $L_{\mathcal{L}(H)}^\infty$ containing both $PC_{\mathcal{L}(H)}$ and $QC_{\mathcal{L}(H)}^s$.
It is easy to see that
\begin{equation} \label{PQCLH}
       PQC_{\mathcal{L}(H)} = \mathrm{clos}_{L_{\mathcal{L}(H)}^\infty}\left\{\sum_i a_ip_i q_i: \enskip a_i\in \mathcal{L}(H), p_i\in PC,\; q_i\in QC\right\}.
\end{equation}
We refer to $PQC_{\mathcal{L}(H)}$ as the class of $\cL(H)$-valued piecewise quasicontinuous functions.

Let $M(QC)$ denote the maximal ideal space of $QC$. It was shown in  \cite{ehrhardt1996symbol}
that $QC_{\mathcal{L}(H)}^s$ is \emph{locally trivial} at all points $\xi \in M(QC)$, i.e., for each $q\in QC_{\mathcal{L}(H)}^s$, there is an $a\in \mathcal{L}(H)$ such that $q-a \in I_{\xi, \mathcal{L}(H)}$, where
        \begin{equation}
        I_{\xi, \mathcal{L}(H)} := \mathrm{clos}\enskip \mathrm{id}_{QC_{\mathcal{L}(H)}^s}\{f: f\in QC, \xi(f)=0\}
        \end{equation}
is the smallest closed ideal in $QC_{\mathcal{L}(H)}^s$ containing the ``scalar ideal'' $\xi$.
Furthermore, $a$ is uniquely determinined and we will write
 \begin{equation}
        \Phi_\xi(q) = a.\label{locallytrivial}
 \end{equation}
Then the map $\Phi_\xi$ is a $^*$-homomorphism from $QC_{\mathcal{L}(H)}^s$ onto $\mathcal{L}(H)$.
As shown in \cite{ehrhardt1996symbol},  $QC_{\mathcal{L}(H)}^s$ is the largest $^*$-subalgebra of $QC_{\mathcal{L}(H)}$
which is locally trivial at each $\xi\in M(QC)$.

\section{A generalization of compactness and of strong convergence}\label{sec:3}
In order to study stability, we will follow a general scheme introduced in \cite{bottcher1990analysis}. Specifically, in the scalar case, this method relies heavily on the compactness of Hankel operators with continuous generating functions and on the notion of strong convergence. However, if $\dim H = \infty$, a Hankel operator of continuous $\mathcal{L}(H)$-valued function fails to be compact in general, and as a consequence we need a modification of ``compactness'' and  ``strong convergence''. The following definitions were given already in \cite{bottcher1996infinite} and were also used in \cite{ehrhardt1996finite}.

Let $\mathcal{K}$ denote the set of all operators $K\in \mathcal{L}(l_H^2(\mathbb{Z}))$ for which
        \[\|KQ_n\|_{\mathcal{L}(l_H^2(\mathbb{Z}))}\to 0, \quad \|Q_nK\|_{\mathcal{L}(l_H^2(\mathbb{Z}))}\to 0\]
as $n\to \infty$. Furthermore, let $\mathcal{A}$ stand for the set of all operators $A\in \mathcal{L}(l_H^2(\mathbb{Z}))$ for which both $AK\in \mathcal{K}$ and $KA\in \mathcal{K}$ whenever $K\in \mathcal{K}$.

The following basic properties were proved in \cite{ehrhardt1996finite}, Section 3.
\begin{prop} \label{prop21}\hfill
\begin{enumerate}[label=(\alph*)]
\item $\mathcal{A}$ is a $^*$-subalgebra of $\mathcal{L}(l_H^2(\mathbb{Z}))$, and $\mathcal{K}$ is a $^*$-ideal of $\mathcal{A}$.
\item $I, P, Q, J\in \mathcal{A}$. Moreover, $L(a)\in \mathcal{A}$ for all $a\in L_{\mathcal{L}(H)}^\infty$.
\item $P_nAP_n \in \mathcal{K}$ for all $A\in \mathcal{L}(l_H^2(\mathbb{Z}))$. In particular, $P_n\in \mathcal{K}$ and $W_n\in \mathcal{K}$.
\end{enumerate}
\end{prop}

The following proposition describes the connection between Hankel operators with quasicontinuous generating functions and the concept of  
``$Q_n$-compactness'' introduced above. It is the immediate consequence of the $\mathcal{L}(H)$-valued version of the Hartman Theorem (see \cite{bottcher1996infinite}, Proposition 3.2, and \cite{gohberg1994projection}).

\begin{prop}\label{prop22}
Let $f\in L_{\mathcal{L}(H)}^\infty$. Then $f\in QC_{\mathcal{L}(H)}$ if and only if $PfQ\in \mathcal{K}$ and $QfP\in \mathcal{K}$.
\end{prop}

Next we introduce the modified version of ``strong convergence" for operators contained in $\mathcal{A}$. Let $(A_n)_{n=1}^\infty$ be a sequence of operators $A_n\in \mathcal{A}$. We say that $A_n$ \emph{converges} $\mathcal{K}$-\emph{strongly} to an operator $A$, if, for all $K\in \mathcal{K}$, both
        \[\|K(A_n-A)\|_{\mathcal{L}(l_H^2(\mathbb{Z}))}\to 0 \enskip \mathrm{and}\enskip \|(A_n-A)K\|_{\mathcal{L}(l_H^2(\mathbb{Z}))}\to 0\]
as $n\to \infty$.
In this case we will also write
\begin{equation}
 A =\slim\limits_{n\to \infty} A_n.
\end{equation}

\begin{prop}\hfill
\begin{enumerate}[label=(\alph*)]
\item If $A_n \to A\enskip \mathcal{K}$-strongly, then $A\in \mathcal{A}$, and
        \[\|A\|_{\mathcal{A}} \le \liminf_{n\to\infty} \|A_n\|_\mathcal{A} \le \sup_{n\ge 1}\|A_n\|_\mathcal{A}<\infty.\]
\item If $A_n \to A$ and $B_n\to B \enskip \mathcal{K}$-strongly, then $A_nB_n\to AB$ and $A_n+B_n \to A+B \enskip\mathcal{K}$-strongly.
\item If $A_n \to A\enskip \mathcal{K}$-strongly and $\lambda\in \mathbb{C}$, then $A_n^*\to A^*$ and $\lambda A_n\to \lambda A \enskip \mathcal{K}$-strongly.
\item We have $Q_n\to 0, \enskip P_n\to I$, and $W_nLW_n\to 0\enskip \mathcal{K}$-strongly provided $L\in \mathcal{K}$.
\end{enumerate}
\end{prop}
\begin{proof}
See \cite{ehrhardt1996finite}, Propositions 3.7 and 3.8.
\end{proof}
\section{\Large \bf Algebraization of stability and lifting theorem}\label{s:4}
In this section we restate the stability problem in an algebraic language. More precisely, we construct a $C^*$-algebra $\mathcal{F}$, such that a sequence of operators is stable if and only if a specifically assigned element of $\mathcal{F}$ is invertible. Note that all the algebras constructed subsequently are $C^*$-algebras. Later, we consider a $^*$-subalgebra $\mathcal{F}_0$ of $\mathcal{F}$, to which we
will apply a  ``lifting theorem''.

To start with, let $\mathcal{F}$ be the set of all sequences $(A_n)_{n=1}^\infty$ of operators $A_n\in \mathcal{L}(l_H^2(\mathbb{Z}_n))$ for which
        \[\|(A_n)\|_\mathcal{F} := \sup_{n\ge 1}\|A_n\|_{ \mathcal{L}(l_H^2(\mathbb{Z}_n))} <\infty.\]
With the above norm and the algebraic operations
        \begin{align*}
        \lambda(A_n)& := (\lambda A_n), &(A_n)+(B_n) & := (A_n+B_n),\\
        (A_n)\cdot(B_n) & := (A_nB_n),&  (A_n)^* & := (A_n^*),
        \end{align*}
$\mathcal{F}$ is a $C^*$-algebra with the unit element $(P_n)_{n= 1}^\infty$.

Let $\mathcal{N}$ be the set of all sequences $(C_n)\in \mathcal{F}$ for which $\|C_n\|_{\mathcal{L}(l_H^2(\mathbb{Z}_n))}\to 0$ as $n\to \infty$. Apparently, $\mathcal{N}$ is a $^*$-ideal of $\mathcal{F}$, and hence the quotient algebra $\mathcal{F}/ \mathcal{N}$ is also a $C^*$-algebra. The following result is well-known (see, e.g., \cite{bottcher1990analysis,ehrhardt1996finite}) and easy to prove.

\begin{prop}\label{prop31}
Let $(A_n)\in \mathcal{F}$. Then $(A_n)$ is stable if and only if $(A_n)+\mathcal{N}$ is invertible in $\mathcal{F}/ \mathcal{N}$.
\end{prop}

Let $\mathcal{F}_0$ be the set of all elements $(A_n)\in \mathcal{F}$ for which the $\mathcal{K}$-strong limits

        \begin{align}
       \label{operatorP} \mathcal{P}(A_n) &:= \slim_{n\to\infty} P_nA_nP_n,\\
       \label{operatorW} \mathcal{W}(A_n) &:= \slim_{n\to\infty} W_nA_nW_n 
        \end{align}
exist. By $\mathcal{J}$ we denote the set
        \begin{equation}
        \mathcal{J} = \{(P_nKP_n + W_nLW_n +C_n)_{n=1}^\infty: K, L \in \mathcal{K}, (C_n)\in \mathcal{N}\}.
        \end{equation}
The following properties can be shown by straightforward computations (see also  \cite{ehrhardt1996finite}).

\begin{prop}\label{prop32}\hfill
\begin{enumerate}[label=(\alph*)]
\item $\mathcal{F}_0$ is a $^*$-subalgebra of $\mathcal{F}$.
\item $\mathcal{P}, \mathcal{W}$ are $^*$-homomorphisms from $\mathcal{F}_0$ into $\mathcal{A}$.
\item $\mathcal{J}$ is a $^*$-ideal of $\mathcal{F}_0$.
\end{enumerate}
\end{prop}

The preceding proposition ensures that $\mathcal{F}_0/\mathcal{J}$ is a $C^*$-algebra. Furthermore, $\mathcal{N}$ is a $^*$-ideal of $\mathcal{F}_0$, the quotient algebra $\mathcal{F}_0/\mathcal{N}$ is a $^*$-subalgebra of $\mathcal{F}/\mathcal{N}$ and hence inverse closed.

\begin{thm} {\bf (Lifting Theorem)}\label{thm33}
Let $(A_n)\in \mathcal{F}_0$. Then the following statements are equivalent:\hfill
\begin{enumerate}[label=(\alph*)]
\item The sequence $(A_n)$ is stable.
\item $(A_n) + \mathcal{N}$ is invertible in $\mathcal{F}_0/ \mathcal{N}$.
\item Both $\mathcal{P}(A_n)$ and $\mathcal{W}(A_n)$ are invertible in $\mathcal{A}$, and $(A_n) + \mathcal{J}$ is invertible in $\mathcal{F}_0/\mathcal{J}$.
\end{enumerate}
\end{thm}
\begin{proof}
(a) $\Leftrightarrow$ (b): This simply follows from Proposition \ref{prop31}. Note that $C^*$-algebras are inverse closed.\\
(b) $\Rightarrow$ (c): See Proposition \ref{prop32}. Observe that $\mathcal{N} \subseteq \ker\mathcal{P}$, $\mathcal{N} \subseteq\ker\mathcal{W}$, and $\mathcal{N} \subseteq\mathcal{J}$.\\
(c) $\Rightarrow$ (b): Let $(A_n)+\mathcal{J}$ be invertible in $\mathcal{F}_0/ \mathcal{J}$, and denote by $(B_n)+\mathcal{J}$ its right inverse, i.e.,
        \[A_nB_n = P_n + P_nKP_n + W_nLW_n + C_n\]
for some $K, L \in \mathcal{K}$ and $(C_n)\in \mathcal{N}$. Define a sequence $(B_n')$ by
        \[B_n' = B_n - P_n\mathcal{P}(A_n)^{-1}KP_n - W_n\mathcal{W}(A_n)^{-1}LW_n.\]
We see that $A_nB_n'$ equals to
        \begin{align*}
         &P_n + P_nKP_n + W_nLW_n + C_n - A_nP_n\mathcal{P}(A_n)^{-1}KP_n - A_nW_n\mathcal{W}(A_n)^{-1}LW_n\\
        &\enskip\enskip\enskip  =\enskip P_n + P_nC_n'P_n + W_nC_n''W_n+C_n,
        \end{align*}
where
	\begin{equation*}
		C_n' = (\mathcal{P}(A_n) - P_nA_nP_n)\mathcal{P}(A_n)^{-1}K,\enskip\enskip C_n''= (\mathcal{W}(A_n) - W_nA_nW_n)\mathcal{W}(A_n)^{-1}L
	\end{equation*}
both converge to zero in the norm. Therefore $(B_n') + \mathcal{N}$ is the right inverse of $(A_n) + \mathcal{N}$ in $\mathcal{F}_0/ \mathcal{N}$. Similarly,  $(A_n) + \mathcal{N}$ is left invertible in $\mathcal{F}_0/ \mathcal{N}$ as well.
\end{proof}

In order to proceed, we restrict again our considerations  to a smaller algebra $\mathcal{F}^J(PQC_{\mathcal{L}(H)})$ for which the final
stability result will be established.

As already mentioned in the introduction, let $\mathcal{S}(PQC_{\mathcal{L}(H)})$ stand for the smallest closed subalgebra of $\mathcal{L}(l_H^2(\mathbb{Z}))$ which contains all Laurent operators $L(a)$ with $a\in PQC_{\mathcal{L}(H)}$ and the operators $P$ and $Q$.
Further, refer to $\mathcal{S}^J(PQC_{\mathcal{L}(H)})$ the smallest closed subalgebra of $\mathcal{L}(l_H^2(\mathbb{Z}))$ containing the same operators and in addition the flip operator $J$.

 By $\mathcal{F}(PQC_{\mathcal{L}(H)})$ we denote the smallest closed subalgebra of $\mathcal{F}$ which includes all sequences
 $(P_nAP_n)$ with $A\in \mathcal{S}(PQC_{\mathcal{L}(H)})$ and all sequences contained in $\cJ$. Furthermore, refer to $\mathcal{F}^J(PQC_{\mathcal{L}(H)})$ as the smallest closed subalgebra of $\mathcal{F}$ which contains all sequences $(P_nAP_n)$ with $A\in \mathcal{S}^J(PQC_{\mathcal{L}(H)})$ and all sequences contained in $\cJ$.

\begin{prop}\label{prop34} $\mathcal{S}(PQC_{\mathcal{L}(H)})$  and $\mathcal{S}^J(PQC_{\mathcal{L}(H)})$ are $^*$-subalgebra of $\mathcal{A}$, and $\mathcal{K}$ is a $^*$-ideal of both $\mathcal{S}(PQC_{\mathcal{L}(H)})$ and  $\mathcal{S}^J(PQC_{\mathcal{L}(H)})$.
\end{prop}
\begin{proof}
Note that by Proposition \ref{prop21}, it suffices to show $\mathcal{K} \subseteq \mathcal{S}(PQC_{\mathcal{L}(H)})$ ($\subseteq \mathcal{S}^J(PQC_{\mathcal{L}(H)})$). Let $K_{ij,a}$ denote the operator for which the $(i,j)$-entry of its matrix representation is equal to $a\in \mathcal{L}(H)$ and elsewhere zero. Since
        \[K_{ij,a} = a\cdot U_i(P-U_1PU_{-1})U_{-j},\]
and $U_k = L(t^k)$, we have $K_{ij,a} \in \mathcal{S}(PQC_{\mathcal{L}(H)})$. Further, every $K\in \mathcal{K}$ can be approximated arbitrarily close by operators $P_nKP_n$, which are finite linear combinations of operators of the form $K_{ij,a}$. It implies $K\in \mathcal{S}(PQC_{\mathcal{L}(H)})$.
\end{proof}

Regarding $\mathcal{S}^J(PQC_{\mathcal{L}(H)})$ we have the following results.

\begin{prop}\label{prop35} For all $A\in \mathcal{S}^J(PQC_{\mathcal{L}(H)})$, the $\mathcal{K}$-strong limit
        \begin{equation}\label{eq4.3}
        \mathcal{U}(A):=\slim_{n\to \infty} \left( \begin{array}{cc}
        U_{-n}PAPU_n & U_{-n}PAQU_{-n}\\
        U_nQAPU_n & U_nQAQU_{-n}\\
        \end{array} \right)
        \end{equation}
exists. Furthermore, the mapping $\mathcal{U}$ is a $^*$-homomorphism from $\mathcal{S}^J(PQC_{\mathcal{L}(H)})$ into the $C^*$-algebra $\mathcal{A}^{2\times2}$. We have $\mathcal{K}\subseteq\ker \mathcal{U}$, and $\mathcal{U}$ acts on the generating elements of $\mathcal{S}^J(PQC_{\mathcal{L}(H)})$ as follows:
        \[ \mathcal{U}(P)=   \left( \begin{array}{cc}
        I & 0\\
        0 & 0\\
        \end{array} \right),\enskip
        \mathcal{U}(f)=   \left( \begin{array}{cc}
        f & 0\\
        0 & f\\
        \end{array} \right),\enskip
        \mathcal{U}(J)=   \left( \begin{array}{cc}
        0 & J\\
        J & 0\\
        \end{array} \right)\]
\end{prop}
\begin{proof}
In \cite[Proposition 5.2]{ehrhardt1996finite}, these statements were shown for $\mathcal{U}$ acting on $\mathcal{S}(PQC_{\mathcal{L}(H)})$.
In view of this, all what remains to be verified is that the $\mathcal{K}$-strong limit \eqref{eq4.3} exists for $A=J$. Note that $PJP = QJQ = 0$, 
and 
\begin{align*}
U_{-n}PJQU_{-n} &=(U_{-n}P U_{n}) J (U_{n} Q U_{-n})\to J,
\\
U_nQJPU_n &= (U_{n} Q U_{-n}) J ( U_{-n} P U_n)\to J
\end{align*}
$\mathcal{K}$-strongly, since both $U_{-n}PU_n=I-Q_nP$ and $U_{n}Q U_{-n}=I-Q_nP$ converge $\mathcal{K}$-strongly
to $I$. From this the desired assertions simply follow.
\end{proof}

The next result is a generalization of Corollary 5.4 of \cite{ehrhardt1996finite}  from the case
 $\mathcal{F}(PQC_{\mathcal{L}(H)})$ to the case
$\mathcal{F}^J(PQC_{\mathcal{L}(H)})$.
\begin{prop}\label{prop36}\hfill
\begin{enumerate}[label=(\alph*)]
\item $\mathcal{F}^J(PQC_{\mathcal{L}(H)})$ is a $^*$-subalgebra of $\mathcal{F}_0$, and $\mathcal{J}$ is a $^*$-ideal of $\mathcal{F}^J(PQC_{\mathcal{L}(H)})$.
\item On the generating elements of $\mathcal{F}^J(PQC_{\mathcal{L}(H)})$, $\mathcal{P}$ and $\mathcal{W}$ act as follows:
        \begin{gather}
        \mathcal{P}(P_nAP_n) = A,\label{P}\\
        \mathcal{W}(P_nAP_n)= (PJ, \; QJ)\enskip\mathcal{U}(A)   \left( \begin{array}{cc}
        JP\\
        JQ\\
        \end{array} \right)
        \end{gather}
for $A\in \mathcal{S}^J(PQC_{\mathcal{L}(H)})$, and
        \begin{equation}
        \mathcal{P}(A_n) = K, \qquad \mathcal{W}(A_n) = L
        \end{equation}
for $(A_n) = (P_nKP_n + W_nLW_n + C_n)\in \mathcal{J}$.
\item $\mathcal{P}$ is a $^*$-homomorphism from $\mathcal{F}^J(PQC_{\mathcal{L}(H)})$ onto $\mathcal{S}^J(PQC_{\mathcal{L}(H)})$.
\item
$\mathcal{W}$ is a $^*$-homomorphism from $\mathcal{F}^J(PQC_{\mathcal{L}(H)})$ into $\mathcal{S}^J(PQC_{\mathcal{L}(H)})$.
\end{enumerate}
\end{prop}
\begin{proof}
The proof is carried out in the same way as the proof of the Corollary 5.4 of \cite{ehrhardt1996finite}. The crucial point is to show (b).
The statements (a), (c) and (d) either follow from it or are obvious. The only non-trivial statement in (b) is the existence and formula for
$\mathcal{W}(P_n A P_n)$. For its proof a formula identical to one in Corollary 5.3 of \cite{ehrhardt1996finite} can be used, by which the 
assertion reduces to the existence of the limit $\mathcal{U}(A)$ for $A\in \mathcal{S}^J(PQC_{\mathcal{L}(H)})$. But this is clear because of
Proposition \ref{prop35}.
\end{proof}

From Theorem \ref{thm33}, Proposition \ref{prop34} and Proposition \ref{prop36}	, and the fact that $C^*$-algebras are inverse closed, we obtain the following:
\begin{thm}{\bf (Lifting Theorem II)}\label{thm37} Let $(A_n)\in \mathcal{F}^J(PQC_{\mathcal{L}(H)})$. The following statements are equivalent:
\begin{enumerate}[label=(\alph*)]
\item The sequence $(A_n)$ is stable.
\item $\mathcal{P}(A_n)$ and $\mathcal{W}(A_n)$ are invertible in $\mathcal{S}^J(PQC_{\mathcal{L}(H)})$, and $(A_n) + \mathcal{J}$ is invertible in $\mathcal{F}^J(PQC_{\mathcal{L}(H)})/\mathcal{J}$.
\end{enumerate}
\end{thm}
\section{Localization}\label{s:5}
In view of Theorem \ref{thm37}, the question now is when is $(A_n) + \mathcal{J}$ invertible. We will investigate the structure of the $C^*$-algebra $\mathcal{F}^J(PQC_{\mathcal{L}(H)})/\mathcal{J}$ to answer this problem. Fortunately, this algebra possesses a sufficiently large center, thus the ``Local Principle'' by  Allan/Douglas (see \cite{douglas1972banach}, Proposition 4.5, or \cite{bottcher1990analysis}) can be used.

\begin{thm}{\bf (Local principle by Allan/Douglas)}\label{thm41}
Let $\mathcal{B}$ be a unital $C^*$-algebra, and let $\mathcal{C}$ be a central subalgebra of $\mathcal{B}$, i.e., a closed $^*$-subalgebra of the center of $\mathcal{B}$ which contains the identity element. For each maximal ideal $\xi$ of $\mathcal{C}$, let $I_\xi$ denote the smallest closed two-sided ideal of $\mathcal{B}$ which contains $\xi$. Then an element $b$ of $\mathcal{B}$ is invertible if and only if the cosets $b+I_\xi$ are invertible in $\mathcal{B}/ I_\xi$ for all $\xi$.
\end{thm}

A proof of the local principle can be found in \cite[Section 1.34]{bottcher1990analysis}. In order to apply it, we still need couple preliminary results. Denote by $\widetilde{QC}$ the space of all even (scalar) quasicontinuous functions, i.e., all $f\in QC$ such that $f=\tilde{f}$, where $\tilde{f}(t) := f(1/t), t\in \T$ (see also \eqref{tilde}).

\begin{lem}\label{lem42} Let $A\in \mathcal{S}^J(PQC_{\mathcal{L}(H)})$. Then\hfill
\begin{enumerate}[label=(\alph*)]
\item $Af-fA\in \mathcal{K}$ for all $f\in \widetilde{QC}$,
\item $(P_nAQ_nfP_n), (P_nfQ_nAP_n)\in \mathcal{J}$ for all $f\in QC_{\mathcal{L}(H)}$.
\end{enumerate}
\end{lem}
\begin{proof}
(a):\ Since $\mathcal{K}$ is an $^*$-ideal of $\mathcal{S}^J(PQC_{\mathcal{L}(H)})$, it suffices to check the assertion for all the generating elements of $\mathcal{S}^J(PQC_{\mathcal{L}(H)})$. For $A = g\in PQC_{\mathcal{L}(H)}$, we have $fg-gf = 0$ since $f$ is scalar. For $A=P$, $Pf-fP = PfQ - QfP\in \mathcal{K}$ by Proposition \ref{prop22}. Finally, since $f = \tilde{f}$ is even, $Jf - fJ = 0$.

\enskip(b): We show that $(P_nAQ_nfP_n)\in \mathcal{J}$. In fact,
        \begin{align*}
        &V_{-n}(PfP+QfQ)W_n \\
        =\enskip& (PU_{-n}P+QU_nQ)(PfP+QfQ)(PU_nJP+QU_{-n}JQ)P_n \\
        =\enskip& (PfJP+QfJQ)P_n.
        \end{align*}
It follows that
        \begin{align*}
        P_nAQ_nfP_n &= P_nAV_nV_{-n}(PfP+QfQ)P_n + P_nAQ_n(PfQ+QfP)P_n \\
        & = W_n(W_nAV_n)(PfQ+QfP)JW_n + P_nAQ_n(PfQ+QfP)P_n.
        \end{align*}
Furthermore, $PfQ+QfP\in \mathcal{K}$, $Q_n\to 0$ and $W_nAV_n\to \tilde{A}$ $\mathcal{K}-$strongly, where $\tilde{A}\in \mathcal{A}$ is a certain operator. Therefore,
        \[P_nAQ_nfP_n = W_nLW_n+ C_n,\]
where $(C_n)\in \mathcal{N}$ and $L = \tilde{A}(PfQ+QfP)J\in \mathcal{K}$. Hence $(P_nAQ_nfP_n)\in \mathcal{J}$. The case for $(P_nfQ_nAP_n)$ can be shown analogously.
\end{proof}

\begin{lem}\label{lem43} The set $\mathcal{C} = \{(P_nfP_n)+\mathcal{J}: f\in \widetilde{QC}\}$ is a $^*$-subalgebra contained in the center of $\mathcal{F}^J(PQC_{\mathcal{L}(H)})/\mathcal{J}$. Moreover, $\mathcal{C}$ is *-isomorphic to $\widetilde{QC}$.
\end{lem}
\begin{proof}
We first show that $\mathcal{C}$ is contained in the center of $\mathcal{F}^J(PQC_{\mathcal{L}(H)})/\mathcal{J}$. 
It suffices to prove that $(P_nfP_n)$ commutes with $(P_nAP_n)$ for all $A\in \mathcal{S}^J(PQC_{\mathcal{L}(H)})$ and $f\in \widetilde{QC}$ modulo $\mathcal{J}$. In fact,
        \[P_nfP_nAP_n - P_nAP_nfP_n = P_n(fA - Af)P_n - P_nfQ_nAP_n + P_nAQ_nfP_n,\]
and therefore $(P_nfP_nAP_n - P_nAP_nfP_n)$ belongs to $\mathcal{J}$ by the previous lemma.

To show the remaining, construct the $^*$-homomorphism mapping
        \begin{equation}
        \Phi: \widetilde{QC}\to \mathcal{C}, \enskip f\mapsto (P_nfP_n)+\mathcal{J}.
        \end{equation}
Following  the proof of Lemma 6.3 in  \cite{ehrhardt1996finite}, we obtain that the
kernel of $\Phi$ is trivial. By definition, the image of $\Phi$ is $\mathcal{C}$ which is therefore a (closed) $^*$-subalgebra of $\mathcal{F}^J(PQC_{\mathcal{L}(H)})/\mathcal{J}$ isomorphic to  $\widetilde{QC}$.
\end{proof}

Let $M(\widetilde{QC})$ stand for the maximal ideal space of $\widetilde{QC}$.
For $\eta\in M(\widetilde{QC})$, we denote by $\mathcal{J}_{\eta, \mathcal{L}(H)}^J$ the smallest closed two-sided ideal of $\mathcal{F}^J(PQC_{\mathcal{L}(H)})$ which contains $\mathcal{J}$ and all sequences $(P_nfP_n)$ with $f\in \widetilde{QC}$ and $\eta(f) = 0$, i.e.,
        \begin{equation}
        \mathcal{J}_{\eta, \mathcal{L}(H)}^J := \mathrm{clos}\enskip \mathrm{id}_{\mathcal{F}^J(PQC_{\mathcal{L}(H)})} (\mathcal{J}\cup \{(P_n f P_n):\enskip f\in \widetilde{QC}, \enskip\eta(f) = 0\}).\label{Jwithaflip}
        \end{equation}
The quotient algebra $\mathcal{F}^J(PQC_{\mathcal{L}(H)})/\mathcal{J}^J_{\eta, \mathcal{L}(H)}$ is denoted by $\mathcal{F}^J_\eta(PQC_{\mathcal{L}(H)})$.

By applying Theorem \ref{thm41} and Lemma \ref{lem43} we immediately obtain the following:

\begin{cor}\label{cor44} Let $(A_n)\in \mathcal{F}^J(PQC_{\mathcal{L}(H)})$. Then $(A_n)+\mathcal{J}$ is invertible in $\mathcal{F}^J(PQC_{\mathcal{L}(H)})/\mathcal{J}$ if and only if $(A_n)+ \mathcal{J}^J_{\eta, \mathcal{L}(H)}$ is invertible in $\mathcal{F}_\eta^J(PQC_{\mathcal{L}(H)})$ for all $\eta \in M(\widetilde{QC})$.
\end{cor}

For $\eta \in M(\widetilde{QC})$, let $\mathcal{I}^J_{\eta, \mathcal{L}(H)}$ refer to the smallest closed ideal of $\mathcal{S}^J(PQC_{\mathcal{L}(H)})$ containing $\mathcal{K}$ and all operators $f\in \widetilde{QC}$ with $\eta(f) = 0$, i.e.,
        \begin{equation}
        \mathcal{I}^J_{\eta, \mathcal{L}(H)} := \mathrm{clos}\enskip \mathrm{id}_{\mathcal{S}^J(PQC_{\mathcal{L}(H)})} (\mathcal{K}\cup \{f:\enskip f\in \widetilde{QC}, \enskip\eta(f) = 0\}).
        \end{equation}
The quotient algebra $\mathcal{S}^J(PQC_{\mathcal{L}(H)})/\mathcal{I}^J_{\eta, \mathcal{L}(H)}$ is denoted by $\mathcal{S}^J_\eta(PQC_{\mathcal{L}(H)})$.

\begin{lem}\label{lem45} $(P_nBP_n)\in \mathcal{J}^J_{\eta, \mathcal{L}(H)}$ whenever $B\in \mathcal{I}^J_{\eta, \mathcal{L}(H)}$ and $\eta\in M(\widetilde{QC})$.
\end{lem}
\begin{proof}
By Lemma \ref{lem42}(a), we see that $B$ can be approximated arbitrarily close by operators of the form $\sum\limits_i A_if_i + K$, where $A_i\in \mathcal{S}^J(PQC_{\mathcal{L}(H)}), f_i\in \widetilde{QC}, \eta(f_i) = 0$ and $K\in \mathcal{K}$. Using Lemma \ref{lem42}(b), it directly follows from computations that
                $$(P_n(\sum_i A_if_i + K)P_n) = \sum_i(P_nA_iP_n)(P_nf_iP_n)\enskip \mathrm{mod}\enskip \mathcal{J}.$$
The right side of the above formula is contained in $\mathcal{J}^J_{\eta, \mathcal{L}(H)}$. Hence $(P_nBP_n)\in \mathcal{J}^J_{\eta, \mathcal{L}(H)}$ by approximation.
\end{proof}

The next step would be to  study the local algebras $\mathcal{F}^J_\eta(PQC_\LH)$ for $\eta\in M(\widetilde{QC})$.
This will be done in Sections \ref{s:7}--\ref{s:9} and is more difficult.
Before we are able to do this, we need some results about the maximal ideal space $M(\widetilde{QC})$.

\section{The maximal ideal space ${M(\widetilde{QC})}$}\label{s:6}

In this section we recall some results about the maximal ideal space $M(\widetilde{QC})$ which have been established recently by the authors in \cite{MevenQC}.
Note that the classical results about the maximal ideals space $M(QC)$ go back to Sarason \cite{sarason1975functions, sarason1977} (see also \cite{bottcher1990analysis}).

If $\frak{B}$ is a $^*$-subalgebra of a commutative $C^*$-algebra $\frak A$, the restriction map from $M({\frak A})$ to $M({\frak B})$ is surjective (see, e.g.,  Proposition 1.26(b) in \cite{bottcher1990analysis}). For $\beta \in M({\frak B})$ the fiber of $M({\frak A})$ over $\beta$ is defined by
$$
M_\beta({\frak A})=\{ \, \alpha \in M({\frak A})\,:\, \alpha|_{\frak B}=\beta\,\}.
$$
The fibers $M_\beta({\frak A})$ are non-empty compact subsets of $M({\frak A})$, and $M({\frak A})$ is the disjoint union of all $M_\beta({\frak A})$.

Since one has corresponding embeddings of various $C^*$-algebras as shown in the first diagram below, what has just been said implies
natural (continuous) restriction maps between their maximal ideal spaces as shown in the second diagram:
\begin{center}
\begin{tikzpicture}[thick,scale=1]
 \matrix (m) [matrix of math nodes,row sep=2em,column sep=4em,minimum width=3em]
  {
     PQC & QC & \tQC& M(PQC) &  M(QC)  & M(\widetilde{QC})\\
     PC & C(\T) & \widetilde{C}(\T)  & M(PC)&  M(C)\cong {\T }& M(\widetilde{C})\cong \overline{\T_+} \\};
  \path[-stealth]
        (m-1-2) edge node [above] {}  (m-1-1)
        (m-1-3) edge node [above] {} (m-1-2)
        (m-2-1) edge node [left]  {} (m-1-1)
        (m-2-2) edge node [left]  {} (m-1-2)
        (m-2-3) edge node [left]  {} (m-1-3)
        (m-2-2) edge node [above] {} (m-2-1)
        (m-2-3) edge node [above] {} (m-2-2) 
        (m-1-4) edge node [above] {}  (m-1-5)
        (m-1-5) edge node [above] {\footnotesize{$\Psi$}} (m-1-6)
        (m-1-4) edge node [left]  {} (m-2-4)
        (m-1-5) edge node [left]  {} (m-2-5)
        (m-1-6) edge node [left]  {} (m-2-6)        
        (m-2-4) edge node [above] {}  (m-2-5)
        (m-2-5) edge node [above] {\footnotesize{$\Psi'$}} (m-2-6);
\end{tikzpicture}
\end{center}
Recall that $C(\T)$ is the $C^*$-algebra of continuous functions on $\T$, and $\widetilde{C}(\T)$ is the $C^*$-algebra
 of all even continous functions. 
Here $\overline{\T_+}:=\{ t\in \T\,:\, \mbox{Im}(t)\ge 0\}$. Trivially, the map $\Psi'$ is defined such that the pre-image
of $\tau\in\overline{ \T_+}$ equals the set $\{\tau,\overline{\tau}\}$, which consists of one or two points.
We also note that $M(PC)$ can be identified with the set $\T\times\{-1,+1\}$.

Corresponding to the `vertical' restriction maps we have the following fibers spaces:
       \begin{align*}
        M_\tau(QC) &= \{\xi\in M(QC): \xi(f) = f(\tau), \forall f\in C(\mathbb{T})\}
        			 \qquad \mbox{ over } \tau \in \T, \\
        M_\tau(\widetilde{QC}) &= \{\eta\in M(\widetilde{QC}): \eta(f) = f(\tau), \forall f\in \widetilde{C}(\mathbb{T}) \}
        			\qquad \mbox{ over } \tau \in \overline{\T_+}.
        \end{align*}
Repeating what was said above in general, these are non-empty compact sets, and 
	\begin{equation}\label{decom of QC1}
    		M(QC) =\bigcup_{\tau\in\T} M_{\tau}(QC),\quad
   		M(\tQC)=\bigcup_{\tau\in\overline{\T_+}} M_\tau (\tQC)
	\end{equation}
are disjoint unions.
Following Sarason, define
$$
M_\tau^\pm(QC) =\Big\{
\, \xi\in M(QC)\;:\; \xi(f)=0\;\mbox{ whenever } \limsup_{t\to \tau\pm 0} |f(t)|=0 \mbox{ and } f\in QC \, \Big\},
$$
which are closed subsets of $M_\tau(QC)$. Sarason also introduced $M^0(QC)$ whose definition requires some preparations.

Let $\frak A$ be a $C^*$-subalgebra of $L^\infty$, $\Lambda:=[1,\infty)$, and let $\{k_\lambda\}_{\lambda\in \Lambda}$ be an approximate identity generated by $K$ in the sense of Section 3.14 in \cite{bottcher1990analysis}. Then, each pair $ (\lambda, \tau )\in \Lambda\times\mathbb{T}$ induces a functional $\delta_{\lambda,\tau}\in \frak A^*$ given by
        \begin{equation}
        \delta_{\lambda,\tau}: {\frak A}\to \mathbb{C}, \enskip a\mapsto (k_\lambda a)(\tau),\label{deltafunction}
        \end{equation}
where
$$
   (k_\lambda a)(e^{i\theta})= \int_{-\infty}^\infty \lambda K(\lambda x) a(e^{i(\theta-x)})\, dx.
$$
Therefore $\Lambda\times \mathbb{T}$ can be regarded as a subset of $A^*$.
Examples of approximate identities in the above sense are the moving average
\begin{equation}\label{mlambda}
        (m_\lambda a)(e^{i\theta}) = \frac{\lambda}{2\pi}\int_{\theta-\pi/\lambda}^{\theta+\pi/\lambda}a(e^{ix})\,dx,
 \end{equation}
or the Poisson kernel
\begin{equation}\label{mlambda}
        (\varrho_r a)(e^{i\theta}) = 
        \frac{1}{2\pi}\int_{0}^{2\pi} \frac{1-r^2}{1+r^2-2r \cos (x-\theta)}\, a(e^{ix})\,dx.
 \end{equation}
To make a connection with $QC$, note that we have  the following result.
Therein, the dual space $QC^*$ is equipped with the weak-$^*$ topology (see \cite[Prop.~3.29]{bottcher1990analysis}).

\begin{prop}\label{p2.1}
$M(QC)= (\clos_{\,QC^*}(\Lambda\times \mathbb{T}))\setminus (\Lambda \times \mathbb{T}).$
\end{prop}

Now we are able to define
\begin{equation}
        M_\tau^0(QC) = M(QC)\cap  \clos_{\,QC^*}(\Lambda \times \{\tau\}) .\label{Mtau0}
\end{equation}
Here any approximate identities (in the sense of Section 3.14 in \cite{bottcher1990analysis}) can be used (see \cite[Lemma~3.31]{bottcher1990analysis}).
Clearly,  $M_\tau^0(QC)$ is a compact subset of the fiber $M_\tau(QC)$. 

The following result was originally proved by Sarason \cite{sarason1975functions}
(see also  \cite[Proposition~3.34]{bottcher1990analysis}).

\begin{prop}{\bf (Sarason)}\label{prop52}
$\enskip$ If $\tau\in \mathbb{T}$, then
        \begin{equation}
        M_\tau^0(QC) = M_\tau^+(QC)\cap M_\tau^-(QC), \quad  M_\tau^+(QC)\cup M_\tau^-(QC) = M_\tau(QC).
        \end{equation}
\end{prop}

Now, define $\chi_+$ (resp., $\chi_-$) as the characteristic function of the upper (resp., lower) semi-circle,
i.e., 
\begin{equation}\label{charfct}
\chi_{\pm}(t) 
=\left\{\begin{array}{ll} 1 & \mbox{ if } |t|=1,\; \pm \mbox{Im}(t)\ge0
\\[.5ex] 
 0 & \mbox{ if } |t|=1,\; \pm \mbox{Im}(t)<0.
\end{array}\right.
\end{equation}
 The following properties of quasicontinuous functions are needed subsequently. We refer to \cite{MevenQC} for its proof.

\begin{prop}\label{QCprop}
Let $q\in QC$.
\begin{enumerate}[label=(\alph*)]
\item If $q$ is an odd function, i.e., $q(t) = -q(1/t)$, then $q|M_{\pm1}^0(QC) = 0$.
\item For $\tau=\pm1$, if  $q|M_\tau^0(QC) = 0$ and  $p\in PC\cap C(\mathbb{T}\setminus\{\tau\})$, 
then $qp\in QC$.
\item
In particular, if $q|M_{1}^0(QC)=0$ and $q|M_{-1}^0(QC)=0$ , then $q\chi_+,q\chi_-\in QC$.
\end{enumerate}
\end{prop}

Let 
	\begin{equation}
		\Psi: M(QC) \to M(\widetilde{QC}), \enskip \xi \mapsto \hat{\xi} := \xi|_{\widetilde{QC}}
	\end{equation}
be the (surjective) map shown in the diagram given above. Now we proceed to analyze the fibers of $M(QC)$ over $\eta\in M(\widetilde{QC})$, which can be defined as
	\begin{equation}
		M^\eta(QC) = \{\xi\in M(QC): \hat{\xi} = \eta\}.
	\end{equation}
To prepare for it, for a given $\xi \in M(QC)$, we define its "conjugate" $\xi'\in M(QC)$ by 
	\[\xi'(q):= \xi (\tilde{q}), \enskip q\in QC.\]
Here recall the definition \eqref{tilde}.	
It is clear that $\hat{\xi} = \hat{\xi'}\in M(\widetilde{QC})$. Furthermore, the following statements are obvious:
\begin{enumerate}
\item[(i)]
If $\xi\in M_\tau(QC)$, then $\xi'\in M_{\bar{\tau}}(QC)$.
\item[(ii)]
If $\xi\in M_\tau^\pm(QC)$, then $\xi'\in M_{\bar{\tau}}^\mp(QC)$.
\item[(iii)]
If $\xi\in M_\tau^0(QC)$, then $\xi'\in M_{\bar{\tau}}^0(QC)$.
\end{enumerate}

Consider functionals $\delta_{\lambda, \tau}\in \widetilde{QC}^*$ associated with the moving average $\{m_\lambda\}$ given by (\ref{mlambda}),
	\[\delta_{\lambda,\tau}: a\in \widetilde{QC} \mapsto (m_\lambda a)(\tau), \enskip (\delta, \tau)\in \Lambda\times \T.\]
In analogy to (\ref{Mtau0}), define
	\begin{equation}\label{eqn-6.10}
        		M_\tau^0(\widetilde{QC}) := 
		 M(\widetilde{QC}) \cap
		\mathrm{clos}_{\widetilde{QC}^*}(\Lambda\times \{\tau\}), \enskip \tau\in \overline{\mathbb{T}_+}.
        	\end{equation}

To procced further, we have to distingish whether $\eta\in M_\tau(\tQC)$  with
$\tau\in\{+1,-1\}$ or with
$\tau\in \T_+:=\{ \tau \in \T\,:\, \mbox{Im}(\tau)> 0\}$.
Note that this relates to (\ref{decom of QC1}). \\

\subsection{Fibers over $\boldsymbol{M_{\tau}(\tQC)}$, $\boldsymbol{\tau\in\{+1,-1\}}$}

For the description of $M^\eta(QC)$ with $\eta\in M_{\pm 1}(\widetilde{QC})$, the following property is crucial.

\begin{prop}\label{p3.1}
If $\xi_1,\xi_2\in M_{\pm 1}^+(QC)$ and $\hat{\xi_1}=\hat{\xi_2}$, then $\xi_1=\xi_2$.
\end{prop}
\begin{proof}
Each $q\in QC$ admits a unique decomposition
$$
q = \frac{q+\tilde{q}}{2} + \frac{q - \tilde{q}}{2} =: q_e + q_o,
$$
where $q_e$ is even and $q_o$ is odd.
By Proposition \ref{QCprop}(c),  we have $q_o\chi_-\in QC$, and
\begin{align*}
        \xi_1(q) 
        &= \xi_1(q_e)+ \xi_1(q_o) = \xi_1(q_e) + \xi_1(q_o - 2q_o\chi_-)\\
        &= \eta(q_e) + \eta(q_o-2q_o\chi_-) \\
        &=  \xi_2(q_e) + \xi_2(q_o - 2q_o\chi_-)= \xi_2(q_e) + \xi_2(q_o) = \xi_2(q).
\end{align*}
Note that $q_o-2q_o\chi_-=q_o(\chi_+-\chi_-)\in \widetilde{QC}$ and that $\lim\limits_{t\to \pm1+0} q_o(t)\chi_-(t)=0$, 
whence $\xi_i(q_o\chi_-)=0$. It follows that $\xi_1=\xi_2$. 
\end{proof}

In other words, for $\eta\in M_{\pm1}(\widetilde{QC})$, there are only two possibilities for $M^\eta(QC)$: \hfill
\begin{enumerate}
\item[(a)] 
$M^\eta(QC) = \{\xi\}$ with $\xi=\xi' \in M_{\pm 1}^0(QC)$, or
\item[(b)] 
$M^\eta(QC) = \{\xi, \xi'\}$ with $$\xi\in M_{\pm 1}^+(QC)\setminus M_{\pm1}^-(QC)
\quad\mbox{and}\quad
\xi'\in M_{\pm 1}^-(QC)\setminus M_{\pm 1}^+(QC).$$
\end{enumerate}
This leads to the following characterization, which was proved in \cite[Theorems 3.2 and 3.3]{MevenQC}.

\begin{thm}\label{thm56}
$M_{\pm 1}^0(\tQC)$ is a closed subset of $M_{\pm 1}(\tQC)$. Moreover, 
\begin{enumerate}
\item[(a)] 
if $\eta\in M_{\pm 1}^0(\tQC)$, then $M^\eta(QC) = \{\xi\}$ with  $\xi=\xi' \in M_{\pm 1}^0(QC)$;
\item[(b)] 
if $\eta\in M_{\pm 1 }(\tQC)\setminus M_{\pm 1}^0(\tQC)$, then 
$M^\eta(QC) = \{\xi, \xi'\}$  such that 
	\[ \xi\in M_{\pm 1}^+(QC)\setminus M_{\pm 1}^-(QC)\enskip \text{and} \enskip \xi'\in M_{\pm 1}^-(QC)\setminus M_{\pm 1}^+(QC).\]
\end{enumerate}
\end{thm}

\subsection{Fibers over $\boldsymbol{M_\tau(\tQC)}$, $\boldsymbol{\tau\in \T_+}$}
Now we consider the fibers of $M^\eta(QC)$ over $\eta\in M_\tau(\tQC)$ with $\tau\in\T_+$.

\begin{prop}\label{prop59} If $\hat{\xi}_1 = \hat{\xi}_2$ for $\xi_1, \xi_2\in M_\tau(QC)$ with $\tau\in\T_+$, then $\xi_1 = \xi_2$.
\end{prop}
\begin{proof}
Otherwise, there exists a $q\in QC$ such that $\xi_1(q)\neq 0$, $\xi_2(q) = 0$. Since $\tau \in  \T_+$, one can choose a smooth function $c_\tau$ such that $c_\tau = 1$ in a neighborhood of $\tau$ and  it vanishes on the lower semi-circle. Define 
$\overline{q} = qc_\tau + \widetilde{qc_\tau}\in\tQC$, and note that  $\overline{q} - q$ is continuous at $\tau$ and vanishes there, hence 
$\xi_1(\overline{q} - q) = \xi_2(\overline{q} - q) = 0$. 
But then, we have	
	\[   0\neq \xi_1(q) = \xi_1(\overline{q}) = \xi_2(\overline{q}) = \xi_2(q) = 0 \]
since $\overline{q}\in\tQC$ and $\hat{\xi}_1 = \hat{\xi}_2$, which is a contradiction.
\end{proof}

Combined with the statements (i)-(iii) listed previously, the above proposition implies that for any $\eta\in M_\tau(\widetilde{QC})$ with $\tau\in \mathbb{T}_+$,  $M^\eta(QC)=\{\xi,\xi'\}$ with some (unique) $\xi\in M_\tau(QC)$.  This suggest to define $M_\tau^\pm(\widetilde{QC})$  as follows:
	\begin{equation}
		M_\tau^\pm(\widetilde{QC}):=\{ \, \hat\xi \,:\, \xi  \in M_\tau^\pm (QC)\,\},\enskip \tau\in \T_+ \,.
	\end{equation}
The following proposition characterizes the structure of $M_\tau(\widetilde{QC})$ in a way similar to Proposition \ref{prop52}. We refer to \cite[Prop.\  3.7]{MevenQC} for its proof.

\begin{prop}\label{prop510}
For $\tau\in \T_+$, we have
$$
M_\tau(\widetilde{QC})=M_\tau^+(\widetilde{QC})\cup M_\tau^-(\widetilde{QC}),\qquad
M_\tau^0(\widetilde{QC})=M_\tau^+(\widetilde{QC})\cap M_\tau^-(\widetilde{QC}).
$$
\end{prop}

\bigskip

To summarize the structure of $M(\widetilde{QC})$, we have the following disjoint unions:
\begin{align}
\label{dec-1}
 M(\widetilde{QC}) &=M_1(\widetilde{QC})\cup M_{-1}(\widetilde{QC})\cup \bigcup_{\tau\in\T_+} M_\tau(\widetilde{QC}),
\\
 M_{\pm1}(\widetilde{QC}) &= \left( M_{\pm1}(\widetilde{QC}) \setminus M_{\pm1}^0(\widetilde{QC}) \right) \cup M_{\pm1}^0(\widetilde{QC}) ,
\\
\label{dec-3}
  M_\tau(\widetilde{QC}) &=  \left(M_\tau^+(\widetilde{QC}) \setminus  M_\tau^0(\widetilde{QC})\right) \cup\left(  M^-_\tau(\widetilde{QC}) \setminus  M^0_\tau(\widetilde{QC})\right) \cup  M_\tau^0(\widetilde{QC}) ,\qquad \tau\in\T_+.
\end{align}
Furthermore, for each $\eta\in M(\widetilde{QC})$ we have either
$$
|M^\eta(QC)|=1 \quad\mbox{ or }\quad |M^\eta(QC)|=2.
$$
The first case happens if and only if $\eta\in M_1^0(\widetilde{QC})\cup M_{-1}^0(\widetilde{QC})$.
Then $M^\eta(QC)=\{\xi\}$ with $\xi=\xi'\in M_{\pm1}^0(QC)$. 
If the second case occurs, $M^\eta(QC)=\{\xi,\xi'\}$ with $\xi\neq \xi'$.

\bigskip
\bigskip
 In order to study the local algebras $\mathcal{F}^J_\eta(PQC_{\mathcal{L}(H)})$, we have to consider two cases separately. Namely, when $\eta\in M_\tau(\widetilde{QC})\setminus M_\tau^0(\widetilde{QC})$ and when $\eta\in M_\tau^0(\widetilde{QC})$. It turns out that, in the first case, for the purpose of establishing 
 the stability result, it is redundant to identity $\mathcal{F}_\eta^J(PQC_{\mathcal{L}(H)})$. For that reason, we will just give the invertibility criterion for the first case without identifying the local structures. The second case is more complicated, and we will analyze the structure of the local algebras in a constructive way. 
 Note that in each of these cases we have to further distinguish whether $\tau = \pm1$ or $\tau\in\T_+$.
 

\section{Invertibility in local algebras for $\eta\in M_\tau(\widetilde{QC})\setminus M_\tau^0(\widetilde{QC})$}
\label{s:7}

We start with the case where $\tau\in \T_+$.

\begin{thm}\label{thm61} Let $\eta\in M_\tau(\widetilde{QC})\setminus M_\tau^0(\widetilde{QC})$ with $\tau\in \T_+.$ If $(A_n)\in \mathcal{F}^J(PQC_{\mathcal{L}(H)})$ and $\mathcal{P}(A_n)$ is invertible in $\mathcal{S}^J(PQC_{\mathcal{L}(H)})$, then $(A_n)+\mathcal{J}^J_{\eta, \mathcal{L}(H)}$ is invertible in $\mathcal{F}^J_\eta(PQC_{\mathcal{L}(H)})$.
\end{thm}
\begin{proof}
\enskip We first define two *-homomorphisms between $\mathcal{F}^J_\eta(PQC_{\mathcal{L}(H)})$ and $\mathcal{S}^J_\eta(PQC_{\mathcal{L}(H)})$. Let $\mathcal{P}'$ and $\Phi'$ be the mappings

        \begin{align}
        \mathcal{P}'&: \mathcal{F}^J_\eta(PQC_{\mathcal{L}(H)}) \to \mathcal{S}^J_\eta(PQC_{\mathcal{L}(H)}), \enskip (A_n)+ \mathcal{J}^J_{\eta, \mathcal{L}(H)}\mapsto \mathcal{P}(A_n) + \mathcal{I}^J_{\eta, \mathcal{L}(H)},\\
        \Phi'&: \mathcal{S}^J_\eta(PQC_{\mathcal{L}(H)}) \to \mathcal{F}^J_\eta(PQC_{\mathcal{L}(H)}), \enskip A+ \mathcal{I}^J_{\eta, \mathcal{L}(H)}\mapsto (P_nAP_n) + \mathcal{J}^J_{\eta, \mathcal{L}(H)}.
        \end{align}
By Proposition \ref{prop36}(c), $\mathcal{P}$ is a *-homomorphism from $\mathcal{F}^J(PQC_{\mathcal{L}(H)})$ onto $\mathcal{S}^J(PQC_{\mathcal{L}(H)})$. Since $\mathcal{P}$ sends all the generating elements of $\mathcal{J}^J_{\eta, \mathcal{L}(H)}$ into $\mathcal{I}^J_{\eta, \mathcal{L}(H)}$, the mapping $\mathcal{P}'$ is a well-defined *-homomorphism. By Lemma \ref{lem45}, $\Phi'$ is correctly defined, linear and symmetric from $\mathcal{S}^J_\eta(PQC_{\mathcal{L}(H)})$ into $\mathcal{F}^J_\eta(PQC_{\mathcal{L}(H)})$.

It remains to show $\Phi'$ is actually multiplicative. For this purpose, we claim that for each $g\in \widetilde{PC}_{\mathcal{L}(H)}$, the class of all $\cL(H)$-valued even piecewise continuous functions, there exists an $a\in \mathcal{L}(H)$ such that $g-a\in \mathcal{I}^J_{\eta, \mathcal{L}(H)}$. Without loss of generality, suppose $\eta \in M_\tau(\widetilde{QC})\setminus M_\tau^+(\widetilde{QC})$. Then, one can find an $f\in \widetilde{QC}$ such that $\eta(f)\neq 0$ and $\limsup\limits_{t\to \tau+0} |f(t)| = 0$. Indeed, by Proposition \ref{prop59} and Proposition \ref{prop510}, there exists a unique $\xi\in M^-_\tau(QC)\setminus M_\tau^+(QC)$ such that $\hat{\xi} = \eta$, and thus there is a function $f'\in QC$ with $\xi(f')\neq 0$ and $\limsup\limits_{t\to \tau+0} |f'(t)| = 0$ by definition. Further, choose a smooth function $c_\tau$ on $\mathbb{T}$ such that $c_\tau=1$ around $\tau$ and vanishes on the lower semi-circle, and define $f:=f'c_\tau + \widetilde{f'c_\tau}$. Then, $f\in QC$ is even, $f-f'$ is continuous at $\tau$ and vanishes there. By an approximation argument, we have $\eta(f) = \xi(f) = \xi(f')\neq0$, hence $f$ has the desired property.

Now, choose $a\in \mathcal{L}(H)$ such that $\limsup\limits_{t\to \tau-0} \|g(t)-a\| = 0$. Thus $f(g-a)$ is continuous at $\pm \tau$ and vanishes there. Since $\eta \in M_\tau(\widetilde{QC})$, by an approximation argument one can conclude that $f(g-a)\in \mathcal{I}^J_{\eta, \mathcal{L}(H)}$. On the other hand, we have $f - \eta(f)e \in \mathcal{I}^J_{\eta, \mathcal{L}(H)}$ and therefore $(f-\eta(f)e)(g-a)\in \mathcal{I}^J_{\eta, \mathcal{L}(H)}$. Since $\eta(f)\neq 0$, we obtain that $g-a\in \mathcal{I}^J_{\eta, \mathcal{L}(H)}$. The case where $\eta \in M_\tau(\widetilde{QC})\setminus M_\tau^-(\widetilde{QC})$ can be treated analogously.

Choose an odd smooth function $b_\tau\in C(\mathbb{T})$ such that $b_\tau = 1$ around $\tau$ and $b_\tau = -1$ around $1/\tau$. Every $p\in PC_{\cL(H)}$ admits a unique decomposition $p = p_e+ p_o$ where $p_e$ is even and $p_o$ is odd. Moreover, it can be decomposed as
        \[p = p_e + (p_ob_\tau)\cdot b_\tau + p_o(1-b_\tau^2).\]
Since $1-b_\tau^2$ is continuous at $\tau$ and $1/\tau$ and vanishes there, $p_o(1-b_\tau^2)\in \mathcal{I}_{\eta,\mathcal{L}(H)}^J$. Further, $p_e, p_o  b_\tau\in \widetilde{PC}_{\cL(H)} \subseteq \mathcal{L}(H) + \mathcal{I}_{\eta,\mathcal{L}(H)}^J$ by the previous argument. In general, for any $p\in PC_{\mathcal{L}(H)}$, by \eqref{PCLH}, it can be approximated arbitrarily close by a finite sum $\sum a_ip^i = \sum a_i(p_e^i + p_o^i)$ where $a_i\in \mathcal{L}(H), p^i\in PC$. Denote by $\mathcal{S}^J(QC_{\mathcal{L}(H)}^s)$ (resp. $\mathcal{S}^J(\widetilde{QC}_{\mathcal{L}(H)}^s)$) the smallest closed algebra containing all operators $f\in QC_{\mathcal{L}(H)}^s$ (resp. $\widetilde{QC}_{\mathcal{L}(H)}^s$) and the operators $P, Q$ and $J$. Based on the above argument, we have
        \[\mathcal{S}^J(PQC_{\mathcal{L}(H)}) = \mathrm{clos}\{A+B: \enskip A\in\mathcal{S}^J(QC_{\mathcal{L}(H)}^s), B\in  \mathcal{I}^J_{\eta, \mathcal{L}(H)}\}.\]

Similarly, for any $q\in QC_{\mathcal{L}(H)}^s$,
        \[q = q_e + q_o = q_e + (q_ob_\tau)\cdot b_\tau + q_o(1-b_\tau^2).\]
Again, by an approximation argument,
        \[\mathcal{S}^J(PQC_{\mathcal{L}(H)}) = \mathrm{clos}\{A+B\cdot b_\tau+C:\enskip A,B\in \mathcal{S}^J(\widetilde{QC}_{\mathcal{L}(H)}^s), C\in \mathcal{I}^J_{\eta, \mathcal{L}(H)}\}.\]
Since $b_\tau\in QC$ is odd, it commutes with every element in $\mathcal{S}(PQC_{\mathcal{L}(H)})$ modulo $\mathcal{K}$ by Lemma \ref{lem42}, and $Jb_\tau = -b_\tau J$. Further
        \begin{align*}
        \mathcal{S}_\eta^J(\widetilde{QC}_{\mathcal{L}(H)}^s) &= \mathrm{clos} \{(q_1P + q_2 Q + q_3 PJ + q_4 QJ) + \mathcal{I}^J_{\eta, \mathcal{L}(H)}: q_i\in \widetilde{QC}_{\mathcal{L}(H)}^s\},\\
        \mathcal{S}^J_\eta(PQC_{\mathcal{L}(H)}) &= \mathrm{clos} \{(q_1P + q_2 Q + q_3 PJ + q_4 QJ) +\\
        &\qquad     (q_1'P + q_2' Q + q_3' PJ + q_4' QJ)b_\tau + \mathcal{I}^J_{\eta, \mathcal{L}(H)}: q_i, q_i'\in \widetilde{QC}_{\mathcal{L}(H)}^s\}.
        \end{align*}
Using Proposition \ref{prop22} and Lemma \ref{lem42}, from the above presentation, we obtain by straightforward computations that $\Phi'$ is multiplicative, and thus we have shown that it is a $^*$-homomorphism.

We claim that the $^*$-homomorphism $\Phi'\circ \mathcal{P}' : \mathcal{F}^J_\eta(PQC_{\mathcal{L}(H)}) \to \mathcal{F}^J_\eta(PQC_{\mathcal{L}(H)})$ is the identity mapping. It suffices to show that $\Phi'\circ \mathcal{P}'$ maps all the generating elements of $\mathcal{F}^J_\eta(PQC_{\mathcal{L}(H)})$ to themselves. Indeed, by Proposition \ref{prop36}(b), $\mathcal{P}(P_nAP_n) = A$ for any $A\in \mathcal{S}^J(PQC_{\mathcal{L}(H)})$, and the assertion simply follows.

To finish the proof, let $(A_n)\in \mathcal{F}^J(PQC_{\mathcal{L}(H)})$, and assume that $\mathcal{P}(A_n)$ is invertible in $\mathcal{S}^J(PQC_{\mathcal{L}(H)})$. Then $\mathcal{P}(A_n) + \mathcal{I}^J_{\eta, \mathcal{L}(H)}$ is invertible in $\mathcal{S}^J_\eta(PQC_{\mathcal{L}(H)})$. On the other hand, by definition, $\mathcal{P}(A_n) + \mathcal{I}^J_{\eta, \mathcal{L}(H)} = \mathcal{P}'((A_n) + \mathcal{J}^J_{\eta, \mathcal{L}(H)})$. Hence
                $$(\Phi'\circ \mathcal{P}')((A_n) + \mathcal{J}^J_{\eta, \mathcal{L}(H)}) = (A_n) + \mathcal{J}^J_{\eta, \mathcal{L}(H)}$$
is invertible in $\mathcal{F}^J_\eta(PQC_{\mathcal{L}(H)})$.
\end{proof}

Next we consider the case where $\eta\in M_\tau(\widetilde{QC})\setminus M_\tau^0(\widetilde{QC})$ with $\tau=\pm1$. By 
Theorem \ref{thm56}(b), 
there exists a unique $\xi \in M_\tau^+(QC)\setminus M_\tau^-(QC)$, such that $\hat{\xi}' = \hat{\xi}' = \eta$. Again, we construct two *-homomorphisms between $\mathcal{F}^J_\eta(PQC_{\mathcal{L}(H)})$ and $\mathcal{S}^J_\eta(PQC_{\mathcal{L}(H)})$ to show that they are inverse to each other. 

\begin{thm}\label{thm62} Let $\eta\in M_\tau(\widetilde{QC})\setminus M_\tau^0(\widetilde{QC})$ with $\tau=\pm1$. If $(A_n)\in \mathcal{F}^J(PQC_{\mathcal{L}(H)})$ and $\mathcal{P}(A_n)$ is invertible in $\mathcal{S}^J(PQC_{\mathcal{L}(H)})$, then $(A_n)+\mathcal{J}^J_{\eta, \mathcal{L}(H)}$ is invertible in $\mathcal{F}^J_\eta(PQC_{\mathcal{L}(H)})$.
\end{thm}
\begin{proof}
For sake of definiteness, let $\tau=1$. The case $\tau=-1$ is analogous.
Similar to Theorem \ref{thm61}, let $\mathcal{P}'$ and $\Phi'$ be the mappings
        \begin{align*}
        \mathcal{P}'&:\enskip \mathcal{F}^J_\eta(PQC_{\mathcal{L}(H)}) \to \mathcal{S}^J_\eta(PQC_{\mathcal{L}(H)}), \enskip (A_n)+ \mathcal{J}^J_{\eta, \mathcal{L}(H)}\mapsto \mathcal{P}(A_n) + \mathcal{I}^J_{\eta, \mathcal{L}(H)},\\
        \Phi'&:\enskip \mathcal{S}^J_\eta(PQC_{\mathcal{L}(H)}) \to \mathcal{F}^J_\eta(PQC_{\mathcal{L}(H)}), \enskip A+ \mathcal{I}^J_{\eta, \mathcal{L}(H)}\mapsto (P_nAP_n) + \mathcal{J}^J_{\eta, \mathcal{L}(H)}
        \end{align*}
as defined before. Again, it remains to show that $\Phi'$ is multiplicative.

With $\xi \in M_1^+(QC)\setminus M_1^-(QC)$ chosen above, there exists an $f\in QC$, such that $\xi(f) = 1$, and $f|M_1^-(QC) = 0$. Define $\overline{f} := f+ \tilde{f}$. Then
        \[\xi(\overline{f}) = \xi(f) + \xi'(f) = \xi(f) = 1 = \eta(\overline{f}), \quad \overline{f}|M_1^0(QC) = 0.\]
Choose an odd function $c\in PC\cap C(\mathbb{T}\setminus\{1\})$ such that
        \[\quad c(e^{ix}) =\begin{cases}
        1   & \text{if}\enskip 0<x<\pi/2 \\
        -1  & \text{if}\enskip -\pi/2<x<0.
        \end{cases}\]
We have
        \[c = c\overline{f} + c(1-\overline{f}) = c\overline{f} \enskip \mathrm{mod} \enskip \mathcal{I}_{\eta,\mathcal{L}(H)}^J.\]
Note that $c\overline{f}\in QC$ by Proposition \ref{QCprop}(c). For any $p\in PC$, there exist constants $\alpha, \beta\in \mathbb{C}$ such that $p- \alpha c-\beta$ is continuous at 1 and vanishes there, thus
        \[p = \alpha c+\beta+ (p-\alpha c-\beta) = \alpha c\overline{f} + \beta \enskip \mathrm{mod} \enskip \mathcal{I}_{\eta,\mathcal{L}(H)}^J.\]
Since every constant can be regarded as a function in $\widetilde{QC}$, it shows that $PC \subseteq \{g+ h\cdot c\overline{f}+ C: g,h\in \widetilde{QC}, C\in \mathcal{I}_{\eta,\mathcal{L}(H)}^J\}$.

Similarly, for any function $q\in QC$,
        \[q = q_e + q_o = q_e + (q_oc)c\overline{f} + q_o(1-c^2\overline{f}).\]
We have $q_oc\in QC$ by Proposition \ref{QCprop}(c), and $QC \subseteq \{g+ h\cdot c\overline{f}+ C: g,h \in \widetilde{QC}, C\in \mathcal{I}_{\eta,\mathcal{L}(H)}^J\}$.
Since every element in $PC_{\mathcal{L}(H)}$ (resp., $QC_{\mathcal{L}(H)}^s$) can be approximated arbitrarily close by a finite sum $\sum a_if_i$, where $a_i\in \mathcal{L}(H), f_i\in PC$ (resp., $f_i\in QC$), an approximation argument shows that
        \[\mathcal{S}^J(PQC_{\mathcal{L}(H)}) = \mathrm{clos} \{A+B\cdot c\overline{f} + C: A, B\in \mathcal{S}^J(\widetilde{QC}_{\mathcal{L}(H)}^s), \enskip C\in \mathcal{I}_{\eta,\mathcal{L}(H)}^J\}.\]
Since $c\overline{f}\in QC$ is odd, it commutes with all the generating elements in $\mathcal{S}(\widetilde{QC}_{\mathcal{L}(H)}^s)$ and $c\overline{f}J = -Jc\overline{f}$. An argument similar to the proof of Theorem \ref{thm61} shows that $\Phi'$ is
multiplicative, and $\Phi', \enskip \mathcal{P}'$ are inverse to each other. Therefore, the invertibility of $(A_n)+\mathcal{J}^J_{\eta, \mathcal{L}(H)}$ in $\mathcal{F}^J_\eta(PQC_{\mathcal{L}(H)})$  follows analogously  that $\mathcal{P}(A_n)$ is invertible in $\mathcal{S}^J(PQC_{\mathcal{L}(H)})$.
\end{proof}

\section{Identification of local algebras ($\tau\in \T_+$)}

In the sequel, we analyze the local algebras $\mathcal{F}^J_\eta(PQC_{\mathcal{L}(H)})$ in the case where $\eta \in M_\tau^0(\widetilde{QC})$ for $\tau\in \T_+$.  To identify these algebras, we will eliminate the flip by doubling up the dimension. In fact, by \cite[Scheme~3.3]{Roch1988SymbolCalwithFlip}, one has the following lemma:

\begin{lem}\label{lem71} Let $\mathcal{X}$ be a $C^*$-algebra with identity e whose center contains a self-adjoint projection $p$. Let $\mathcal{Y}$ be generated by $\mathcal{X}$ and a self-adjoint flip j with the properties that $j\mathcal{X}j \subseteq \mathcal{X}$ and, in particular, $jpj = e-p$. Then any element of $\mathcal{Y}$ can be written uniquely as a sum $a_1+a_2j$ with $a_1, a_2\in \mathcal{X}$, and the mapping $L : \mathcal{Y}\to [p\mathcal{X}p]^{2\times 2}$ which maps $a = a_1+a_2 j$ to
        \[ \left(
        \begin{array}{cc}
        pa_1p & pa_2p\\
        p\tilde{a}_2p & p\tilde{a}_1p\end{array}
        \right),\]
where  $\tilde{a}:= jaj$, is an isometric *-isomorphism.
\end{lem}

Denote by 
\begin{align*}
\Phi_{\eta, \mathcal{S}} &: 
\mathcal{S}^J(PQC_{\mathcal{L}(H)}) \to 
\mathcal{S}_\eta^J(PQC_{\mathcal{L}(H)}),
\\
 \Phi_{\eta, \mathcal{F}} &: 
 \mathcal{F}^J(PQC_{\mathcal{L}(H)})\to \mathcal{F}_\eta^J(PQC_{\mathcal{L}(H)})
\end{align*}
the quotient maps, and put
\begin{equation}
        \mathcal{X}_\mathcal{S} = \Phi_{\eta, \mathcal{S}}(\mathcal{S}(PQC_{\mathcal{L}(H)})),  \label{X_SF}
        \qquad\quad
           \mathcal{X}_\mathcal{F} = \Phi_{\eta, \mathcal{F}}(\mathcal{F}(PQC_{\mathcal{L}(H)})),
 \end{equation}
which are *-subalgebras of $\mathcal{S}_\eta^J(PQC_{\mathcal{L}(H)})$ and $\mathcal{F}^J_\eta(PQC_{\mathcal{L}(H)})$, respectively.
Let $c_\tau$ be a continuous function which equals to $1$ around $\tau$ and vanishes on the lower semi-circle, and define
\begin{alignat*}{3}
        p_\cS &= \Phi_{\eta, \mathcal{S}}(L(c_\tau)),   &\qquad  
        j_\cS &= \Phi_{\eta, \mathcal{S}}(J),                &\qquad 
        e_\cS &= \Phi_{\eta, \mathcal{S}}(I),\\
        p_\cF &= \Phi_{\eta, \mathcal{F}}(P_n L(c_\tau)P_n),  & 
        j_\cF  &= \Phi_{\eta, \mathcal{F}}(P_nJP_n),      &
         e_\cF &= \Phi_{\eta, \mathcal{F}}(P_n).
\end{alignat*}
We can then apply the preovious lemma in the following two settings,
$$
p=p_\cS,\quad j=j_\cS, \quad e=e_\cS,\quad \mathcal{X}=\mathcal{X}_\cS, \quad\mathcal{Y}=\cS_\eta^J(PQC_{\mathcal{L}(H)}),
$$
and 
$$
p=p_\cF,\quad j=j_\cF,\quad e=e_\cF,\quad \mathcal{X}=\mathcal{X}_\cF,\quad \mathcal{Y}=\cF_\eta^J(PQC_{\mathcal{L}(H)}).
$$
Indeed, $p$ is a projection which commutes with all the generating elements of $\mathcal{X}$ by Proposition \ref{prop22}, and $jpj=e-p$. Furthermore, 
$\cS_\eta^J(PQC_{\mathcal{L}(H)})$ is generated by $j_\cS$ and $\mathcal{X}_\cS=\Phi_{\eta,\cS}(\cS(PQC_{\mathcal{L}(H)}))$
on the one hand, and $\cF_\eta^J(PQC_{\mathcal{L}(H)})$ is generated by $j_\cF$ and $\mathcal{X}_\cF=\Phi_{\eta,\cF}(\cS(PQC_{\mathcal{L}(H)}))$ on the other hand.

As a consequence, each $a\in \mathcal{S}_\eta^J(PQC_{\mathcal{L}(H)})$ can be uniquely written as $a = a_1 + a_2j_\cS$ with $a_1,a_2\in \mathcal{X}_\mathcal{S}$, and each $a\in \mathcal{F}_\eta^J(PQC_{\mathcal{L}(H)})$ can be uniquely written as $a = a_1 + a_2j_\cF$ with $a_1,a_2\in \mathcal{X}_\mathcal{F}$. Furthermore we have the following isometric $^*$-isomorphisms 
\begin{equation}\label{eq:8.2}
L_{\eta,\mathcal{S}}: \cS_\eta^J(PQC_{\mathcal{L}(H)}) \to [p_\cS\mathcal{X}_\mathcal{S}p_\cS]^{2\times 2},
\qquad
L_{\eta,\mathcal{F}}: \cF_\eta^J(PQC_{\mathcal{L}(H)}) \to [p_\cF\mathcal{X}_\mathcal{F}p_\cF]^{2\times 2}.
\end{equation}
In view of our goal of identifying the local algebras
$\cS_\eta^J(PQC_{\mathcal{L}(H)})$ and $\cF_\eta^J(PQC_{\mathcal{L}(H)})$, we have reduced the problem to identifying the 
$C^*$-algebras $p_\cS \mathcal{X}_\cS p_\cS$ and $p_\cF \mathcal{X}_\cF p_\cF$, respectively.  Note that whenever $\mathcal{X}$ is unital  $C^*$-algebra with  a self-adjoint projection $p$, $p\mathcal{X} p$ is a unital $C^*$-algebra with unit element $p$.

Denote by $\mathcal{S}(PC_{\mathcal{L}(H)})$ the smallest closed subalgebra of $\mathcal{L}(l_H^2(\mathbb{Z}))$ which contains all Laurent operators $L(f)$ with $f\in PC_{\mathcal{L}(H)}$, the operators $P$ and $Q$, and the ideal $\mathcal{K}$. Define $\mathcal{F}(PC_{\mathcal{L}(H)})$ as the smallest closed subalgebra of $\mathcal{F}$ which includes the ideal $\mathcal{J}$ and all sequences $(P_nAP_n)$ with $A\in \mathcal{S}(PC_{\mathcal{L}(H)})$. Then, $\mathcal{S}(PC_{\mathcal{L}(H)})$ is a *-subalgebra of $\mathcal{S}(PQC_{\mathcal{L}(H)})$ containing $\mathcal{K}$, and $\mathcal{F}(PC_{\mathcal{L}(H)})$ is a *-subalgebra of $\mathcal{F}(PQC_{\mathcal{L}(H)})$ including $\mathcal{J}$.

\begin{lem}\label{lem72} Let $\eta\in M_\tau^0(\widetilde{QC})$ with $\tau\neq \pm1$. Then
        \begin{equation}\label{f.pX_Sp}    
        p_\cS\mathcal{X}_\mathcal{S}p_\cS = p_\cS\Phi_{\eta, \mathcal{S}}(\mathcal{S}(PC_{\mathcal{L}(H)}))p_\cS    
        \end{equation}
   and
         \begin{equation}\label{f.pX_Fp} 
        p_\cF\mathcal{X}_\mathcal{F}p_\cF = p_\cF\Phi_{\eta, \mathcal{F}}(\mathcal{F}(PC_{\mathcal{L}(H)}))p_\cF.   
        \end{equation}
\end{lem}
\begin{proof} 
The inclusions ``$\supseteq$'' are trivial. Note that 
$p_\cS\Phi_{\eta, \mathcal{S}}(\mathcal{S}(PC_{\mathcal{L}(H)}))p_\cS\label{pX_Sp}$ is a (closed) *-subalgebra of $\cS_\eta^J(PQC_{\mathcal{L}(H)})$ with (different) unit element $p_\cS$, and $p_\cF\Phi_{\eta, \mathcal{F}}(\mathcal{F}(PC_{\mathcal{L}(H)}))p_\cF$ is a (closed) *-subalgebra of $\cF_\eta^J(PQC_{\mathcal{L}(H)})$ with (different) unit element $p_\cF$.

Let us focus on the first identity. In order to show the inclusion ``$\subseteq$'', it suffices to show that 
each element of the form $p_\mathcal{S} a_1\cdots a_N p_\cS$ belongs to the right hand side where $a_1,\dots,a_N$
are generating elements of $\mathcal{X}_\cS$. Since $p_\cS$ commutes with every $a_k$, it is enough
to show that $p_\mathcal{S}  a p_\mathcal{S}$ belongs to the right hand side for every generating element $a$.
Note that this is trivially the case if $a\in \Phi_{\eta, \mathcal{S}}(\mathcal{S}(PC_{\mathcal{L}(H)}))$. Therefore we are left with showing that 
$$
p_\cS \Phi_{\eta,\cS}(f) p_\cS\in p_\cS\Phi_{\eta,\cS}(\cS(PC_{\cL(H)})) p_\cS
$$
for $f\in QC^s_{\cL(H)}$. It suffices to consider $f=q a$ with $a\in \LH$ and $q\in QC$.

Let $q\in QC$ be given. By Proposition \red{\ref{prop59}} and Proposition \red{\ref{prop510}}, there is a unique $\xi \in M_\tau^0(QC)$ such that $\hat{\xi} = \eta$. Define $\overline{q} := c_\tau q + \widetilde{c_\tau q}-\xi(q)$ where $c_\tau$ was introduced above. Then $\overline{q}\in \widetilde{QC}$  and 
	\[\eta(\overline{q})=\xi(\overline{q})= \xi( c_\tau q + \widetilde{c_\tau q})-\xi(q)=0. \]
 Furthermore,
	\[c_\tau(q-\xi(q)-\overline{q})=c_\tau((1-c_\tau) q - \widetilde{c_\tau q})\]
is continuous at $\tau$ and $1/\tau$ and vanishes there, and  it belongs to the ideal in $QC$ generated by 
$\eta$ by an approximation argument. It follows that 
	\[c_\tau (q-\xi(q)) \in \mathrm{clos} \, \mathrm{id}_{QC_{\LH}^s}\{\: f\,:\, f\in \widetilde{QC},\; \eta(f)=0\;\}\]
which implies $L(c_\tau (q-\xi(q))  \in \mathcal{I}_{\eta,\LH}^J$, and $p_\cS \Phi_{\eta,\cS}(q) p_\cS= \xi(q)\cdot p_\cS$ for scalar $q\in QC$.

In general, for $f\in QC^s_{\LH}$, 
\begin{equation}
   p_\cS \Phi_{\eta,\cS}(f) p_\cS = \Phi_\xi(f) \cdot p_{\cS} \subseteq \Phi_{\eta,\cS}(\cS(PC_\LH))
\end{equation}   
followed by an approximation argument. Here $\Phi_\xi(f)\in\LH$ is given by \eqref{locallytrivial}, i.e., 
we are using the fact that $f$ is locally trivial at $\xi$. This completes the proof of \eqref{f.pX_Sp}.

For the inclusion ``$\subseteq$'' of \eqref{f.pX_Fp}, we argue as in the first lines of above
and observe that it is sufficient to prove that 
$$
p_\cF \Phi_{\eta,\cF}(P_n A P_n )  p_\cF\in  p_\cF \Phi_{\eta,\cF}( PC_\LH) p_\cF
$$
for $A\in \cS(PQC_\LH)$, noting that $(P_nA P_n)$ are the generating elements of $\cF(PQC_\LH)$.

Thus let $A\in \cS(PQC_\LH)$ be given. Using the identity \eqref{f.pX_Sp}, it follows that there exists
$B\in \cS(PC_\LH)$ such that $c_\tau(A-B)c_\tau\in \mathcal{I}^J_{\eta,\cS}$. Using Lemma \ref{lem42}(b) and 
Lemma \ref{lem45} it follows that 
$$
(P_nc_\tau P_n (A-B) P_n c_\tau P_n)\cong  (P_n c_\tau(A-B) c_\tau P_n)\cong  0  \mod \cJ_{\eta,\LH}^J.
$$
Hence $p_\cF \Phi_{\eta,\cF}(P_nA P_n)p_\cF = p_\cF \Phi_{\eta,\cF}(P_n B P_n) p_\cF$ with the latter belonging to 
$p_\cF \Phi_{\eta,\cF}(\cF(PC_\LH))p_\cF$.
\end{proof}

Now we are going to make connections between the $C^*$-algebras $p_\cS \mathcal{X}_\cS p_\cS$ and $p_\cF \mathcal{X}_\cF p_\cF$ 
on the one hand, and between $C^*$-algebras $\cS_1(PC_\LH)$ and $\cF_1(PC_\LH)$ on the other hand. The latter $C^*$-algebras already occurred in 
\cite{ehrhardt1996finite} and are being introduced next. Define the following $*$-ideals
        \begin{align}
        \mathcal{I}_{1, \mathcal{L}(H)}&:=\mathrm{clos}\, \mathrm{id}_{\mathcal{S}(PC_{\mathcal{L}(H)})}(\mathcal{K}\cup \{f:f\in C(\mathbb{T}), f(1)=0\}),\\
        \mathcal{J}_{1, \mathcal{L}(H)}&:=\mathrm{clos}\, \mathrm{id}_{\mathcal{F}(PC_{\mathcal{L}(H)})}(\mathcal{J}\cup \{(P_nfP_n):f\in C(\mathbb{T}), f(1)=0\}),
        \end{align}
and let 
\begin{equation}
\mathcal{S}_1(PC_{\mathcal{L}(H)})
  := \mathcal{S}(PC_{\mathcal{L}(H)})/ \mathcal{I}_{1, \mathcal{L}(H)},
  \quad
\mathcal{F}_1(PC_{\mathcal{L}(H)})
   :=  \mathcal{F}(PC_{\mathcal{L}(H)})/ \mathcal{J}_{1, \mathcal{L}(H)}
\end{equation}
be the corresponding quotient $C^*$-algebras. 
Furthermore, for $\tau\in \mathbb{T}$, define the operators $Y_\tau$ by
        \begin{equation}
        Y_\tau:  l_H^2(\mathbb{Z}) \to l_H^2(\mathbb{Z}),\qquad (x_n)_{n=-\infty}^\infty \mapsto (\tau^{-n}x_n)_{n=-\infty}^\infty.
        \end{equation}
Evidently, $Y_\tau^{-1} = Y_{\tau^{-1}} = Y_\tau^*$.

\begin{prop} \label{prop74}Let  $\eta \in M_\tau^0(\widetilde{QC})$ with $\tau\in \T_+$ .
 Then the mappings
        \begin{alignat}{2}
        \label{PhitauS}& \Phi_{\tau, \mathcal{S}}: \mathcal{S}_1(PC_{\mathcal{L}(H)}) \to p_\cS\mathcal{X}_\mathcal{S}p_\cS,
        \enskip\quad&
         A+ \mathcal{I}_{1, \mathcal{L}(H)} & \mapsto p_\cS\Phi_{\eta, \mathcal{S}}(Y_\tau AY_\tau^{-1})p_\cS,\\
        \label{PhitauF}& \Phi_{\tau, \mathcal{F}}: \mathcal{F}_1(PC_{\mathcal{L}(H)}) \to p_\cF\mathcal{X}_\mathcal{F}p_\cF,
         \enskip&
         (A_n)+ \mathcal{J}_{1, \mathcal{L}(H)} &\mapsto p_\cF\Phi_{\eta, \mathcal{F}}(Y_\tau A_nY_\tau^{-1})p_\cF
        \end{alignat}
are well-defined, surjective *-homomorphisms.
\end{prop}
\begin{proof} 
Note that the algebras $\mathcal{S}(PC_{\mathcal{L}(H)})$ and $\mathcal{F}(PC_{\mathcal{L}(H)})$ as well as the ideal  $\mathcal{K}$ and $\mathcal{J}$ are rotation invariant. It can be easily verified since $(Y_\tau P_nY_\tau^{-1}) = (P_n),\enskip Y_\tau P Y_\tau^{-1} = P,\enskip Y_\tau L(a)Y_\tau^{-1} = L(a_\tau)$, where $a_\tau(t) := a(t/\tau), t \in \mathbb{T}$, for the generating elements of the algebras mentioned above. In particular,  $Y_\tau (P_n A P_n) Y_{\tau}^{-1} = ( P_n Y_\tau A Y_\tau^{-1} P_n)$.

Define the mappings 
        \begin{alignat*}{2}
        & \Phi^0_{\tau, \mathcal{S}}: \mathcal{S}(PC_{\mathcal{L}(H)}) \to p_\cS\mathcal{X}_\mathcal{S}p_\cS,
        \enskip\quad&
         A & \mapsto p_\cS\Phi_{\eta, \mathcal{S}}(Y_\tau AY_\tau^{-1})p_\cS,\\
        & \Phi^0_{\tau, \mathcal{F}}: \mathcal{F}(PC_{\mathcal{L}(H)}) \to p_\cF\mathcal{X}_\mathcal{F}p_\cF,
         \enskip&
         (A_n) &\mapsto p_\cF\Phi_{\eta, \mathcal{F}}(Y_\tau A_nY_\tau^{-1})p_\cF.
        \end{alignat*}
These are *-homomorphisms, for which the multiplicativity follows from the fact that $p_\cS$ and $p_\cF$ are 
commuting with all elements of $\mathcal{X}_\cS$ and $\mathcal{X}_\cF$, respectively. To show that $\Phi_{\tau,\cS}$ and $\Phi_{\tau,\cF}$ are well-defined, consider $f\in C(\mathbb{T})$ with $f(1) = 0$.  Then the first map sends $A=L(f)$ into $p_\cS \Phi_{\eta,\mathcal{S}}(f_\tau) =  \Phi_{\eta,\mathcal{S}}(c_\tau f_\tau)=0$ since $c_\tau f_\tau$ is continuous at $\tau$ and $1/\tau$ and vanishes there. Indeed, $f_\tau(\tau)=f(1)=0$ and $c_\tau(1/\tau)=0$ by  the choice of $c_\tau$.
Therefore, the generating elements of $\mathcal{I}_{1,\LH}$ are mapped to zero. Using rotation invariance, it follows that all of $\mathcal{I}_{1,\LH}$ are mapped to zero by \eqref{PhitauS}, and thus
$\Phi_{\tau,\cS}$ is a well-defined *-homomorphism. Along the same lines it can be shown that $\cJ_{1,\LH}$ is sent to zero by \eqref{PhitauF}, which implies that the *-homorphism $\Phi_{\tau,\cF}$ is well-defined as well.

The statement that $\Phi_{\tau,\cS}$ and $\Phi_{\tau,\cF}$ are surjective follows immediately from Lemma  \ref{lem72}, together with 
the rotation invariance of $\mathcal{S}(PC_{\mathcal{L}(H)})$ and $\mathcal{F}(PC_{\mathcal{L}(H)})$.
\end{proof}

The next step is to show that the `local algebras'  $\cS_1(PC_\LH)$ and $\cF_1(PC_\LH)$ are *-isomorphic to certain 
$C^*$-algebras of operators $\Sigma_\LH$ and $\Xi_\LH$, respectively. These are algebras of  operators acting on the space $L^2_H(\R)$ and $L^2_H([-1,1])$, respectively, i.e., the spaces of square-integrable $H$-valued functions 
on the real line and on the interval $[-1,1]$, respectively. 
This identification has already been done in 
\cite{ehrhardt1996finite}.

Let $\Sigma_{\mathcal{L}(H)}$ denote the smallest closed subalgebra of $\mathcal{L}(L_H^2(\mathbb{R}))$ which contains the operator 
$\chi_{[0,\infty)}I$ of multiplication  by the characteristic function of the positive half-axis, considered as operator  $f(x)\in L^2_H(\R) \mapsto \chi_{[0,\infty)}(x) f(x) \in L^2_H(\R)$, the singular integral operator $S_\mathbb{R}$ on the real line, 
$$
(S_\R f)(x) = \frac{1}{\pi i}\int_\R \frac{f(y)}{x-y}\, dy,
$$
and all constants $a\in \mathcal{L}(H)$, considered as operators $f(x)\in L^2_H(\R) \mapsto a f(x) \in L^2_H(\R)$.
In brief notation, 
        \begin{equation}
        \Sigma_{\mathcal{L}(H)} := \mathrm{alg}\{\chi_{[0,\infty)}I, S_\mathbb{R}, a\in \mathcal{L}(H)\}.\label{SigmaLH}
        \end{equation}
Furthermore, denote by $\Xi_{\mathcal{L}(H)}$ the smallest closed subalgebra of $\mathcal{L}(L_H^2([-1,1]))$ containing all operators $\chi_{[-1,1]}A\chi_{[-1,1]}$ with $A\in \Sigma_{\mathcal{L}(H)}$, i.e., 
        \begin{equation}
        \Xi_{\mathcal{L}(H)} :=\text{alg}\{\chi_{[-1,1]}A\chi_{[-1,1]}:\enskip A\in \Sigma_{\mathcal{L}(H)}\}.\label{XiLH}
        \end{equation}
Here $\chi_{[-1,1]}$ stands for the operator of multiplitation with the characteristic function. Note that       
$\Sigma_{\mathcal{L}(H)}$ and $\Xi_{\mathcal{L}(H)}$ are $C^*$-algebras.

Let $E_n$ and $E_{-n}$ stand for the bounded linear operators given by ($n\ge 1$)
        \begin{align}
        E_n&: l_H^2(\mathbb{Z})\to L_H^2(\mathbb{R}),\enskip (x_i)_{i=-\infty}^\infty \mapsto \sqrt{n}\sum_{i=-\infty}^\infty x_i\chi_{[\frac{i}{n}, \frac{i+1}{n}]},\label{EN}\\
        E_{-n}&: L_H^2(\mathbb{R})\to l_H^2(\mathbb{Z}),\enskip f \mapsto \left(\sqrt{n}\int_{-\infty}^\infty f(x)\chi_{[\frac{i}{n}, \frac{i+1}{n}]}(x)dx\right)_{i=-\infty}^\infty.\label{E-N}
        \end{align}
Evidently, $E^*_{-n} = E_n$ and $E_{-n}E_n = I$. 

\begin{prop}\label{prop75}\hfill
\begin{enumerate}[label=(\alph*)]
\item There exist two *-homomorphisms $\mathcal{E}_{\mathcal{S}, \mathcal{L}(H)}$ and $\mathcal{E}_{\mathcal{F}, \mathcal{L}(H)}$
        \begin{align}\label{eq:8.16}
        \mathcal{E}_{\mathcal{S}, \mathcal{L}(H)}&: \mathcal{S}(PC_{\mathcal{L}(H)}) \to \Sigma_{\mathcal{L}(H)} \quad \mathrm{with}\quad \mathcal{I}_{1, \mathcal{L}(H)} \subseteq \ker \mathcal{E}_{\mathcal{S}, \mathcal{L}(H)},\\
        \label{eq:8.17}
        \mathcal{E}_{\mathcal{F}, \mathcal{L}(H)}&: \mathcal{F}(PC_{\mathcal{L}(H)}) \to \Xi_{\mathcal{L}(H)} \quad \mathrm{with}\quad \mathcal{J}_{1, \mathcal{L}(H)} \subseteq \ker \mathcal{E}_{\mathcal{F}, \mathcal{L}(H)}
        \end{align}
which map onto the corresponding algebras, and for which
        \begin{align} 
        \label{eq:8.18}
        \mathcal{E}_{\mathcal{S}, \mathcal{L}(H)}(P) &= \chi_{[0,\infty)}I ,\\
        \mathcal{E}_{\mathcal{S}, \mathcal{L}(H)}(a) &= a(1+0)\frac{I-S_\mathbb{R}}{2}+ a(1-0)\frac{I+S_\mathbb{R}}{2},\enskip (a\in PC_{\mathcal{L}(H)})\\
        \mathcal{E}_{\mathcal{F}, \mathcal{L}(H)}(P_nAP_n) &= \chi_{[-1,1]}\mathcal{E}_{\mathcal{S}, \mathcal{L}(H)}(A)\chi_{[-1,1]}. \qquad\qquad (A\in \mathcal{S}(PC_{\mathcal{L}(H)})) \label{eq:8.20}
        \end{align}
\item Let $H = \mathbb{C}$. Then  $\mathcal{E}_{\mathcal{S}, \mathcal{L}(H)}$ and $\mathcal{E}_{\mathcal{F}, \mathcal{L}(H)}$ are given 
by the following strong limits: 
        \begin{alignat}{2}\label{eq:8.21}
        \mathcal{E}_\mathcal{S}(A) &:= \slim_{n\to \infty} E_n A E_{-n},   &\qquad\qquad
        A &\in \mathcal{S}(PC),\\  \label{eq:8.22}
        \mathcal{E}_\mathcal{F}(A_n) &:= \slim_{n\to \infty} E_n A_n E_{-n} &
        (A_n) &\in  \mathcal{F}(PC).
         \end{alignat}
\item
$\ker \mathcal{E}_{\mathcal{S}, \mathcal{L}(H)} = \mathcal{I}_{1, \mathcal{L}(H)}$ and $\ker \mathcal{E}_{\mathcal{F}, \mathcal{L}(H)} = \mathcal{J}_{1, \mathcal{L}(H)}$.
\item
The induced mappings
     \begin{alignat}{2}
        \tilde{\mathcal{E}}_{\mathcal{S}, \mathcal{L}(H)} & : \mathcal{S}_1(PC_{\mathcal{L}(H)})\to \Sigma_{\mathcal{L}(H)},
        &\qquad\quad
       A+\mathcal{I}_{1,\LH} &\mapsto \mathcal{E}_\cS(A) ,   \label{PCiso1}
 \\
       \tilde{\mathcal{E}}_{\mathcal{F}, \mathcal{L}(H)} &: \mathcal{F}_1(PC_{\mathcal{L}(H)}) \to \Xi_{\mathcal{L}(H)}
       &
         (A_n)+\mathcal{J}_{1,\LH} &\mapsto \mathcal{E}_\cF(A) 
       \label{PCiso2}
        \end{alignat}
are  well-defined $^*$-isomorphisms. The inverses are given by
\begin{alignat}{2}
     \tilde{\mathcal{E}}_{\mathcal{S}, \mathcal{L}(H)}^{-1} & :  \Sigma_{\mathcal{L}(H)} \to \mathcal{S}_1(PC_{\mathcal{L}(H)}),
        &\qquad\quad
       B &\mapsto E_{-1} B E_{1}+\mathcal{I}_{1,\LH}  ,   \label{PCiso1x}
 \\
       \tilde{\mathcal{E}}_{\mathcal{F}, \mathcal{L}(H)}^{-1} &:  \Xi_{\mathcal{L}(H)} \to \mathcal{F}_1(PC_{\mathcal{L}(H)}) 
       &
      B &\mapsto   (E_{-n} B  E_n)+\mathcal{J}_{1,\LH}
       \label{PCiso2x}
        \end{alignat}
\end{enumerate}
\end{prop}
\begin{proof}
Parts (a) and (b) are stated and proved in Proposition 7.2 of  \cite{ehrhardt1996finite}.
The fact that   $\tilde{\mathcal{E}}_{\mathcal{S}}$  and $\tilde{\mathcal{E}}_{\mathcal{F}}$ are well-defined follows from the 
inclusions of the kernels stated in \eqref{eq:8.16} and \eqref{eq:8.17}. The definitions of $\Sigma_\LH$ and $\Xi_\LH$ together with 
\eqref{eq:8.18}--\eqref{eq:8.20} imply that the *-homomophisms are surjective. The statement that $\tilde{\mathcal{E}}_{\mathcal{S}}$  and $\tilde{\mathcal{E}}_{\mathcal{F}}$
are *-isomorphisms will thus follow from (c). 

Statement (c) is proved in Proposition 9.7 of \cite{ehrhardt1996finite} for the scalar case $H=\mathbb{C}$.
For the general Hilbert space case, it follows immediately by a general tensorization argument for short exact sequences, which is stated in 
Proposition 9.8 of \cite{ehrhardt1996finite} (see also Corollaries 9.9 and 9.10 therein).
The formulas for the inverses in (d) also follow from the just-mentioned results.
\end{proof}

As an aside to the previous proposition let us point out that we will not use part (b). Its purpose is rather to serve as a construction of the
*-homomorphisms $\mathcal{E}_\cS$ and $\mathcal{E}_\cF$ in the scalar case. In the general case, they can be obtained
via a tensorization argument (as done in \cite{ehrhardt1996finite}). Alternatively, one could also define them via \eqref{eq:8.21} and \eqref{eq:8.22}  as `generalized' strong limits
in the same manner as it has been done in Section 3. Notice that the underlying space is $L^2_H(\R)$ rather than $\ell^2_H(\Z)$.
We will not pursue the details at this point.

\medskip

Let us summarize what kinds of *-homomorphism and *-isomorphims  between the various
$C^*$-algebras we have constructed so far.
We can illustrate this in the following diagram:
\begin{center}
\begin{tikzpicture}[thick,scale=1]
 \matrix (m) [matrix of math nodes,row sep=4em,column sep=6em,minimum width=3em]
  {
     p_\cS\mathcal{X}_\mathcal{S}p_\cS & \mathcal{S}_1(PC_{\mathcal{L}(H)})&\Sigma_{\mathcal{L}(H)}\\
     p_\cF\mathcal{X}_\mathcal{F}p_\cF & \mathcal{F}_1(PC_{\mathcal{L}(H)})&\Xi_{\mathcal{L}(H)}\\};
  \path[-stealth]
        (m-1-2) edge node [above] {\footnotesize{$\Phi_{\tau,\mathcal{S}}$}} (m-1-1)
        (m-1-2) edge node [above] {\footnotesize{$\widetilde{\mathcal{E}}_{\mathcal{S},\mathcal{L}(H)}$}} (m-1-3)
        (m-2-1) edge node [left]  {\footnotesize{$\mathcal{P}_\eta$}} (m-1-1)
        (m-2-2) edge node [above] {\footnotesize{$\Phi_{\tau,\mathcal{F}}$}} (m-2-1)
        (m-2-2) edge node [left]  {\footnotesize{$\mathcal{P}_1$}} (m-1-2)
        (m-2-2) edge node [above] {\footnotesize{$\widetilde{\mathcal{E}}_{\mathcal{F},\mathcal{L}(H)}$}} (m-2-3)
        (m-2-3) edge node [left]  {\footnotesize{$\Phi_Z$}} (m-1-3);
\end{tikzpicture}
\end{center}
By the previous proposition, we already know that  $\widetilde{\mathcal{E}}_{\mathcal{S},\mathcal{L}(H)}$ and $\widetilde{\mathcal{E}}_{\mathcal{S},\mathcal{L}(H)}$ are $^*$-isomorphism, and our next goal is to show that also $\Phi_{\tau,\mathcal{S}}$ and $\Phi_{\tau,\mathcal{F}}$
are $^*$-isomorphisms. 

The `vertical' mappings were mentioned only for completeness sake:
\begin{alignat*}{1}
 	\mathcal{P}_\eta &: (A_n) + \mathcal{J}_{\eta,\LH} \mapsto \mathcal{P}(A_n) +  \mathcal{I}_{\eta,\LH}, 
	\\
	\mathcal{P}_1 &: (A_n) + \mathcal{J}_{1,\LH} \mapsto \mathcal{P}(A_n) +  \mathcal{I}_{1,\LH},
\end{alignat*}	
while $\Phi_Z$ is defined in (87) of \cite{ehrhardt1996finite}. These mappings are all *-homomorphisms.
We will not make use of them.

\bigskip

In view of Proposition \ref{prop74} we need to prove that ther kernels of $\Phi_{\tau, \mathcal{S}}$ and $\Phi_{\tau, \mathcal{F}}$
are trivial in order to show that these mappings are *-isomorphism. We follow essentially the method used in 
\cite{ehrhardt1996finite} with necessary modification. We need a lemma and the following definition.

For a function $f\in L^\infty$, let $\sigma_nf$ denote the \emph{Fejer-Cesaro} mean
                $$(\sigma_nf)(e^{ix}) := \sum_{k=-n}^n \Big(1-\frac{|k|}{n+1}\Big)f_ke^{ikx}, \enskip x\in [0,2\pi).$$

\begin{lem}\label{lem76} Let $\tau\in \mathbb{T}$ and $\xi\in M_\tau^0(QC)$. Then, for all $p, q\in QC$, there is a sequence $\{k_n\}_{n=1}^\infty \subseteq \mathbb{Z}^+$ with $k_n\to\infty$ such that $(\sigma_{2k_n-1}q)(\tau)\to \xi(q)$ and $(\sigma_{2k_n-1}p)(\tau)\to \xi(p)$ as $n\to \infty$.
\end{lem}
\begin{proof} For the proof, see Lemma 9.11 in \cite{ehrhardt1996finite}. Note that here we have two quasicontinuous functions $p, q$.
Therefore, one has to choose the corresponding neighborhood $U_\epsilon$ of $\xi$ in $QC^*$ as
                $$U_\varepsilon = \{\rho\in QC^*: |\rho(p)- \xi(p)|+ |\rho(q) - \xi(q)| + |\rho(\chi_1) - \xi(\chi_1)|<\varepsilon\},$$
where $\chi_1$ is the function $\chi_1(t) = t, \enskip t\in \mathbb{T}$.
\end{proof}

\begin{prop}\label{prop77}
The kernels of $\Phi_{\tau, \mathcal{S}}$ and $\Phi_{\tau, \mathcal{F}}$ are trivial.
\end{prop}
\begin{proof}
First we consider $\Phi_{\tau, \mathcal{F}}$. Assume that $\ker \Phi_{\tau, \mathcal{F}}\neq\{0\}$. By  \eqref{PCiso2}, $\mathcal{F}_1(PC_{\mathcal{L}}(H))$ and $\Xi_{\mathcal{L}(H)}$ are $^*$-isomorphic via $\tilde{\mathcal{E}}_{\mathcal{F}, \mathcal{L}(H)}$, and the ideal $\ker \Phi_{\tau, \mathcal{F}}$ corresponds to a non-trivial ideal $\mathfrak{J}$ of $\Xi_{\mathcal{L}(H)}$. For $x,y\in L_H^2([-1,1])$, define the operator $K_{x,y}\in \mathcal{K}(L_H^2([-1,1]))$ by
        \begin{equation}\label{Kxy}
        K_{x,y}z := x\cdot\langle y,z\rangle_{L_H^2([-1,1])}, \enskip z\in L_H^2([-1,1]).
        \end{equation}
Note that $K_{x,y}\in \mathcal{K}(L^2([-1,1]))\otimes \mathcal{K}(H)\subseteq \Xi\otimes\mathcal{L}(H) \cong \Xi_{\mathcal{L}(H)}$.

Now take $A\in \mathfrak{J}\setminus\{0\}$, and choose $x_1,x_2\in L^2_H([-1,1])$ such that $x_2=A x_1\neq0$. Put
$x = \chi_{[-1,1]}\in L^2([-1,1])$ and choose any $h\in H, \enskip \|h\| = 1$. Then
	\begin{equation*}
		K_{xh,xh} \cdot \langle x_2, A x_1\rangle = K_{xh,x_2} A K_{x_1,xh}\in \mathfrak{J}, 
	\end{equation*}
and hence $K_{xh,xh}\in \mathfrak{J}$.  
From \eqref{PCiso2x} it follows that 
        \[\tilde{\mathcal{E}}^{-1}_{\mathcal{F}, \mathcal{L}(H)}(K_{xh,xh}) = (E_{-n}K_{xh,xh}E_n) + \mathcal{J}_{1,\mathcal{L}(H)}\in \ker \Phi_{\tau, \mathcal{F}},\]
and by applying the homomorphism $ \Phi_{\tau, \mathcal{F}}$ we get
        \begin{equation}
        (K_n) := (P_nc_\tau P_nY_\tau E_{-n}K_{xh,xh}E_nY_{\tau}^{-1}P_n c_\tau P_n)\in \mathcal{J}_{\eta, \mathcal{L}(H)}^J.
        \end{equation}
We are going to prove that this leads to a contradiction.        
        
Let $0<\varepsilon<1/5$. By the definition of $\mathcal{J}_{\eta, \mathcal{L}(H)}^J$, there is a sequence $(A_n)$ of the form
	\begin{equation*}
		(A_n) = \sum_{i=1}^k (A_n^{(i)})(P_nf_iP_n) + (B_n'),
	\end{equation*}
where $(A_n^{(i)})\in \mathcal{F}^J({PQC_{\mathcal{L}(H)}}), \enskip f_i\in \widetilde{QC},\enskip \eta(f_i) = 0$ and $(B_n')\in \mathcal{J}$, such that 
        \[\|(A_n) - (K_n)\|_\mathcal{F}\le\epsilon.\]
Consider the open neighborhood $U$ of $\eta$ in $M(\widetilde{QC})$,
        \[U := \{\zeta\in M(\widetilde{QC}): |\zeta(f_i)|<\epsilon k^{-1}\|(A_n^{(i)})\|_\mathcal{F}^{-1}, \enskip \forall i\; \}.\]
By the Gelfand-Naimark Theorem, there exists an $f\in \widetilde{QC}$ for which $\eta(f) = 1,\enskip\zeta(f)\in [0,1]$ if $\zeta\in U$, and $\zeta(f) = 0$ if $\zeta\notin U$. Hence $\|f\| = 1$, and $\sum_{i=1}^k \|(A_n^{(i)})\|\cdot \|ff_i\|<\epsilon$. Therefore, 
multiplying the above equation with $(P_nf P_n)$ from the left, we get with
 some $(B_n'')\in \mathcal{J}$ (see Lemma \ref{lem42}(b))
        \begin{gather*}
        (A_n)(P_nfP_n) = \sum_{i=1}^k (A_n^{(i)})(P_nf_ifP_n) + (B_n''),\\
        \|(A_n)(P_nfP_n) - (B_n'')\|_\mathcal{F} \le \epsilon.
        \end{gather*}
Put $(K_n') := (P_nc_\tau P_nY_\tau E_{-n}K_{xh,xh}E_nY_{\tau}^{-1}P_n c_\tau f P_n)$. Then
        \begin{align*}
        (K_n') &= (K_nP_nfP_n) - (B_n''')\\
        &= (A_nP_nfP_n) - ((A_n-K_n)P_nfP_n) - (B_n''')\\
        &= (A_nP_nfP_n - B_n'') - ((A_n-K_n)P_nfP_n) + (B_n),
        \end{align*}
where $(B_n) = (B_n'') - (B_n''') = (P_nKP_n+W_nLW_n + C_n')\in \mathcal{J}$ for some $K, L\in \mathcal{K}$ and $(C_n')\in \mathcal{N}$. 
It follows that 
$$
\|(K_n')-(B_n)\|_{\cF}\le 2\varepsilon.
$$
Now, choose a sufficiently large $\kappa$ such that $\|Q_{\kappa}K\|_\mathcal{A}\le \epsilon$ and $\|Q_{\kappa}L\|_\mathcal{A}\le \epsilon$, and define
        \begin{equation}
        (R_n) := (P_nQ_\kappa P_n)(W_n Q_\kappa W_n) = (W_n Q_\kappa W_n)(P_nQ_\kappa P_n).
        \end{equation}
We have
	\begin{align*}
		R_n(B_n-C_n') &= R_n (P_nK P_n +W_nL W_n)\\
				   &= (W_n Q_\kappa W_n)  (P_n Q_\kappa  K P_n)  +(P_n Q_\kappa P_n)(W_n  Q_\kappa  L W_n),
	\end{align*}
which implies that 
$$
\| (R_n)(B_n-C_n')\|_{\cF}\le 2 \varepsilon.
$$
Observing that $\|(R_n)\| = 1$ and estimating
$$
\|R_n K_n'\| \le \|R_n (K_n'-B_n)\|+\|R_n(B_n- C_n')\|+\|R_n C_n' \|,
$$
we obtain
        \begin{equation}
        \limsup_{n\to \infty}\|R_nK_n'\|_{\mathcal{L}(l_H^2(\mathbb{Z}_n))}\le 4\epsilon.\label{contra}
        \end{equation}
Next let $z_n$ and $z_n^*$ be the bounded linear operators ($n\ge 1$)
        \begin{align}
        \label{zzz}z_n&:\enskip \mathbb{C}\to l_H^2(\mathbb{Z}_n), \quad \lambda\mapsto \left(\lambda h/\sqrt{2n}\right)_{i=-n}^{n-1},\\
        \label{zzzz}z_{-n}&:\enskip l_H^2(\mathbb{Z}_n) \to\mathbb{C} , \quad (x_i)_{i=-n}^{n-1}\mapsto \sum\limits_{i=-n}^{n-1}\langle h,x_i\rangle_H/\sqrt{2n}.
        \end{align}
Then $\|z_n\| = \|z_n^*\| = 1$ and 
     $$
        (E_{-n}K_{xh,xh}E_n ) = ( z_nz_n^*).
     $$
 Therefore, for $n\ge 2\kappa$, we have
        \begin{align*}
        \|R_nK_n'\| &\ge |z_n^*Y_\tau^{-1}R_nK_n'Y_\tau z_n|\\
        &=     |z_n^*Y_\tau^{-1}R_n    P_n c_\tau P_n Y_\tau E_{-n} K_{xh,xh} E_n Y_{\tau}^{-1} P_n c_\tau f  P_n Y_\tau z_n|\\
        &= |z_n^*Y_\tau^{-1}R_nP_nc_\tau P_n Y_\tau z_n|\cdot |z_n^*Y_\tau^{-1}P_nc_\tau fP_nY_\tau z_n|.
        \end{align*}
Denote the $k$-th Fourier coefficients of $c_\tau$ by $(c_\tau)_k$. Observe that
        \begin{align*}
        |z_n^*Y_\tau^{-1}(P_n-R_n)P_nc_\tau P_n Y_\tau z_n| 
        &\le\frac{1}{2n}\left|\sum_{l\in \Lambda} \sum_{k=-n}^{n-1}(c_\tau)_{l-k}\tau^{l-k}\right|\\
        &\le\frac{4\kappa}{2n}\sum_{k=-\infty}^\infty|(c_\tau)_k|\\
        &= \frac{4\kappa C}{2n} \to 0
        \end{align*}
as $n\to\infty$, where $C = \sum\limits_{k=-\infty}^\infty |(c_\tau)_k| <\infty$, and
        \[\Lambda = \{l\in \mathbb{Z}: -n\le l\le -n+\kappa -1 \text{ or } -\kappa\le l\le \kappa-1 \text{ or } n-\kappa \le l \le n-1\}.\]
Therefore,
        \begin{align*}
        |z_n^*Y_\tau^{-1}R_n P_n c_\tau P_n Y_\tau P_n|
        & \ge|z_n^*Y_\tau^{-1}P_nc_\tau P_nY_\tau z_n| - |z_n^*Y_\tau^{-1}(P_n-R_n)P_nc_\tau P_n Y_\tau z_n|\\
        & \ge  \left|\sum_{i=-2n+1}^{2n-1} (c_\tau)_i \tau^i \frac{2n-|i|}{2n}\right| - \frac{4\kappa C}{2n}\\
        & =  |(\sigma_{2n-1}c_\tau)(\tau)| - \frac{4\kappa C}{2n}.
        \end{align*}
Analogously, 
   $$        
        |z_n^*Y_\tau^{-1}P_nc_\tau fP_nY_\tau z_n| = |(\sigma_{2n-1}( c_\tau f))(\tau)|.
    $$    
For $\eta\in M^0_\tau(\widetilde{QC})$, there exists a unique $\xi\in M_\tau^0(QC)$ such that $\hat{\xi} = \eta$. Since $\xi(c_\tau) = \xi(c_\tau f) = 1$, we obtain from Lemma \ref{lem76} and \eqref{contra} that
        \begin{eqnarray*}
       4\varepsilon \ge  \limsup_{n\to\infty}\|R_nK_n'\| &\ge& \limsup_{n\to\infty} \left(|(\sigma_{2n-1}c_\tau)(\tau)| - \frac{4\kappa C}{2n}\right)\cdot
         |(\sigma_{2n-1}(c_\tau f))(\tau)|\\
        &\ge& |\xi(c_\tau)\cdot \xi(c_\tau f)| = 1,
        \end{eqnarray*}
which contradicts the choice of $\varepsilon$. This completes the proof that $\ker \Phi_{\tau, \mathcal{F}}$ is trivial.

Now we consider the kernel of $\Phi_{\tau, \mathcal{S}}$. Let $A+\mathcal{I}_{1,\mathcal{L}(H)}\in \ker\Phi_{\tau, \mathcal{S}}$ 
where $A\in \mathcal{S}(PC_{\mathcal{L}(H)})$. In view of Proposition \ref{prop74} we have
$p_{\cS}\Phi_{\eta,\cS}(Y_\tau AY_\tau^{-1})p_{\cS}=
c_\tau Y_\tau A Y_\tau^{-1} c_\tau\in \mathcal{I}_{\eta,\mathcal{L}(H)}^J$, 
 and thus $(P_nc_\tau  Y_\tau A Y_\tau^{-1}c_\tau P_n)\in \mathcal{J}_{\eta,\mathcal{L}(H)}^J$
 by  Lemma \ref{lem45}. Furthermore, by
 Lemma \ref{lem42}  we get 
 $$
  (P_nc_\tau P_n Y_\tau P_nA P_nY_\tau^{-1} P_nc_\tau P_n)\in \mathcal{J}_{\eta,\mathcal{L}(H)}^J
  $$
   Since $(P_nAP_n)\in \mathcal{F}(PC_{\mathcal{L}(H)})$, we have (by the definition of $\Phi_{\tau,\cF}$)
        \[(P_nAP_n) + \mathcal{J}_{1,\mathcal{L}(H)}\in \ker \Phi_{\tau, \mathcal{F}},\]
which implies $(P_nAP_n)\in \mathcal{J}_{1,\mathcal{L}(H)}$ by what we just proved. Now we apply the homomorphism $\mathcal{P}$ to this sequence (see \eqref{P}). As $\mathcal{P}$ maps $\mathcal{J}_{1,\mathcal{L}(H)}$ into $\mathcal{I}_{1,\mathcal{L}(H)}$, we obtain that $\mathcal{P}(P_nAP_n) = A\in \mathcal{I}_{1,\mathcal{L}(H)}$. Thus  we have shown that  also the kernel of $\Phi_{\tau, \mathcal{S}}$ is trivial.
\end{proof}

Now, consider the induced mappings $\Phi'_{\tau, \mathcal{S}}, \Phi'_{\tau, \mathcal{F}}, \widetilde{\mathcal{E}}'_{\mathcal{S},\mathcal{L}(H)}$ and $\widetilde{\mathcal{E}}'_{\mathcal{F},\mathcal{L}(H)}$ as defined below
        \begin{align}
        \Phi'_{\tau, \mathcal{S}}: \mathcal{S}_1(PC_{\mathcal{L}(H)})^{2\times 2} \to& [p\mathcal{X}_\mathcal{S}p]^{2\times 2},\enskip
        \left( \begin{array}{cc}
        A & B\\
        C & D\\
        \end{array} \right)
        \mapsto
        \left( \begin{array}{cc}
        \Phi_{\tau, \mathcal{S}}(A) & \Phi_{\tau, \mathcal{S}}(B)\\
        \Phi_{\tau, \mathcal{S}}(C) & \Phi_{\tau, \mathcal{S}}(D)\\
        \end{array} \right),\\
        \Phi'_{\tau, \mathcal{F}}: \mathcal{F}_1(PC_{\mathcal{L}(H)})^{2\times 2} \to& [p'\mathcal{X}_\mathcal{F}p']^{2\times 2}, \enskip
        \left( \begin{array}{cc}
        A & B\\
        C & D\\
        \end{array} \right)
        \mapsto
        \left( \begin{array}{cc}
        \Phi_{\tau, \mathcal{F}}(A) & \Phi_{\tau, \mathcal{F}}(B)\\
        \Phi_{\tau, \mathcal{F}}(C) & \Phi_{\tau, \mathcal{F}}(D)\\
        \end{array} \right),\\
        \widetilde{\mathcal{E}}'_{\mathcal{S},\mathcal{L}(H)}: \mathcal{S}_1(PC_{\mathcal{L}(H)})^{2\times 2} \to& \Sigma_{\mathcal{L}(H)}^{2\times 2},
        \left( \begin{array}{cc}
        A & B\\
        C & D\\
        \end{array} \right)
        \mapsto
        \left( \begin{array}{cc}
        \widetilde{\mathcal{E}}_{\mathcal{S},\mathcal{L}(H)}(A) & \widetilde{\mathcal{E}}_{\mathcal{S},\mathcal{L}(H)}(B)\\
        \widetilde{\mathcal{E}}_{\mathcal{S},\mathcal{L}(H)}(C) & \widetilde{\mathcal{E}}_{\mathcal{S},\mathcal{L}(H)}(D)\\
        \end{array} \right),\\
        \widetilde{\mathcal{E}}'_{\mathcal{F},\mathcal{L}(H)}: \mathcal{F}_1(PC_{\mathcal{L}(H)})^{2\times 2} \to& \Xi_{\mathcal{L}(H)}^{2\times 2},
        \left( \begin{array}{cc}
        A & B\\
        C & D\\
        \end{array} \right)
        \mapsto
        \left( \begin{array}{cc}
        \widetilde{\mathcal{E}}_{\mathcal{F},\mathcal{L}(H)}(A) & \widetilde{\mathcal{E}}_{\mathcal{F},\mathcal{L}(H)}(B)\\
        \widetilde{\mathcal{E}}_{\mathcal{F},\mathcal{L}(H)}(C) & \widetilde{\mathcal{E}}_{\mathcal{F},\mathcal{L}(H)}(D)\\
        \end{array} \right).
        \end{align}
Evidently, they are $^*$-isomorphisms by Proposition \ref{prop77} and Proposition \ref{prop75}(c). Further, define the $^*$-isomorphisms
        \begin{equation}
        \Psi'_{\eta,\mathcal{S}} =\widetilde{\mathcal{E}}'_{\mathcal{S},\mathcal{L}(H)} \circ(\Phi'_{\tau, \mathcal{S}})^{-1}\circ L_{\eta,\mathcal{S}},
        \enskip
        \Psi'_{\eta,\mathcal{F}} =\widetilde{\mathcal{E}}'_{\mathcal{F},\mathcal{L}(H)} \circ(\Phi'_{\tau, \mathcal{F}})^{-1}\circ L_{\eta,\mathcal{F}}.\label{Psi notpm1}
        \end{equation}

\begin{cor}\label{cor78} For $\eta\in M_\tau^0(\widetilde{QC})$ with $\tau\in \T_+$,\hfill
\begin{enumerate}[label=(\alph*)]
\item $\mathcal{S}^J_\eta(PQC_{\mathcal{L}(H)})$ is $^*$-isomorphic to $\Sigma_{\mathcal{L}(H)}^{2\times 2}$ via $\Psi'_{\eta,\mathcal{S}}$, and\\
\item $\mathcal{F}^J_\eta(PQC_{\mathcal{L}(H)})$ is $^*$-isomorphic to $\Xi_{\mathcal{L}(H)}^{2\times 2}$ via $\Psi'_{\eta,\mathcal{F}}$.
\end{enumerate}
\end{cor}
\section{Identification of local algebras ($\tau = \pm 1$)}
\label{s:9}

Finally, we consider $S^J_\eta(PQC_{\mathcal{L}(H)})$ where $\eta\in M_1^0(\widetilde{QC})$, as the case where $\eta\in M_{-1}^0(\widetilde{QC})$ follows in a similar way. Note that, without the flip operator, the local structures of $\mathcal{S}_\xi(PQC_{\mathcal{L}(H)})$ and $\mathcal{F}_\xi(PQC_{\mathcal{L}(H)})$ with $\xi\in M_\tau^0(QC)$ have been studied explicitly in the previous paper.

Let $\mathcal{S}^J(PC_{\mathcal{L}(H)})$ be the smallest closed subalgebra of $\mathcal{L}(l_H^2(\mathbb{Z}))$ which contains all Laurent operators $L(f)$ with $f\in PC_{\mathcal{L}(H)}$ and the operators $P, Q$ and $J$. Further, denote by $\mathcal{F}^J(PC_{\mathcal{L}(H)})$ the smallest closed subalgebra of $\mathcal{F}$ containing the ideal $\mathcal{J}$ and all sequences $(P_nAP_n)$ with $A\in \mathcal{S}^J(PC_{\mathcal{L}(H)})$. Then, $\mathcal{S}^J(PC_{\mathcal{L}(H)})$ is a $^*$-subalgebra of $\mathcal{S}^J(PQC_{\mathcal{L}(H)})$ containing $\mathcal{K}$, and $\mathcal{F}^J(PC_{\mathcal{L}(H)})$ is a $^*$-subalgebra of $\mathcal{F}^J(PQC_{\mathcal{L}(H)})$ including $\mathcal{J}$. Similar to $\mathcal{I}_{1,\mathcal{L}(H)}$ and $\mathcal{J}_{1,\mathcal{L}(H)}$, we define the $^*$- ideals

        \begin{align*}
        \mathcal{I}^J_{1,\mathcal{L}(H)} :=& \mathrm{clos}\enskip \mathrm{id}_{\mathcal{S}^J(PC_{\mathcal{L}(H)})}\left( \mathcal{K}\cup \{f:f\in \widetilde{C}(\mathbb{T}), f(1) = 0\}\right)\\
        =& \mathrm{clos}\enskip \mathrm{id}_{\mathcal{S}^J(PC_{\mathcal{L}(H)})}\left( \mathcal{K}\cup \{f:f\in C(\mathbb{T}), f(1) = 0\}\right),\\
        \mathcal{J}^J_{1,\mathcal{L}(H)} :=& \mathrm{clos}\enskip \mathrm{id}_{\mathcal{F}^J(PC_{\mathcal{L}(H)})}\left( \mathcal{J}\cup \{(P_nfP_n):f\in \widetilde{C}(\mathbb{T}), f(1) = 0\}\right)\\
        =& \mathrm{clos}\enskip \mathrm{id}_{\mathcal{F}^J(PC_{\mathcal{L}(H)})}\left( \mathcal{J}\cup \{(P_nfP_n):f\in C(\mathbb{T}), f(1) = 0\}\right).
        \end{align*}
We denote the quotient algebras $\mathcal{S}^J(PC_{\mathcal{L}(H)})/\mathcal{I}^J_{1,\mathcal{L}(H)}$ and $\mathcal{F}^J(PC_{\mathcal{L}(H)})/\mathcal{J}^J_{1,\mathcal{L}(H)}$ by $\mathcal{S}_1^J(PC_{\mathcal{L}(H)})$ and $\mathcal{F}_1^J(PC_{\mathcal{L}(H)})$, respectively.
Define $c = \chi_+ - \chi_-\in PC$. Note that, $f-\eta(f) \in \mathcal{I}^J_{\eta, \mathcal{L}(H)}$ whenever $f\in \widetilde{QC}$, and for any $q\in QC$,
                $$q = q_e + (q_o\cdot c)c = \eta(q_e)+ \eta(q_o\cdot c)c+ C$$
with some $C\in \mathcal{I}^J_{\eta, \mathcal{L}(H)}$ by Proposition \ref{QCprop}. Hence, by an approximation argument,
        \begin{align}
        \label{decomp 91} \mathcal{S}^J(PQC_{\mathcal{L}(H)}) &= \text{clos} \{A + B: A\in \mathcal{S}^J(PC_{\mathcal{L}(H)}), B\in \mathcal{I}^J_{\eta, \mathcal{L}(H)}\},\\
        \label{decomp 92} \mathcal{F}^J(PQC_{\mathcal{L}(H)}) &= \text{clos} \{(A_n) + (B_n): (A_n)\in \mathcal{F}^J(PC_{\mathcal{L}(H)}), (B_n)\in \mathcal{J}^J_{\eta, \mathcal{L}(H)}\}
        \end{align}
for $\eta\in M_1^0(\widetilde{QC})$. Thus we obtain the following, which is the counter-part to Proposition \ref{prop74}.
\begin{prop} \label{prop81}Let $\eta\in M_1^0(\widetilde{QC})$. Then the mappings
        \begin{align}
        & \Phi_{1, \mathcal{S}}: \mathcal{S}_1^J(PC_{\mathcal{L}(H)}) \to \mathcal{S}_\eta^J(PQC_{\mathcal{L}(H)}), A+ \mathcal{I}^J_{1, \mathcal{L}(H)} \mapsto A + \mathcal{I}^J_{\eta,\mathcal{L}(H)},\\
        & \Phi_{1, \mathcal{F}}: \mathcal{F}_1^J(PC_{\mathcal{L}(H)}) \to \mathcal{F}_\eta^J(PQC_{\mathcal{L}(H)}), (A_n)+ \mathcal{J}^J_{1, \mathcal{L}(H)} \mapsto (A_n) + \mathcal{J}^J_{\eta,\mathcal{L}(H)}\qquad\qquad
        \end{align}
are surjective $^*$-homomorphisms.
\end{prop}

Similar to \eqref{SigmaLH}, let $\Sigma'_{\mathcal{L}(H)}$ denote the smallest closed subalgebra of $\mathcal{L}(L_H^2(\mathbb{R}))$ which contains the multiplication operator $\chi_{[0,\infty)}$, the singular integral operator $S_\mathbb{R}$, all constants $a\in \mathcal{L}(H)$ and the (continuous) flip operator $\hat{J}$, i.e.,
        \begin{equation}
        \Sigma'_{\mathcal{L}(H)} := \text{alg}\{\chi_{[0,\infty)}, S_\mathbb{R}, \hat{J}, a\in \mathcal{L}(H)\}.\label{SigmaLH'}
        \end{equation}
We remark that $\hat{J}$ is defined by $(\hat{J}f)(x) = f(-x)$.

Analogously, denote by $\Xi'_{\mathcal{L}(H)}$ the smallest closed subalgebra of $\mathcal{L}(L_H^2([-1,1]))$ containing all operators $\chi_{[-1,1]}A\chi_{[-1,1]}$ with $A\in \Sigma'_{\mathcal{L}(H)}$, i.e.,
        \begin{equation}
        \Xi'_{\mathcal{L}(H)} :=\text{alg}\{\chi_{[-1,1]}A\chi_{[-1,1]}:\enskip A\in \Sigma'_{\mathcal{L}(H)}\}.\label{XiLH'}
        \end{equation}
Evidently, $\Sigma'_{\mathcal{L}(H)}$ and $\Xi'_{\mathcal{L}(H)}$ are both $C^*$-algebras. 

\begin{prop}\label{prop82}\hfill
\begin{enumerate}[label=(\alph*)]
\item There exist two $^*$-homomorphisms $\mathcal{E}^J_{\mathcal{S}, \mathcal{L}(H)}$ and $\mathcal{E}^J_{\mathcal{F}, \mathcal{L}(H)}$
        \begin{align}
        \mathcal{E}^J_{\mathcal{S}, \mathcal{L}(H)}&: \mathcal{S}^J(PC_{\mathcal{L}(H)}) \to \Sigma'_{\mathcal{L}(H)} \quad \mathrm{with}\quad \mathcal{I}^J_{1, \mathcal{L}(H)} \subseteq \ker \mathcal{E}^J_{\mathcal{S}, \mathcal{L}(H)},\\
        \mathcal{E}^J_{\mathcal{F}, \mathcal{L}(H)}&: \mathcal{F}^J(PC_{\mathcal{L}(H)}) \to \Xi'_{\mathcal{L}(H)} \quad \mathrm{with}\quad \mathcal{J}^J_{1, \mathcal{L}(H)} \subseteq \ker \mathcal{E}^J_{\mathcal{F}, \mathcal{L}(H)}
        \end{align}
which map onto the corresponding algebras, and for which
        \begin{align*}
        \mathcal{E}^J_{\mathcal{S}, \mathcal{L}(H)}(P) &= \chi_{[0,\infty)},\\
        \mathcal{E}^J_{\mathcal{S}, \mathcal{L}(H)}(a) &= a(1+0)\frac{I-S_\mathbb{R}}{2}+ a(1-0)\frac{I+S_\mathbb{R}}{2}, \enskip(a\in PC_{\mathcal{L}(H)})\\
        \mathcal{E}^J_{\mathcal{S}, \mathcal{L}(H)}(J) &= \hat{J},\\
        \mathcal{E}^J_{\mathcal{F}, \mathcal{L}(H)}(P_nAP_n) &= \chi_{[-1,1]}\mathcal{E}_{\mathcal{S}, \mathcal{L}(H)}(A)\chi_{[-1,1]}. \qquad\qquad (A\in \mathcal{S}^J(PC_{\mathcal{L}(H)}))
        \end{align*}
\item Let $H = \mathbb{C}$. For $A\in \mathcal{S}^J(PC)$ and $(A_n)\in \mathcal{F}^J(PC)$, the limits
        \begin{align*}
        \mathcal{E}^J_\mathcal{S}(A) &:= \slim_{n\to \infty} E_n A E_{-n},\\
        \mathcal{E}^J_\mathcal{F}(A_n) &:= \slim_{n\to \infty} E_n A_n E_{-n}
        \end{align*}
exist. Moreover, the mappings $\mathcal{E}^J_\mathcal{S}$ and $\mathcal{E}^J_\mathcal{F}$ satisfy the conditions of (a).
\end{enumerate}
\end{prop}
\begin{proof} Without the flip operator $J$, it is simply just Proposition \ref{prop75} and was proved in \cite{ehrhardt1996finite}. In particular, in the scalar case, the mappings are exactly the ones introduced in Proposition \ref{prop75}. Furthermore, the proof used previously still applies to this case, as $E_{-n}E_n = I,\enskip E_{-n}^* = E_n$, and
                $$\mathcal{E}_\mathcal{S}(J) = \slim_{n\infty} E_nJE_{-n} = \hat{J}.$$
\end{proof}

Based on the proceeding proposition, the mappings
        \begin{align}
        \widetilde{\mathcal{E}}^J_{\mathcal{S}, \mathcal{L}(H)}&: \mathcal{S}_1^J(PC_{\mathcal{L}(H)}) \to \Sigma'_{\mathcal{L}(H)}, \quad A+\mathcal{I}^J_{1, \mathcal{L}(H)}\mapsto  \mathcal{E}^J_{\mathcal{S}, \mathcal{L}(H)}(A),\\
        \widetilde{\mathcal{E}}^J_{\mathcal{F}, \mathcal{L}(H)}&: \mathcal{F}_1^J(PC_{\mathcal{L}(H)}) \to \Xi'_{\mathcal{L}(H)}, \quad (A_n)+ \mathcal{J}^J_{1, \mathcal{L}(H)} \mapsto  \mathcal{E}^J_{\mathcal{F}, \mathcal{L}(H)}(A_n)
        \end{align}
are onto $^*$-homomorphisms.

In order to identify the local algebras at $\eta\in M_1^0(\widetilde{QC})$, it remains to show the kernels of $\Phi_{1, \mathcal{S}}, \enskip \Phi_{1, \mathcal{F}}, \enskip \widetilde{\mathcal{E}}^J_{\mathcal{S}, \mathcal{L}(H)}$ and $\widetilde{\mathcal{E}}^J_{\mathcal{F}, \mathcal{L}(H)}$ are trivial. For that purpose, we start investigating with the case when $H = \mathbb{C}$. Again, we will employ some of the results from \cite{ehrhardt1996finite} with necessary modifications.
Denote by $\omega$ the mapping $f\mapsto (f_1,f_2)^T$ with $f_1(x) = f(x)$ and $f_2(x) = f(-x),\enskip \forall x\in \mathbb{R}^+$. Obviously, $\omega$ is an isometry from $L^2(\mathbb{R})$ onto $L^2(\mathbb{R}^+)\oplus L^2(\mathbb{R}^+)$, as well as from $L^2([-1,1])$ onto $L^2([0,1])\oplus L^2([0,1])$. Therefore, the $^*$-isomorphism $\Phi_\omega$ given by
        \begin{equation}
        \Phi_\omega: \enskip A\mapsto \omega A \omega^{-1}
        \end{equation}
maps $\mathcal{L}(L^2(\mathbb{R}))$ onto $\mathcal{L}(L^2(\mathbb{R}^+))^{2\times 2}$ and $\mathcal{L}(L^2([-1,1]))$ onto $\mathcal{L}(L^2([0,1]))^{2\times 2}$. Furthermore, let $\widetilde{P}$ denote the operator $\chi_{[0,1]}$ of multiplication, let $S = S_{\mathbb{R}^+}$ stand for the singular integral operator on the positive semi-axis, and $N$ be the Hankel operator:
        \begin{align*}
        (Sf)(x) &= \frac{1}{\pi i}\int_0^\infty \frac{f(y)}{y-x}dy,\quad x\in \mathbb{R}^+,\\
        (Nf)(x) &= \frac{1}{\pi i}\int_0^\infty \frac{f(y)}{y+x}dy,\quad x\in \mathbb{R}^+.
        \end{align*}
Then we have (with $I$ refers to the identity operator on $L^2(\mathbb{R}^+)$)
        \begin{align*}
        \Phi_\omega(\chi_{[0,\infty)}) &=
        \left( \begin{array}{cc}
        I & 0\\
        0 & 0\\
        \end{array} \right),
        \quad \Phi_\omega(S_\mathbb{R}) =
        \left( \begin{array}{cc}
        S & -N\\
        N & -S\\
        \end{array} \right),\\
        \Phi_\omega(\chi_{[-1,1]}) &=
        \left( \begin{array}{cc}
        \widetilde{P} & 0\\
        0 & \widetilde{P}\\
        \end{array} \right),
        \quad \Phi_\omega(J) =
        \left( \begin{array}{cc}
        0 & I\\
        I & 0\\
        \end{array} \right).
        \end{align*}

Define the $C^*$-algebras $\Sigma_2'$ and $\Xi_2'$ by

        \begin{align}
        \Sigma_2' &:= \text{alg}_{\mathcal{L}(L^2(\mathbb{R}^+))^{2\times2}}\left\{
        \left( \begin{array}{cc}
        S & -N\\
        N & -S\\
        \end{array} \right),
        \left( \begin{array}{cc}
        I & 0\\
        0 & 0\\
        \end{array} \right),
        \left( \begin{array}{cc}
        0 & I\\
        I & 0\\
        \end{array} \right)
        \right\},\label{Sigma2'}\\
        \Xi_2' &:= \text{alg}_{\mathcal{L}(L^2([0,1]))^{2\times2}}\left\{
        \left( \begin{array}{cc}
        \widetilde{P} & 0\\
        0 & \widetilde{P}\\
        \end{array} \right)
        A
        \left( \begin{array}{cc}
        \widetilde{P} & 0\\
        0 & \widetilde{P}\\
        \end{array} \right):
        \enskip A\in \Sigma_2'
        \right\}.\label{Xi2'}
        \end{align}
A straightforward computation shows the following.
\begin{cor}\label{cor83} The $^*$-isomorphism $\Phi_\omega$ maps $\Sigma'$ onto $\Sigma_2'$ and $\Xi'$ onto $\Xi_2'$.
\end{cor}

Therefore, we will analyze $\Sigma_2'$ and $\Xi_2'$ instead of $\Sigma'$ and $\Xi'$. Recall that the Mellin transform $M$ and its inverse $M^{-1}$ are defined by
        \begin{align*}
        (Mf)(z) &= \int_0^\infty x^{-iz-\frac{1}{2}}f(x)dx,\quad z\in \mathbb{R}, \quad\qquad M:\enskip L^2(\mathbb{R}^+)\to L^2(\mathbb{R}),\\
        (M^{-1}f)(x) &= \frac{1}{2\pi}\int_{-\infty}^\infty x^{iz-\frac{1}{2}}f(z)dz,\quad x\in \mathbb{R}^+, \quad M^{-1}:\enskip L^2(\mathbb{R})\to L^2(\mathbb{R}^+).
        \end{align*}
For a multiplication operator $b\in L^\infty(\mathbb{R})$, let $M^0(b)$ denote the \emph{Mellin convolution operator}
        \begin{equation}
        M^0(b):=  M^{-1}bM, \enskip L^2(\mathbb{R}^+)\to L^2(\mathbb{R}^+).
        \end{equation}
Evidently, $\|M^0(b)\| = \|b\|,\enskip M^0(b)^* = M^0(b^*)$ and $M^0(b_1b_2) = M^0(b_1)M^0(b_2)$. For more properties,we refer to \cite{bottcher1990analysis}.

Let $PC_\infty(\mathbb{R})$ stand for the set of all continuous functions $f$ on $\mathbb{R}$ for which the limits $\lim\limits_{x\to \infty}f(x)$ and $\lim\limits_{x\to -\infty}f(x)$ exist and are finite, and let $C_\infty^0(\mathbb{R})$ be the set of all continuous functions $f$ on $\mathbb{R}$ for which $\lim\limits_{x\to \pm\infty}f(x) = 0$. Obviously, $PC_\infty(\mathbb{R})$ is a $C^*$-algebra with the $^*$-ideal $C_\infty^0(\mathbb{R})$.

As introduced in \cite{ehrhardt1996finite}, we define
                $$\Sigma_1 := \{\alpha I+ \beta S + M^0(b):\alpha,\beta\in \mathbb{C}, b\in C_\infty^0(\mathbb{R})\}.$$
A straightforward computation shows that $\Sigma_1$ is a $C^*$-algebra. By considering  the generating elements, 
we obtain the following.

\begin{prop}\label{prop84}
        \[
        \enskip \Sigma_2' := \left\{
        \left( \begin{array}{cc}
        A & B\\
        C & D\\
        \end{array} \right):
        \enskip A, B, C, D\in \Sigma_1
        \right\}.\qquad\qquad\qquad\]
\end{prop}

Similarly, define
                $$\Xi_1 :=\enskip \text{alg}_{\mathcal{L}(L^2([-1,1]))}\{\widetilde{P}A\widetilde{P}:\enskip A\in \Sigma_1\}.$$
As a consequence of Proposition \ref{prop84} and \eqref{Xi2'}, we immediately get the following result.
\begin{cor}\label{cor85}
\[
\enskip \Xi_2' := \left\{
\left( \begin{array}{cc}
A & B\\
C & D\\
\end{array} \right):
\enskip A, B, C, D\in \Xi_1
\right\}.\qquad\qquad\qquad\]
\end{cor}

Now we are ready to construct the inverses of $\widetilde{\mathcal{E}}_\mathcal{S}^J$ and $\widetilde{\mathcal{E}}_\mathcal{F}^J$. Note that, the mappings $E_n$ and $E_{-n}$ defined in \eqref{EN} and \eqref{E-N} can also be considered as acting on
        \[E_n:\enskip l^2(\mathbb{Z}^+)\to L^2(\mathbb{R}^+),\quad E_{-n}:\enskip L^2(\mathbb{R}^+)\to l^2(\mathbb{Z}^+).\]
For $b\in L^\infty(\mathbb{R})$, the operator $G(b)\in \mathcal{L}(l^2(\mathbb{Z}^+))$ given by
        \begin{equation}
        G(b):=\enskip E_{-1}M^0(b)E_1
        \end{equation}
is called \emph{discretized Mellin convolution operator}. The basic facts about $G(b)$ are introduced in \cite{hagen2012spectral}.

Let $\mathcal{F}_+$ stand for the set of all sequences $(A_n)\in \mathcal{F}$ with $(PA_nP) = (A_n)$. Regard to $\Sigma_1$ and $\Xi_1$, we define the bounded linear operators
        \begin{align*}
        \mathcal{E}_\mathcal{S}^0&:  \Sigma_1 \to \mathcal{L}(l^2(\mathbb{Z}^+)),\quad M^0(b)\mapsto G(b),\\
        \mathcal{E}_\mathcal{F}^0&:  \Xi_1 \to \mathcal{F}_+,\qquad\qquad A\mapsto (E_{-n}AE_n).
        \end{align*}

\begin{lem}\label{lem86}\hfill
\begin{enumerate}[label=(\alph*)]
\item $\mathcal{E}_\mathcal{S}^0(AB) - \mathcal{E}_\mathcal{S}^0(A)\mathcal{E}_\mathcal{S}^0(B)\in \mathcal{I}_1^J$ if $A\in \Sigma_1$ and $B\in \Sigma_1$.
\item $\mathcal{E}_\mathcal{F}^0(AB) - \mathcal{E}_\mathcal{F}^0(A)\mathcal{E}_\mathcal{F}^0(B)\in \mathcal{J}_1^J$ if $A\in \Xi_1$ and $B\in \Xi_1$.
\end{enumerate}
\end{lem}
\begin{proof} For $A\in \Sigma_1$ and $B\in \Sigma_1$, it is shown in Lemma 9.2 in \cite{ehrhardt1996finite} that $\mathcal{E}_\mathcal{S}^0(AB) - \mathcal{E}_\mathcal{S}^0(A)\mathcal{E}_\mathcal{S}^0(B)\in \mathcal{I}_1$, where
        \[\mathcal{I}_1 = \mathrm{clos}\enskip \mathrm{id}_{\mathcal{S}(PC)}\left(\mathcal{K}\cup \{f: f\in C(\mathbb{T}, f(1) = 0)\}\right).\]
Hence part (a) immediately follows since $\mathcal{I}_1 \subseteq \mathcal{I}_1^J$.

Part (b) can be treated analogously by using Lemma 9.3 in \cite{ehrhardt1996finite}.
\end{proof}
Similar to Lemma \ref{lem45}, we have
\begin{lem}\label{lem87} $(P_nBP_n)\in \mathcal{J}_1^J$ whenever $B\in \mathcal{I}_1^J$.
\end{lem}

Denote by $\mathcal{E}_\mathcal{S}'':\Sigma_2'\to \mathcal{L}(l^2(\mathbb{Z}^+))^{2\times 2}$ and $\mathcal{E}_\mathcal{F}'':\Xi_2'\to \mathcal{F}_+$ the mappings
        \begin{equation}
        \mathcal{E}_\mathcal{S}'':
        \left( \begin{array}{cc}
        A & B\\
        C & D\\
        \end{array} \right)\mapsto
        \left( \begin{array}{cc}
        \mathcal{E}_\mathcal{S}^0(A) & \mathcal{E}_\mathcal{S}^0(B)\\
        \mathcal{E}_\mathcal{S}^0(C) & \mathcal{E}_\mathcal{S}^0(D)\\
        \end{array} \right),\enskip
        \mathcal{E}_\mathcal{F}'':
        \left( \begin{array}{cc}
        A & B\\
        C & D\\
        \end{array} \right)\mapsto
        \left( \begin{array}{cc}
        \mathcal{E}_\mathcal{F}^0(A) & \mathcal{E}_\mathcal{F}^0(B)\\
        \mathcal{E}_\mathcal{F}^0(C) & \mathcal{E}_\mathcal{F}^0(D)\\
        \end{array} \right).
        \end{equation}
Further, we introduce the $^*$-homomorphisms $\Phi_{J,\mathcal{S}}$ and $\Phi_{J,\mathcal{F}}$ defined by
        \begin{align}
        \Phi_{J,\mathcal{S}}&: \mathcal{L}(l^2(\mathbb{Z}^+))^{2\times 2}\to \mathcal{L}(l^2(\mathbb{Z})),
        \left( \begin{array}{cc}
        A & B\\
        C & D\\
        \end{array} \right)\mapsto
        A+BJ+JC+JDJ,\qquad\\
        \Phi_{J,\mathcal{F}}&: \mathcal{F}_+^{2\times 2}\to \mathcal{F},
        \left( \begin{array}{cc}
        (A_n) & (B_n)\\
        (C_n) & (D_n)\\
        \end{array} \right)\mapsto
        (A_n+B_nJ+JC_n+JD_nJ).
        \end{align}
Recall that, the mappings we have constructed so far act between the following $C^*$-algebras:
        \begin{align}
        &\mathcal{S}^J(PC)\xrightarrow[]{\enskip\mathcal{E}_\mathcal{S}\enskip}\Sigma'\xrightarrow[]{\enskip\Phi_\omega\enskip}\Sigma_2'
        \xrightarrow[]{\enskip\mathcal{E}''_\mathcal{S}\enskip}\mathcal{L}(l^2(\mathbb{Z}^+))^{2\times2}\xrightarrow[]{\enskip\Phi_{J,\mathcal{S}}\enskip}
        \mathcal{L}(l^2(\mathbb{Z})),\\
        &\mathcal{F}^J(PC)\xrightarrow[]{\enskip\mathcal{E}_\mathcal{F}\enskip}\Xi'\xrightarrow[]{\enskip\Phi_\omega\enskip}\Xi_2'
        \xrightarrow[]{\enskip\mathcal{E}''_\mathcal{F}\enskip}\quad\enskip \mathcal{F}_+^{2\times 2}\enskip\enskip\quad\xrightarrow[]{\enskip\Phi_{J,\mathcal{F}}\enskip}
        \mathcal{F}.
        \end{align}
Finally, define the mappings $\mathcal{E}_\mathcal{S}'$ and $\mathcal{E}_\mathcal{F}'$ by
        \begin{equation}
        \mathcal{E}_\mathcal{S}':= \Phi_{J,\mathcal{S}}\circ \mathcal{E}_\mathcal{S}''\circ \Phi_\omega: \Sigma' \to \mathcal{L}(l^2(\mathbb{Z})),\enskip \mathcal{E}_\mathcal{F}':= \Phi_{J,\mathcal{F}}\circ \mathcal{E}_\mathcal{F}''\circ \Phi_\omega: \Xi' \to \mathcal{F}.
        \end{equation}
Simple computation shows that they are $^*$-homomorphisms, and
        \begin{align*}
        \mathcal{E}_\mathcal{S}'(A)&= E_{-n}AE_n,  \enskip \enskip\enskip\forall n\ge 1, A\in \Sigma',\\
        \mathcal{E}_\mathcal{F}'(A)&= (E_{-n}AE_n), \enskip \forall A\in \Xi'.
        \end{align*}

By Proposition \ref{prop75}, Corollary \ref{cor83} and Lemma \ref{lem86}, we have
\begin{lem}\label{lem88} Let $\rho_\mathcal{S} := \mathcal{E}_\mathcal{S}'\circ \mathcal{E}_\mathcal{S}:\enskip \mathcal{S}^J(PC)\to \mathcal{L}(l^2(\mathbb{Z}))$. Then\hfill
\begin{enumerate}[label=(\alph*)]
\item $\rho_\mathcal{S}(AB) - \rho_\mathcal{S}(A)\rho_\mathcal{S}(B)\in \mathcal{I}_1^J$ for all $A,B\in \mathcal{S}^J(PC)$.
\item $\rho_\mathcal{S}(A) - A \in \mathcal{I}_1^J$ for all $A\in \mathcal{S}^J(PC)$.
\end{enumerate}
\end{lem}
In addition, with Lemma \ref{lem87} and Lemma \ref{lem88}(a),
\begin{lem}\label{lem89} Let $\rho_\mathcal{F} := \mathcal{E}_\mathcal{F}'\circ \mathcal{E}_\mathcal{F}:\enskip \mathcal{F}^J(PC)\to \mathbb{F}$. Then\hfill
\begin{enumerate}[label=(\alph*)]
\item $\rho_\mathcal{F}(A_nB_n) - \rho_\mathcal{F}(A_n)\rho_\mathcal{F}(B_n)\in \mathcal{J}_1^J$ for all $(A_n),(B_n)\in \mathcal{F}^J(PC)$.
\item $\rho_\mathcal{F}(A_n) - (A_n) \in \mathcal{J}_1^J$ for all $(A_n)\in \mathcal{F}^J(PC)$.
\end{enumerate}
\end{lem}

In order to state the result, we use the notion of \emph{short exact sequences} and \emph{continuous cross-sections}. Recall that,  $0\xrightarrow{\enskip} \mathfrak{A}\xrightarrow{\alpha} \mathfrak{B}\xrightarrow{\beta} \mathfrak{C}\xrightarrow{\enskip} 0$ is a short exact sequence if $\alpha, \beta$ are $^*$-homomorphisms satisfying $\ker \alpha = \{0\}$, $\ker\beta = \mathrm{Im}\alpha$, and $ \mathfrak{C} = \mathrm{Im}\beta$. Further, $\rho$ is a continuous cross-section of $\beta$ if $\rho$ is a linear and continuous mapping $\mathfrak{C}\to \mathfrak{B}$ with $\beta \circ \rho = \mathrm{id}$. 

\begin{prop}\label{prop810}\hfill
\begin{enumerate}[label=(\alph*)]
\item The sequence $0\xrightarrow{\quad} \mathcal{I}_1^J \xrightarrow{\enskip \mathrm{id}\enskip}  \mathcal{S}^J(PC)\xrightarrow{\enskip\mathcal{E}_\mathcal{S}\enskip} \Sigma'\xrightarrow{\quad} 0$ is short exact, and $\mathcal{E}_\mathcal{S}'$ is a continuous cross-section of $\mathcal{E}_\mathcal{S}$.
\item The sequence $0\xrightarrow{\quad} \mathcal{J}_1^J \xrightarrow{\enskip \mathrm{id}\enskip} \mathcal{F}^J(PC)\xrightarrow{\enskip\mathcal{E}_\mathcal{F}\enskip} \Xi'\xrightarrow{\quad} 0$ is short exact, and $\mathcal{E}_\mathcal{F}'$ is a continuous cross-section of $\mathcal{E}_\mathcal{F}$.
\end{enumerate}
\end{prop}

The previous proposition implies that $\mathcal{F}_1^J(PC)$ is $^*$-isomorphic to $\Xi'$, which means in the scalar-valued case, the $^*$-homomorphism $\widetilde{\mathcal{E}}_\mathcal{F}$ is actually a $^*$-isomorphism. In general, in the $\mathcal{L}(H)$-valued case, we need the following result on tensor products from \cite{douglas1972banach}. A proof can be found in \cite{douglas1971C}.

\begin{prop}\label{prop811}If $0\xrightarrow{\enskip} \mathfrak{A}\xrightarrow{\alpha} \mathfrak{B}\xrightarrow{\beta} \mathfrak{C}\xrightarrow{\enskip} 0$ is a short exact sequence of $C^*$-algebras such that $\beta$ has a continuous cross-section $\rho$ and $\mathfrak{D}$ is a $C^*$-algebra, then the sequence
        \[0\xrightarrow{\quad} \mathfrak{A}\otimes\mathfrak{D} \xrightarrow{\alpha\otimes\text{id}} \mathfrak{B}\otimes\mathfrak{D} \xrightarrow{\beta\otimes\text{id}} \mathfrak{C}\otimes\mathfrak{D}\xrightarrow{\quad} 0\]
is short exact.
\end{prop}

\begin{cor}\label{cor812} The sequence
        \[0\xrightarrow{\quad} \mathcal{J}_{1,\mathcal{L}(H)}^J\xrightarrow{\enskip \mathrm{id}\enskip} \mathcal{F}^J(PC_{\mathcal{L}(H)})\xrightarrow{ \mathcal{E}_{\mathcal{F},\mathcal{L}(H)}} \Xi_{\mathcal{L}(H)}'\xrightarrow{\quad} 0\]
is short exact. Furthermore, the mapping $\mathcal{E}_{\mathcal{F},\mathcal{L}(H)}':\Xi_{\mathcal{L}(H)}'\to \mathcal{F}^J(PC_{\mathcal{L}(H)})$ defined by $\mathcal{E}_{\mathcal{F},\mathcal{L}(H)}'(A) :=(E_{-n}AE_n)$ is a continuous cross-section of $\mathcal{E}_{\mathcal{F},\mathcal{L}(H)}$.
\end{cor}
\begin{proof}
Note that $\mathcal{F}^J(PC_{\mathcal{L}(H)}) \cong \mathcal{F}^J(PC) \otimes \mathcal{L}(H)$ and $\mathcal{J}^J_{1,\mathcal{L}(H)} \cong \mathcal{J}_1^J\otimes \mathcal{L}(H)$. Therefore, with $\mathcal{E}_{\mathcal{F},\mathcal{L}(H)}' \cong \mathcal{E}_{\mathcal{F}}'\otimes id$, the mapping $\mathcal{E}_{\mathcal{F},\mathcal{L}(H)}'$ maps into $\mathcal{F}^J(PC_{\mathcal{L}(H)})$ and is a continuous cross-section of $\mathcal{E}_{\mathcal{F},\mathcal{L}(H)}$. By Proposition \ref{prop810}(b) and Proposition \ref{prop811}, the sequence
                $$0\xrightarrow{\quad} \mathcal{J}_1^J\otimes\mathcal{L}(H) \xrightarrow{\mathrm{id}\otimes\mathrm{id}} \mathcal{F}^J(PC)\otimes\mathcal{L}(H) \xrightarrow{\mathcal{E}_\mathcal{F}\otimes\mathrm{id}} \Xi'\otimes\mathcal{L}(H)\xrightarrow{\quad} 0$$
is short exact, which proves the assertion. For more details, see Proposition 7.2(a) and also Corollary 9.10 from \cite{ehrhardt1996finite}.
\end{proof}
\begin{cor}\label{cor813}$\widetilde{\mathcal{E}}_{\mathcal{F},\mathcal{L}(H)}$ is a $^*$-isomorphism from $\mathcal{F}_1^J(PC_{\mathcal{L}(H)})$ onto $\Xi'_{\mathcal{L}(H)}$, and
                $$\widetilde{\mathcal{E}}_{\mathcal{F},\mathcal{L}(H)}^{-1}:\enskip A\mapsto \mathcal{E}_{\mathcal{F},\mathcal{L}(H)}'(A)+ \mathcal{J}_{1,\mathcal{L}(H)}^J.$$
\end{cor}

Now, by applying the same argument used in the proof of Proposition \ref{prop77}, we can show that
\begin{prop}\label{prop814} The kernels of $\Phi_{1,\mathcal{S}}$ and $\Phi_{1,\mathcal{F}}$ defined in Proposition \ref{prop81} are trivial.
\end{prop}
\begin{proof}
Here we provide a short proof. Assume $\ker \Phi_{1, \mathcal{F}}\neq\{0\}$. Since $\mathcal{F}_1^J(PC_{\mathcal{L}}(H))$ and $\Xi'_{\mathcal{L}(H)}$ are $^*$-isomorphic via $\tilde{\mathcal{E}}_{\mathcal{F}, \mathcal{L}(H)}$, the ideal $\ker \Phi_{1, \mathcal{F}}$ corresponds to a non-trivial ideal $\mathfrak{J}$ of $\Xi'_{\mathcal{L}(H)}$. Again, we define the operators $K_{x,y}$ as in
\eqref{Kxy}, and $K_{x,y}\in \Xi_{\mathcal{L}(H)} \subseteq \Xi'_{\mathcal{L}(H)}$. Similarly, with $x = \chi_{[-1,1]}\in L^2([-1,1])$, there exists an $h\in H, \enskip \|h\| = 1$ such that $K_{xh,xh}\in \mathfrak{J}$ is non-trivial. Therefore, we obtain
        \[\tilde{\mathcal{E}}^{-1}_{\mathcal{F}, \mathcal{L}(H)}(K_{xh,xh}) = (E_{-n}K_{xh,xh}E_n) + \mathcal{J}^J_{1,\mathcal{L}(H)}\in \ker \Phi_{1, \mathcal{F}},\]
and
        \begin{equation}
        (K_n) := (E_{-n}K_{xh,xh}E_n)\in \mathcal{J}_{\eta, \mathcal{L}(H)}^J.
        \end{equation}
Given $\epsilon >0$, by the definition of $\mathcal{J}_{\eta, \mathcal{L}(H)}^J$, there is a sequence $(A_n)$ for which
        \begin{gather*}
        \|(A_n) - (K_n)\|_\mathcal{F}\le\epsilon,\\
        (A_n) = \sum_{i=1}^k (A_n^{(i)})(P_nf_iP_n) + (B_n'),
        \end{gather*}
where $(A_n^{(i)})\in \mathcal{F}^J({PQC_{\mathcal{L}(H)}}), \enskip f_i\in \widetilde{QC},\enskip \eta(f_i) = 0$ and $(B_n')\in \mathcal{J}$. Similarly, there is an $f\in \widetilde{QC}$ with $\eta(f)= 1, \|f\|= 1$ and $\sum\limits_{i=1}^k \|A_n^{(i)}\|\cdot\|ff_i\|\le\epsilon$. Therefore,
        \[\|(A_n)(P_nfP_n) - (B_n)\|\le \epsilon,\]
where $(B_n)=(P_nKP_n+W_nLW_n + C_n')\in \mathcal{J}$ for some $K, L\in \mathcal{K}$ and $(C_n')\in \mathcal{N}$. Now, choose a $\kappa$ sufficiently large such that $\|Q_{\kappa}K\|_\mathcal{A}\le \epsilon$ and $\|Q_{\kappa}L\|_\mathcal{A}\le \epsilon$, and set
        \begin{equation}
        (R_n) := (P_nQ_\kappa P_n)(W_n Q_\kappa W_n) = (W_n Q_\kappa W_n)(P_nQ_\kappa P_n).
        \end{equation}
A simple computation shows that there exists some $(C_n)\in \mathcal{N}$ such that
        \[\|(R_n)(B_n) - (C_n)\|_\mathcal{F} \le 2\epsilon.\]
Observing that $\|(R_n)\| = 1$, $\|f\| = 1$, we obtain
        \begin{equation}
        \limsup_{n\to \infty}\|R_nK_nP_nfP_n\|_{\mathcal{L}(l_H^2(\mathbb{Z}_n))}\le 4\epsilon.\label{contra2}
        \end{equation}
On the other hand, let $z_n$ and $z_n^*$ be the bounded linear operators introduced in \eqref{zzz} and \eqref{zzzz}. Then $\|z_n\| = \|z_n^*\| = 1$, and it yields that
        \begin{equation}
        (E_{-n}K_{xh,xh}E_n) = (z_nz_n^*).
        \end{equation}
Therefore, for $n\ge 2\kappa$,
        \begin{align*}
        \|R_nK_nP_nfP_n\| &\ge |z_n^*R_nK_nP_nfP_nz_n|\\
        &\ge |z_n^*R_nz_n|\cdot|z_n^*P_nfP_nz_n|\\
        &= \frac{2n-4\kappa}{2n} \cdot |\sigma_{2n-1}f(1)|.
        \end{align*}
For $\eta\in M^0_1(\widetilde{QC})$, there exists a unique $\xi\in M_1^0(QC)$ such that $\hat{\xi} = \eta$. Since $\xi(f) = 1$, we obtain from Lemma \ref{lem76} that
        \begin{equation*}
        \limsup_{n\to\infty}\|R_nK_nP_nfP_n\| \ge \limsup_{n\to\infty} |\sigma_{2n-1}f(1)|\ge |\xi(f)| = 1,
        \end{equation*}
which contradicts \eqref{contra2} when $\epsilon$ is chosen sufficiently small. Hence $\ker \Phi_{\tau, \mathcal{F}}$ is trivial. Analogously, one can show that $\ker \Phi_{\tau, \mathcal{S}}$ is also trivial by following the proof of Proposition \ref{prop77}.
\end{proof}
To summarize what we have so far, the next result directly follows from Proposition  \ref{prop81}, Proposition \ref{prop82}, Corollary \ref{cor813} and  Proposition \ref{prop814}
\begin{cor}\label{cor915}
For $\eta\in M_1^0(\widetilde{QC})$, the map
	\begin{equation}
		\Phi^{-1}_{1, \mathcal{F}} \circ \tilde{\mathcal{E}}_{\mathcal{F}, \mathcal{L}(H)}: \mathcal{F}^J_\eta(PQC_{\mathcal{L}(H)}) \to \Xi_{\mathcal{L}(H)}' 
	\end{equation}
is a $^*$-isomorphism. 
\end{cor}

\begin{rmk}\label{rmk916} Let us consider the case where $\eta \in M_{-1}^0(\widetilde{QC})$. Similarly, \eqref{decomp 91} and \eqref{decomp 92} still hold, and one can define the mappings
        \begin{align}
        & \Phi_{-1, \mathcal{S}}: \mathcal{S}_1^J(PC_{\mathcal{L}(H)}) \to \mathcal{S}_\eta^J(PQC_{\mathcal{L}(H)}), A+ \mathcal{I}^J_{1, \mathcal{L}(H)} \mapsto Y_{-1}AY_{-1}^{-1} + \mathcal{I}^J_{\eta,\mathcal{L}(H)},\\
        & \Phi_{-1, \mathcal{F}}: \mathcal{F}_1^J(PC_{\mathcal{L}(H)}) \to \mathcal{F}_\eta^J(PQC_{\mathcal{L}(H)}), (A_n)+ \mathcal{J}^J_{1, \mathcal{L}(H)} \mapsto (Y_{-1}A_nY_{-1}^{-1}) + \mathcal{J}^J_{\eta,\mathcal{L}(H)},
        \end{align}
which are $^*$-isomorphisms as well. By Proposition \ref{prop814}, the mapping
	\begin{equation}
		\Phi^{-1}_{-1, \mathcal{F}} \circ \tilde{\mathcal{E}}_{\mathcal{F}, \mathcal{L}(H)}: \mathcal{F}^J_\eta(PQC_{\mathcal{L}(H)}) \to \Xi_{\mathcal{L}(H)}' 
	\end{equation}
is a $^*$-isomorphism as desired.
\end{rmk}
\section{Summary and main results}\label{s:10}
In summary, for $\eta\in M_\tau^0(\widetilde{QC})$ where $\tau\neq\pm1$, we have the following diagram constructed:
\begin{center}
\begin{tikzpicture}[thick,scale=0.5]
 \matrix (m) [matrix of math nodes,row sep=2em,column sep=3em,minimum width=2em]
  {
     \mathcal{S}^J(PQC_{\mathcal{L}(H)}) & \mathcal{S}_\eta^J(PQC_{\mathcal{L}(H)})&(p\mathcal{X}_\mathcal{S}p)^{2\times2}&\mathcal{S}_1(PC_{\mathcal{L}(H)})^{2\times2}& \Sigma_{\mathcal{L}(H)}^{2\times 2}\\
     \mathcal{F}^J(PQC_{\mathcal{L}(H)}) & \mathcal{F}_\eta^J(PQC_{\mathcal{L}(H)})&(p'\mathcal{X}_\mathcal{F}p')^{2\times2}&\mathcal{F}_1(PC_{\mathcal{L}(H)})^{2\times2}& \Xi_{\mathcal{L}(H)}^{2\times 2} \\};
  \path[-stealth]
        (m-1-1) edge node [above] {\footnotesize{$\Phi_{\eta,\mathcal{S}}$}} (m-1-2)
        (m-1-2) edge node [above] {\footnotesize{$L_{\eta,\mathcal{S}}$}} (m-1-3)
        (m-1-4) edge node [above] {\footnotesize{$\Phi'_{\tau,\mathcal{S}}$}} (m-1-3)
        (m-1-4) edge node [above] {\footnotesize{$\widetilde{\mathcal{E}}'_{\mathcal{S,\mathcal{L}(H)}}$}} (m-1-5)
        (m-2-1) edge node [left]  {} (m-1-1)
                edge node [above] {\footnotesize{$\Phi_{\eta,\mathcal{F}}$}} (m-2-2)
        (m-2-2) edge node [left]  {} (m-1-2)
                edge node [above] {\footnotesize{$L_{\eta,\mathcal{F}}$}} (m-2-3)
        (m-2-3) edge node [left]  {} (m-1-3)
        (m-2-4) edge node [left]  {} (m-1-4)
                edge node [above] {\footnotesize{$\Phi'_{\tau,\mathcal{F}}$}} (m-2-3)
                edge node [above] {\footnotesize{$\widetilde{\mathcal{E}}'_{\mathcal{F,\mathcal{L}(H)}}$}} (m-2-5)
        (m-2-5) edge node [left]  {} (m-1-5);
\end{tikzpicture}
\end{center}
By Proposition \ref{prop75}, the mappings $\widetilde{\mathcal{E}}'_{\mathcal{S,\mathcal{L}(H)}}$ and $\widetilde{\mathcal{E}}'_{\mathcal{F,\mathcal{L}(H)}}$ are $^*$-isomorphisms. Further, we see that $\Phi_{\tau,\mathcal{S}}'$ and $\Phi_{\eta,\mathcal{F}}'$ are $^*$-isomorphisms as well by Proposition \ref{prop74} and Proposition \ref{prop77}.

When $\eta\in M_{\pm1}^0(\widetilde{QC})$, the $^*$-homomorphisms act between $C^*$-algebras in the following diagram:
\begin{center}
\begin{tikzpicture}[thick,scale=0.5]
 \matrix (m) [matrix of math nodes,row sep=2em,column sep=3em,minimum width=2em]
  {
     \mathcal{S}^J(PQC_{\mathcal{L}(H)}) & \mathcal{S}_\eta^J(PQC_{\mathcal{L}(H)})&\mathcal{S}^J_1(PC_{\mathcal{L}(H)})& \Sigma'_{\mathcal{L}(H)}\\
     \mathcal{F}^J(PQC_{\mathcal{L}(H)}) & \mathcal{F}_\eta^J(PQC_{\mathcal{L}(H)})&\mathcal{F}^J_1(PC_{\mathcal{L}(H)})& \Xi'_{\mathcal{L}(H)}\\};
  \path[-stealth]
        (m-1-1) edge node [above] {\footnotesize{$\Phi_{\eta,\mathcal{S}}$}} (m-1-2)
        (m-1-3) edge node [above] {\footnotesize{$\Phi_{\pm 1,\mathcal{S}}$}} (m-1-2)
        (m-1-3) edge node [above] {\footnotesize{$\widetilde{\mathcal{E}}_{\mathcal{S},\mathcal{L}(H)}$}} (m-1-4)
        (m-2-1) edge node [left]  {} (m-1-1)
                edge node [above] {\footnotesize{$\Phi_{\eta,\mathcal{F}}$}} (m-2-2)
        (m-2-2) edge node [left]  {} (m-1-2)
        (m-2-3) edge node [above] {\footnotesize{$\Phi_{\pm 1,\mathcal{F}}$}} (m-2-2)
        (m-2-3) edge node [left]  {} (m-1-3)
                edge node [above] {\footnotesize{$\widetilde{\mathcal{E}}_{\mathcal{F},\mathcal{L}(H)}$}} (m-2-4)
        (m-2-4) edge node [left]  {} (m-1-4);
\end{tikzpicture}
\end{center}
Combining Proposition \ref{prop82} and Corollary \ref{cor813}, we see $\widetilde{\mathcal{E}}_{\mathcal{F,\mathcal{L}(H)}}$ (and $\widetilde{\mathcal{E}}_{\mathcal{S,\mathcal{L}(H)}}$) are $^*$-isomorphisms. And it follows from Proposition \ref{prop81}, Proposition \ref{prop814} and Remark \ref{rmk916} that $\Phi_{\pm 1,\mathcal{S}}$ as well as $\Phi_{\pm 1,\mathcal{F}}$ are also $^*$-isomorphisms. 

Now, we are able to define $^*$-homomorphisms from $\mathcal{S}^J(PQC_{\mathcal{L}(H)})$ and $\mathcal{F}^J(PQC_{\mathcal{L}(H)})$ onto the desired $C^*$-algebras respectively. Recall that, for each $\eta\in M_\tau^0(\widetilde{QC})$, there is a unique $\xi\in M_\tau^0(QC)$ such that $\hat{\xi} = \eta$. For $q\in QC_{\mathcal{L}(H)}^s$, since $QC_{\mathcal{L}(H)}^s$ is locally trivial at $\xi$, there is a uniquely determined $a\in \mathcal{L}(H)$, such that $q-a\in I_{\xi,\mathcal{L}(H)}$ and $\Phi_\xi(q)$ = a (see \eqref{locallytrivial} and \cite{ehrhardt1996symbol}). As elements in $\mathcal{L}(H)$ can be regarded as constant functions in $L^\infty_{\mathcal{L}(H)}$, we have
\begin{thm}\label{thm91}Let $\eta\in M_\tau^0(\widetilde{QC})$. Then there is a unique $\xi\in M_\tau^0(QC)$ such that $\hat{\xi} = \eta$.\hfill
\begin{enumerate}[label=(\alph*)]
\item If $\tau\neq \pm1$, then the mappings
        \begin{align}
        \Psi_{\eta,\mathcal{S}}&: \mathcal{S}^J(PQC_{\mathcal{L}(H)})\to \Sigma_{\mathcal{L}(H)}^{2\times2}, A\mapsto
        \Psi'_{\eta,\mathcal{S}}(A+\mathcal{I}^J_{\eta,\mathcal{L}(H)}),\\
        \Psi_{\eta,\mathcal{F}}&: \mathcal{F}^J(PQC_{\mathcal{L}(H)})\to \Xi_{\mathcal{L}(H)}^{2\times2}, (A_n)\mapsto
        \Psi'_{\eta,\mathcal{F}}((A_n)+\mathcal{J}^J_{\eta,\mathcal{L}(H)})
        \end{align}
are $^*$-homomorphisms, where $\Psi'_{\eta,\mathcal{S}} = \widetilde{\mathcal{E}}'_{\mathcal{S,\mathcal{L}(H)}}\circ \Phi_{\tau,\mathcal{S}}^{-1}\circ L_{\eta,\mathcal{S}}$ and $\Psi'_{\eta,\mathcal{F}} = \widetilde{\mathcal{E}}'_{\mathcal{F,\mathcal{L}(H)}}\circ \Phi_{\tau,\mathcal{F}}^{-1}\circ L_{\eta,\mathcal{F}}$ are defined in \eqref{Psi notpm1}. Furthermore, $\mathcal{J}\subseteq \ker \Psi_{\eta,\mathcal{F}}$, and
        \begin{align}
        \Psi_{\eta,\mathcal{S}}(P) &=
        \left( \begin{array}{cc}
        \chi_{[0,\infty)} & 0\\
        0 & \chi_{(-\infty,0]}\\
        \end{array} \right),\\
        \Psi_{\eta,\mathcal{S}}(a) &=
        \left( \begin{array}{cc}
        a(\tau+0)\frac{I-S_\mathbb{R}}{2}+ a(\tau-0)\frac{I+S_\mathbb{R}}{2} & 0\\
        0 & \tilde{a}(\tau+0)\frac{I-S_\mathbb{R}}{2}+ \tilde{a}(\tau-0)\frac{I+S_\mathbb{R}}{2}\\
        \end{array} \right),\\
        \Psi_{\eta,\mathcal{S}}(q) &=
        \left( \begin{array}{cc}
        \Phi_\xi(q) & 0\\
        0 & \Phi_\xi(\tilde{q})\\
        \end{array} \right),\\
        \Psi_{\eta,\mathcal{S}}(J) &=
        \left( \begin{array}{cc}
        0 & I\\
        I & 0\\
        \end{array} \right),\\
        \Psi_{\eta,\mathcal{F}}(P_nAP_n) &= \chi_{[-1,1]}\Psi_{\eta,\mathcal{S}}(A)\chi_{[-1,1]},
        \end{align}
where $a\in PC_{\mathcal{L}(H)}, q\in QC_{\mathcal{L}(H)}^s$ and $A\in \mathcal{S}^J(PQC_{\mathcal{L}(H)})$.
\item If $\tau = \pm1$, the mappings
        \begin{align}
        \Psi^1_{\eta,\mathcal{S}}&: \mathcal{S}^J(PQC_{\mathcal{L}(H)})\to \Sigma'_{\mathcal{L}(H)}, A\mapsto
        (\widetilde{\mathcal{E}}_{\mathcal{S,\mathcal{L}(H)}}\circ \Phi_{\tau,\mathcal{S}}^{-1})(A+\mathcal{I}^J_{\eta,\mathcal{L}(H)}),\\
        \Psi^1_{\eta,\mathcal{F}}&: \mathcal{F}^J(PQC_{\mathcal{L}(H)})\to \Xi'_{\mathcal{L}(H)}, (A_n)\mapsto
        (\widetilde{\mathcal{E}}_{\mathcal{F,\mathcal{L}(H)}}\circ \Phi_{\tau,\mathcal{F}}^{-1})((A_n)+\mathcal{J}^J_{\eta,\mathcal{L}(H)})\qquad
        \end{align}
are $^*$-homomorphisms. Moreover, $\mathcal{J}\subseteq \ker \Psi^1_{\eta,\mathcal{F}}$, and
        \begin{align}
        \Psi^1_{\eta,\mathcal{S}}(P) &= \chi_{[0,\infty)},\\
        \Psi^1_{\eta,\mathcal{S}}(a) &= a(\tau+0)\frac{I-S_\mathbb{R}}{2}+ a(\tau-0)\frac{I+S_\mathbb{R}}{2}, \enskip(a\in PC_{\mathcal{L}(H)})\\
        \Psi^1_{\eta,\mathcal{S}}(q) &= \Phi_\xi(q), \enskip\quad\qquad\qquad\qquad\qquad\quad\qquad (q\in QC_{\mathcal{L}(H)}^s)\\
        \Psi^1_{\eta,\mathcal{S}}(J) &= \tau\hat{J},\\
        \Psi^1_{\eta,\mathcal{F}}(P_nAP_n) &= \chi_{[-1,1]}\Psi^1_{\eta,\mathcal{S}}(A)\chi_{[-1,1]}. \qquad\qquad (A\in \mathcal{S}^J(PQC_{\mathcal{L}(H)}))
        \end{align}
\end{enumerate}
\end{thm}

With the help from the $^*$-homomorphisms $\mathcal{P}$ and $\mathcal{W}$ defined via (\ref{operatorP}) and (\ref{operatorW}), we are finally able to state the main result.
\begin{thm}\label{thm92} Let $(A_n)\in \mathcal{F}^J(PQC_{\mathcal{L}(H)})$. Then the sequence $(A_n)$ is stable if and only if the following statements all hold:\hfill
\begin{enumerate}[label=(\alph*)]
\item $\mathcal{P}(A_n)$ is invertible in $\mathcal{S}^J(PQC_{\mathcal{L}(H)})$,
\item $\mathcal{W}(A_n)$ is invertible in $\mathcal{S}^J(PQC_{\mathcal{L}(H)})$,
\item $\Psi_{\eta,\mathcal{F}}(A_n)$ is invertible in $\Xi_{\mathcal{L}(H)}^{2\times 2}$ for all $\eta\in M_\tau^0(\widetilde{QC})$ with $\tau\in \mathbb{T}\setminus\{\pm1\}$, and
\item $\Psi^1_{\eta,\mathcal{F}}(A_n)$ is invertible in $\Xi'_{\mathcal{L}(H)}$ for all $\eta\in M_\tau^0(\widetilde{QC})$ with $\tau = \pm1$.
\end{enumerate}
\end{thm}
\begin{proof}
Using  Theorem \ref{thm37} and Corollary \ref{cor44} it follows that the sequence $(A_n)$ is stable if and only if 
$\mathcal{P}(A_n)$ and $\mathcal{W}(A_n)$ are invertible, and if $(A_n)+\mathcal{J}^J_{\eta,\mathcal{L}(H)}$ is invertible in $\mathcal{F}_\eta^J(PQC_{\mathcal{L}(H)})$  for each $\eta\in M(\widetilde{QC})$. 
Based on the decomposition of $M(\widetilde{QC})$ obtained from Section 6, it suffices to consider the invertiblity in the following three cases, where the following holds:

 It follows directly from Theorem \ref{thm37}, Corollary \ref{cor44}, Theorem \ref{thm61}, Theorem \ref{thm62} and the fact that the mappings $\Psi'_{\eta,\mathcal{F}}$ (for $\eta\in M_\tau^0(\widetilde{QC})$ where $\tau\neq\pm1$), $\widetilde{\mathcal{E}}_{\mathcal{F,\mathcal{L}(H)}}$ (maps $\mathcal{F}_{\pm1}^J(PC_{\mathcal{L}(H)})$ onto $\Xi'_{\mathcal{L}(H)}$) and $\Phi_{\tau,\mathcal{F}}^{-1}$ ($\tau = \pm1$) are $^*$-isomorphisms. The latter means that $(A_n)+\mathcal{J}^J_{\eta,\mathcal{L}(H)}$ is invertible in $\mathcal{F}_\eta^J(PQC_{\mathcal{L}(H)})$ if and only if $\Psi_{\eta,\mathcal{F}}(A_n)$ or $\Psi^1_{\eta,\mathcal{F}}(A_n)$ is invertible in the corresponding $C^*$-algebras.
\end{proof}

\red{another short version of the proof}
\begin{proof}
	By Theorem \ref{thm41}, the sequence $(A_n)$ is stable if and only if $(A_n)+\mathcal{J}^J_{\eta,\mathcal{L}(H)}$ is invertible in $\mathcal{F}_\eta^J(PQC_{\mathcal{L}(H)})$  for each $\eta\in M(\widetilde{QC})$. Based on the decomposition of $M(\widetilde{QC})$ obtained from Section 6, it suffices to consider the invertiblity in the following three cases:
\begin{enumerate}
\item[(i)]
For $\eta \in M_\tau(\widetilde{QC})\setminus M^0_\tau(\widetilde{QC})$ with $\tau\in \overline{\mathbb{T}_+}$,  $(A_n)+\mathcal{J}^J_{\eta,\mathcal{L}(H)}$ is invertible if $\mathcal{P}(A_n)$ is invertible by Theorem \ref{thm61} and Theorem \ref{thm62}.
\item[(ii)]
 For $\eta\in  M^0_\tau(\widetilde{QC})$ with $\tau\in \mathbb{T}_+$, $(A_n)+\mathcal{J}^J_{\eta,\mathcal{L}(H)}$ is invertible if and only if $\Psi_{\eta,\mathcal{F}}(A_n)$ is invertible in $\Xi_{\mathcal{L}(H)}^{2\times 2}$ by Corollary \ref{cor78} and Theorem \ref{thm91}(a).
\item[(iii)]
 For $\eta\in  M^0_\tau(\widetilde{QC})$ with $\tau=\pm1$, $(A_n)+\mathcal{J}^J_{\eta,\mathcal{L}(H)}$ is invertible if and only if $\Psi_{\eta,\mathcal{F}}^1(A_n)$ is invertible in $\Xi_{\mathcal{L}(H)}'$ 
 by Corollary \ref{cor915}, Remark \ref{rmk916} and Theorem \ref{thm91}(b).
 \end{enumerate}
This implies the stability criterion.
\end{proof}

\nocite{*}


\end{document}